
\input epsf

\outer \def \Picture #1. #2 #3\par {{\parindent=0pt
\parskip=0pt 
\vfil
\centerline {\epsfbox {#2}}
\bigskip
{\narrower {\bold Figure #1.} \sl #3\par}
\vfil}}

 \outer\def\Figure #1. #2\par {\vfil \eject \vglue7.5in {\narrower
\noindent {\bold Figure #1. }{\sl#2}\par} \vfil \eject}
 \abovedisplayskip = 12pt plus 3pt minus 3pt
 \abovedisplayshortskip = 6pt plus 3pt minus 3pt
 \belowdisplayskip = 12pt plus 3pt minus 3pt
 \belowdisplayshortskip = 12pt plus 3pt minus 3pt
 \hsize = 5.0in
 \hoffset = 0.8125in
 \parskip = 6pt plus 1pt minus 1pt
 \topskip=18pt

  \font\bigbold=cmbx18
 \font\medbold=cmbx14
 \font\smallbold=cmbx12
 
 \font\bold=cmbx10
 \font\helvetica=cmss10


 \headline={\hfil}
 
 \def\thismonth{\ifcase \month \or JANUARY\or FEBRUARY\or
MARCH\or APRIL\or MAY\or JUNE\or JULY\or AUGUST\or
SEPTEMBER\or OCTOBER\or NOVEMBER\or DECEMBER\fi}

 \def\footremark{{\ }}
 \def\pagenumber{{\tenrm \the \pageno}}


 \footline=
{\ifnum \pageno=1{\hfil}\else 
{\ifnum \pageno=2{\hfil}\else 
{\ifnum \pageno=5{\noindent \hfil \footremark \hfil {\tenrm 5/6}}\else
{\ifnum \pageno=6{\hfil}\else 
{\ifodd \pageno {\noindent \hfil \footremark \hfil \pagenumber}\else
{\noindent \pagenumber \hfil \footremark \hfil}\fi}\fi}\fi}\fi}\fi}

 \outer\def\beginsection#1\par{\vskip0pt plus2.0in\penalty
-250\vskip0pt plus-2.0in\bigskip \bigskip \message {#1}\leftline
{\medbold #1}\nobreak \smallskip \noindent}

 \outer\def\beginsubsection#1\par{\vskip0pt plus2.0in\penalty
-250\vskip0pt plus-2.0in\bigskip \message{#1}\leftline {\smallbold
#1}\nobreak \smallskip \noindent}

 \outer\def\Proclaim #1. #2\par {\medbreak {\narrower \noindent
{\bold #1. }{\sl #2}\par} \ifdim \lastskip <\medskipamount
\removelastskip \penalty55 \medskip \fi}

 \outer\def\table #1. #2 #3\par 
 {{\vglue 0.1in\narrower \noindent
 {\bold Table #1. }{\sl#2}\par
 \bigskip
 {\noindent \hfil{#3}}
 }}

 \outer\def\Table #1. #2 #3\par 
 {{\vfil \eject \ \vfil \narrower \noindent
 {\bold Table #1. }{\sl#2}\par
 \bigskip
 {\hfil#3\hfil} 
 \vglue1.0in 
 \vfil \eject}}

 \def \Ddots{\mathinner {\mkern 1mu \raise5pt \vbox {\kern3pt
\hbox{.}} \mkern 2mu \raise 3pt \hbox{.} \mkern 2mu \raise1pt
\hbox{.} \mkern1mu}}

 \def \halmos {{\vrule height5pt width3pt depth2pt}}

 \def \ocm (#1 #2){oc$\,(#1,#2)$}
 \def \Ocm (#1 #2){Oc$\,(#1,#2)$}

 \def\picture #1 by #2 (#3){$$\vbox to #2{\hrule width #1 height
0pt depth 0pt\vfill \special{picture #3}}$$}

 \def \Span #1{\hbox {span} \, \left\{ #1 \right\}}

 \def \tableau#1{\null \,\vcenter {\normalbaselines \ialign {\vrule
depth1.5ex height2ex width0em \hfil $##$&& \hfil $\quad
##$\crcr \mathstrut \crcr \noalign {\kern-\baselineskip} #1\crcr
\mathstrut \crcr \noalign {\kern-\baselineskip}}}\,}

 \def \transpose #1{{#1}^*}

 \def\Vdots{\lower2pt \vbox {\baselineskip=4pt
\lineskiplimit=0pt \kern2pt \hbox{.} \hbox{.} \hbox{.}}}

 \def \Aitken {1} 
 \def \Arnoldi {2}
 \def \Ashby {3} 
 \def \Brown {4}
 \def \Chronopoulosa {5}
 \def \Chronopoulosb {6}
 \def \Chronopoulosc {7}
 \def \Chronopoulosd {8}
 \def \Concusa {9}
 \def \Concusb {10}
 \def \Craig {11}
 \def \Eisenstata {12} 
 \def \Eisenstatb {13} 
 \def \Elmana {14} 
 \def \Elmanb {15} 
 \def \Engeli {16} 
 \def \Fabera {17} 
 \def \Faberb {18} 
 \def \Faberc {19} 
 \def \Forsythe {20}
 \def \Frankel {21}
 \def \Goluba {22} 
 \def \Golubb {23} 
 \def \Golubc {24} 
 \def \Golubd {25} 
 \def \Golube {26 } 
 \def \Greenbaum {27}
 \def \Hestenesa {28} 
 \def \Hestenesb {29} 
 \def \Householder {30}
 \def \Jea {31}
 \def \Joubert {32}
 \def \Luenbergera {33} 
 \def \Luenbergerb {34} 
 \def \Manteuffela {35} 
 \def \Manteuffelb {36} 
 \def \Manteuffelc {37} 
 \def \Manteuffeld {38} 
 \def \Milne {39}
 \def \Niethammer {40}  
 \def \Opfer {41} 
 \def \Reid {42} 
 \def \Richardson {43}  
 \def \Saad {44}
 \def \Saada {45} 
 \def \Saadb {46} 
 \def \Saadc {47} 
 \def \Shah {48}
 \def \Stiefela {49} 
 \def \Stiefelb {50} 
 \def \Trefethen {51}
 \def \Vargaa {52} 
 \def \Vargab {53} 
 \def \Vinsome {54} 
 \def \Voevodin {55}
 \def \Walkera {56} 
 \def \Walkerb {57} 
 \def \Widlund {58}
 \def \Younga {59} 
 \def \Youngb {60} 
 \def \Youngc {61} 
 \def \Younge {62} 
 \def \Youngd {63} 





 \pageno = 3
 
 {\parindent=0pt \parskip=0pt 

 \line {\hphantom {UC??}\hfil SAND89-8691\hfil \hphantom {UC??}}
 \centerline {Unlimited Release}
 \centerline {Printed November 1989}

 \vglue1.25in
 \centerline {\bigbold Operator Coefficient Methods}
 \vglue0.1in
 \centerline {\bigbold for Linear Equations*\footnote
{}{\noindent \tenrm *~Invited for publication by SIAM Journal on
Matrix Analysis and Applications.}}
 \vglue0.5in
 {\baselineskip=13pt \bold
 \centerline {Joseph F. Grcar} 
 \centerline {Scientific Computing and Applied Math Division}
 \centerline {Sandia National Laboratories}
 \centerline {Livermore, CA 94551-0969 USA}}

 \vglue 1.0in
 \centerline {\medbold Abstract}
 \vglue 0.25in 
{\parindent=20pt New iterative methods for solving linear
equations are presented that are easy to use, generalize good
existing methods, and appear to be faster. The new algorithms mix
two kinds of linear recurrence formulas. Older methods have either
high order recurrence formulas with scalars for coefficients, as in
truncated orthomin, or have 1st order  recurrence formulas with
matrix polynomials for coefficients, as in restarted gcr/gmres.
The new methods include both: high order recurrence formulas and
matrix polynomials for coefficients. These methods provide a
trade-off between recurrence order and polynomial degree that
can be exploited to achieve greater efficiency. Convergence
results are obtained for both constant coefficient and varying
coefficient methods.}

 \vfil \ \eject}

 {

 \ \vfil 
 \centerline {\medbold Acknowledgements}
 \vglue 0.33in

I thank the program and organizing committees of the SIAM
Conference on Sparse Matrices at Salishan Lodge in May 1989 for
choosing an iconoclastic abstract for the keynote address.

I thank Dr.~S.~F.~Ashby for many conversations which clarified the
current literature, Prof.~A.~T.~Chronopoulos and
Prof.~L.~N.~Trefethen for discussing their work, and
Prof.~H.~C.~Elman and Prof.~P.~E.~Saylor for help with
bibliographic matters.

Thanks also to Dr.~S.~F.~Ashby, Dr.~L.~A.~Bertram,
Dr.~J.~M.~Harris, Dr.~J.~F. Lathrop, Mr.~S.~L. Lee, Dr.~J.~C.~Meza,
Dr.~L.~R.~Petzold and Prof.~P.~E.~Saylor for reading the
manuscript. Their patience much improved this paper.

 \vfil \vglue 2.0in \ \eject}

 {\parskip=0pt \parindent=0pt

 \ \vfil 
 \centerline {\medbold Contents}
 \vglue 0.5in
 \def \filled (#1) (#2){\hbox to \hsize {#1\leaders \hbox to
2em{\hss.\hss}\hfill #2}}
 \hfil \vbox {\hsize=4.5in \openup 2\jot
 \filled (1. Introduction) (7)
 \filled (2. Survey) (8)
 \filled (\qquad 2a. Incompletely Specified Methods) (8)
 \filled (\qquad 2b. Completely Specified, Terminating Methods)
(11)
 \filled (\qquad 2c. Completely Specified, Non-Terminating
Methods) (16)
 \filled (3. Simplification) (18)
 \filled (4. Inhomogeneous Methods) (20)
 \filled (5. Operator Coefficient Methods) (23)
 \filled (6. Implementations) (27)
 \filled (\qquad 6a. Some Existing Implementations) (27)
 \filled (\qquad 6b. An Implementation) (30)
 \filled (7. Varying Coefficients) (31)
 \filled (8. Constant Coefficients) (35)
 \filled (References) (43)
 \filled (Appendix 1. Proofs) (47)
 \filled (Appendix 2. Figure Explanations) (53)
 }
 \vglue 2.0in \vfil \eject}

 \ \vfil \eject

\beginsection {1. Introduction}

This paper develops a new class of iterative methods for solving
linear equations. The new algorithms are easy to use, they
generalize good existing methods, and they appear to be faster.

The new class of algorithms combines the essential features of
several methods, in the following way. When many known
algorithms are defined very simply, their sequences of
approximate solutions can be seen to satisfy recurrence formulas
of two types. Either the recurrence formulas have high orders and
scalar coefficients, as in truncated orthomin, or the recurrence
formulas have order $1$ and polynomials of matrices for
coefficients, as in restarted gcr/gmres. The new methods are the
generalization to formulas that mix the two types of recurrences.

Methods in the new class build solutions from linear combinations
of vectors using operators for coefficients, that is, using
polynomials of matrices. The solution sequences therefore are
vector linear recurrences whose coefficients are linear
transformations, whence the name {\sl operator coefficient
methods}. The choice of coefficients leads to many, mostly
unexplored variations. 

A convenient, and in older methods, a frequent choice of
coefficients repetitively solves a simple minimization problem. A
new convergence result proves this choice of coefficients often
results in convergence and establishes upper bounds on the
convergence rates. An examination of the parameters that govern
the convergence rates and some numerical experiments suggest
that the new methods are faster than those currently in use.
Moreover, it is suggested that both old and new algorithms may
be better implemented by means of well-established,
least-squares procedures. 

These are the paper's major results. First, a simple
characterization of iterative methods is proposed that unifies
many algorithms. Similar descriptions are known for some
algorithms, but they have not been systematically applied to
others. Second, a spectrum of new iterative methods is found to
lie between those of high recurrence order such as truncated
orthomin, and those of high polynomial degree such as restarted
gcr/gmres. Third, many of the algorithms are observed to have
identical convergence rates. For the same convergence rate,
there results a trade-off between degree and order that can be
exploited to optimize efficiency. Fourth, new convergence results
are proved for both constant coefficient and varying coefficient
methods.

The paper's organization follows the steps which led to discovery
of the new methods. First, iterative algorithms are surveyed from
the historical point of view. Even those familiar with the subject
may find this survey interesting. Then, some algorithms are
restated in a very simple form. It is suggested this form has
pedagogical and practical advantages. Two new generalizations
come from it. One is a further simplification to a new class of
inhomogeneous methods. Their convergence cannot be explained
by existing analytical tools, nor is it analyzed here. The second
and more important generalization is to operator coefficient
methods. These have both homogeneous and inhomogeneous
forms. A straightforward implementation for one choice of
operator coefficients is proposed, and sufficient conditions for
convergence are found. Finally, necessary and sufficient
conditions for convergence of constant coefficient methods are
found. To improve readability, appendices contain the proofs of
theorems and descriptions of numerical experiments. 

Chronopoulos and Gear [\Chronopoulosa ] [\Chronopoulosb ]
[\Chronopoulosd ] have been led by other considerations to
derive some algorithms in the new class of methods, and their
work-in-progress examines more [\Chronopoulosc ]. Their work
should be consulted for additional insights.

\beginsection {2. Survey}

This section surveys the iterative methods to be analyzed in the
sequel. The survey is mostly historical and phenomenological, with
several omissions and some rearrangement. Those who expect
to read the paper in one sitting or who are familiar with the
subject may prefer to begin at Section~3 and to consult this
section as needed.

{\sl Iterative methods} are prescriptions for building sequences
$x_0$ $x_1$ $x_2$ $\ldots$ $x_n$ $\ldots$ that converge to the
solution of $Ax=y$. If the prescriptions differ but the sequences
are the same then the methods are the same, so some
prescriptions might be better than others for the same method.
The word {\sl prescription} is new in this context.

What appears to be the natural progression of ideas for iterative
algorithms involves the manner of choosing various coefficients
and parameters to build the sequences. The oldest
iterative methods are {\sl incomplete} because their parameters
must be selected from ancillary information. Algorithms in the
second phase of development make automatic choices that are
{\sl globally correct} in some sense for some matrices. Methods in
the third phase choose coefficients that are {\sl suboptimal} but
serviceable for more matrices.

The survey is bounded as follows. The first limitation is to {\sl
polynomial} methods. Each term of the sequence equals a linear
combination of other vectors with coefficients that are
polynomials of $A$. The second limitation is to polynomials {\sl
only} of $A$. If preconditioners are used, then $A$ must be the
result of any preconditioning matrix multiplications. The third
limitation is to amounts of space and time per algorithm step
limited {\sl independent} of the step number. This means not all
the preceding solutions, residuals, or whatever are available to
build the next approximate solution.

The following notational conventions are used throughout. The {\sl
exact solution} of $Ax = y$ is $x_*$. The {\sl error} in $x_n$ is $e_n
= x_* - x_n$. The {\sl residual} of $x_n$ is $r_n = y - Ax_n$. Note
that $A e_n = r_n$. Symbols with negative subscripts equal zero,
and matrix-vector norms are the $2$-norm, unless stated
otherwise. Complex numbers are assumed. The {\sl Hermitian part}
of $A$ is $(\transpose A + A)/2$. The {\sl set of eigenvalues} of the
matrix $H$ is $\lambda (H)$.

\beginsubsection {2a. Incompletely Specified Methods}

\Proclaim Richardson's 1st Order Method, 1910. Iterate from
$x_0$. $$x_{n+1} = \alpha_n r_n + x_n$$

Richardson's paper [\Richardson ] makes interesting reading from
the turn of the century. {\sl Residual polynomials} $$P_n(X) =
\prod_{j=1}^n \, (1 - \alpha_j X)$$ have been used to understand
this and other methods because they provide formulas for the
residuals and the errors. $$r_n = P_n(A) r_0 \qquad e_n = P_n(A)
e_0$$ The many papers such as [\Opfer ] on the proper choice of
coefficients are outside the scope of this survey. If all the
coefficients $\alpha_0$, $\alpha_1$, $\alpha_2$, $\ldots$ equal the
same $\alpha$, then the method converges for all $y$ and $x_0$
exactly when all the eigenvalues of $A$ lie strictly inside the
circle through $0$ around $1/\alpha$ in the complex plane.

The history of iterative methods swings like a pendulum between
two extremes of fashion. On the one side are the basic iterations 
surveyed here, on the other side are things now called
preconditioners. {\sl Preconditioning} replaces $Ax=y$ by
$BAx=By$ where $B$ is an approximate inverse for $A$. The
pendulum started with Richardson's 1st order method, and the
first reversal probably occurred when interest reverted to
relaxation schemes, some of which are much older if their
appellations can be believed. 

The classic preconditioners were
developed for use with the simplest 1st order meth\-od, for
$\alpha_n = 1$ and $x_{n+1} = r_n + x_n$. They {\sl split} the
matrix $A = L + D + U$ into its diagonal and triangular parts, and
choose $B$ as follows. $$\vcenter {\openup2\jot \halign {\hfil
$#\quad$& #\hfil \cr D^{-1}& Jacobi, 1845\cr (D + L)^{-1}&
Gauss-Seidel, 1873\cr ({1 \over \omega} D + L)^{-1}& SOR,
Successive Overrelaxation, 1950\cr {2 - \omega \over \omega} ({1
\over \omega} D + U)^{-1} D ({1 \over \omega} D + L)^{-1}& SSOR,
Symmetric SOR, 1950\cr}}$$ G.~E.~Forsythe said Gauss-Seidel
was not known to Gauss and not recommended by Seidel
[\Householder ], so the very old references in this and [\Shah ]
[\Vargab ] [\Younge ] must be consulted to see how preconditioners
predating Richardson's method, or at least Richardson's
description of his method, were conceived. These and other {\sl
relaxation methods} are now viewed as a class of preconditioners.
SOR was invented independently by Frankel [\Frankel ] and Young
[\Younga ] [\Youngb ], and SSOR by Aitken [\Aitken ]. In the 1950's
and 60's many matrices and preconditioners were found that make
Richardson's method converge. They are described by Golub
[\Goluba ], Golub and Varga [\Golubb ] [\Golubc ], Varga [\Vargab ],
Young [\Youngc ] and many others.

\Proclaim 2nd Order Method,  (1950) 1958. Iterate
from $x_0$ $$x_{n+1} = \alpha_n r_n + \beta_n x_n + \gamma_n
x_{n-1} $$ in which $\beta_0 = 1$ and $\beta_n + \gamma_n = 1$.

Frankel [\Frankel ] invented this method in 1950 and named it
after Richardson. He was deferential to a fault because in the same
paper he invented SOR and named it after Liebmann. Frankel and
many others omit the name of the third coefficient. Stiefel
[\Stiefelb ] makes $\gamma_n = 1 - \beta_n$ appear to be a
natural consequence of normalization. The conditions $\beta_0 =
1$ and $\beta_n + \gamma_n = 1$ enable the following analysis.
With them the familiar residual polynomials exist $$r_n = P_n(A)
r_0 \qquad e_n = P_n(A) e_0$$ and can be built from $$P_0(X) = 1
\qquad P_{n+1}(X) = - \alpha_n X P_n(X) + \beta_n P_n(X) +
\gamma_n P_{n-1}(X)$$ just like normalized orthogonal
polynomials. Stiefel may have been the first to describe the
method's possibilities when he suggested consulting the theory of
orthogonal polynomials to find appropriate coefficients.

\Proclaim Chebyshev Iteration,  (1957) 1975. Iterate from $x_0$
$$x_{n+1} = \alpha_n r_n + \beta_n x_n + (1 - \beta_n) x_{n-1}$$ in
which $$\vcenter {\vglue2ex \openup2\jot \halign {$\displaystyle
#$\hfil \qquad& $\displaystyle #$\hfil \cr \alpha_0 = 1/d& \beta_0
= 1\cr \alpha_n = {2 \, T_n (d/c) \over c \, T_{n+1} (d/c)}& \beta_n
= {2 d \, T_n (d/c) \over c \, T_{n+1} (d/c)}\cr}}$$ where $T_n$ is the
$n^{th}$ Chebyshev polynomial of the first kind.

The Chebyshev iteration is the 2nd order method with residual
polynomials $$P_n(X) = { T_n {\displaystyle \Big( {d - X \over c
\vphantom {|}} \Big)} \over \lower1.75ex \hbox {$T_n
{\displaystyle \Big( {d \over c \vphantom {|}} \Big)}$} } .$$ Varga
[\Vargaa ] derived the method differently by building rapidly
converging sequences from more slowly converging ones.
Manteuffel [\Manteuffela ] [\Manteuffelb ] [\Manteuffelc ] extended
the method beyond symmetric positive definite matrices and
coined the present name. The iteration converges for all $y$ and
$x_0$ exactly when all the eigenvalues of $A$ lie strictly inside
the ellipse through $0$ with foci $d \pm c$ in the complex plane.

\Proclaim Stationary 2nd Order Method, 1982. Iterate
from $x_0$ and~$x_{-1}$. $$x_{n+1} = \alpha r_n + \beta x_n + (1 -
\beta) x_{n-1}$$ 

Iterative methods and linear recurrences are {\sl stationary} when
the coefficients are independent of $n$. The handful of papers on
the stationary 2nd order method seek coefficients that optimize
convergence for a given matrix. The answer to the simpler inverse
question---which matrices converge for a given pair of
coefficients?---can be obtained from [\Manteuffeld ] and
a few napkins. Convergence occurs for all $y$ and $x_0$ exactly
when all the eigenvalues of $A$ lie strictly inside the ellipse
through $0$ with foci $${\beta \over \alpha} \pm {2 \sqrt
{\vphantom {|} \beta - 1} \over \alpha} .$$ Moreover, the Chebyshev
iteration's coefficients for these foci converge to the stationary
coefficients.

\beginsubsection {2b. Completely Specified, Terminating Methods}

By the time 2nd order methods were completely understood, the
pendulum had already swung toward iterations with completely
specified coefficients. The next method is the namesake for the
entire class.

\Proclaim Conjugate Gradient Algorithm, (1952) 1971. Iterate
from $x_0$. $$p_n = r_n - {\transpose {p_{n-1}} A r_n \over
\transpose {p_{n-1}} A p_{n-1}} \, p_{n-1}$$ $$x_{n+1} = x_n +
{\transpose {p_n} r_n \over \transpose {p_n} A p_n} \, p_n$$

Hestenes and Stiefel [\Hestenesa ] drew this method from
optimization theory. If $A$ is Hermitian and positive definite,
then the method searches from $x_n$ along the {\sl direction
vector} $p_n$ for the $x_{n+1}$ that minimizes $$\| e_{n+1} \|_A =
\sqrt {\transpose {e_{n+1}} A e_{n+1}} \; .$$ It happens that the
$p$'s are $A$-orthogonal, the $r$'s are orthogonal, and $x_{n+1}$ is
the global minimizer within $$\eqalign {&x_0 + \Span {p_0 \;\; p_1
\;\; p_2 \;\; \ldots \;\; p_n} = {}\cr &\quad x_0 + \Span {r_0 \;\; r_1
\;\; r_2 \;\; \ldots \;\; r_n} = {}\cr &\qquad x_0 + \Span {r_0 \;\;
Ar_0 \;\; A^2 r_0 \;\; \ldots \;\; A^n r_0}.}$$ When the subspaces
stop growing the last iterate is the exact solution. Hestenes
[\Hestenesb ] derives many algebraic identities including the few
needed to establish global optimality and alternate expressions
for the coefficients. Golub and O'Leary [\Golubd ] provide an
excellent annotated bibliography for the huge corpus. 

Elaborate formulas were a disadvantage on early computers so
for many years the conjugate gradient algorithm was seen as a
freakish alternative to Gaussian elimination [\Householder ]. By
1971 technological improvements enabled Reid [\Reid ] to view
the algorithm as an iterative method and to obtain acceptable
solutions after comparatively few steps. A curious tribute to Reid
is that his paper is no longer read because his ideas are so
completely accepted.

The matrices for which the conjugate gradient iteration finds an
exact solution for all $y$ and $x_0$ are the {\sl terminating class}.
This terminology is new. Preconditioning by $\transpose {A}$ yields
$\transpose {A} A x = \transpose {A} y$ with $\transpose {A} A$ in
the terminating class but with squared condition number. Faber and
Manteuffel [\Fabera ] [\Faberb ] [\Faberc ] answered a challenge of
G.~H.~Golub and found the terminating classes for many polynomial
methods based on $A$ alone. The Russian literature contains a
related announcement at about the same time [\Voevodin ].
Unfortunately, all the classes are severely restricted. If the
Hermitian part of $A$ is positive definite, then Joubert and Young
[\Joubert ] show from the work of Faber and Manteuffel that the
conjugate gradient algorithm terminates for all $y$ and $x_0$
exactly when either $\transpose {A} = P(A)$ where $P(X)$ is a
polynomial of degree at most $1$, or $P(A) = 0$ where $P(X)$ is a
nonzero polynomial of degree at most $2$. The direction vectors'
$A$-orthogonality must be reinterpreted in the non-Hermitian
case.

\Proclaim Conjugate Residual Algorithm, (1955) 1970. Iterate from
$x_0$. $$p_n = r_n - {\transpose {p_{n-1}} \transpose A A r_n \over
\transpose {p_{n-1}} \transpose A A p_{n-1}} \, p_{n-1}$$ $$x_{n+1}
= x_n + {\transpose {p_n} \transpose A r_n \over \transpose {p_n}
\transpose A A p_n} \, p_n$$

This algorithm has an interesting genealogy. Hestenes and Stiefel
allude to it [\Hestenesa ], but Stiefel describes it fully without
naming it [\Stiefela ], and Luenberger finally names it when he
reinvents it [\Luenbergera ]. The algorithm is one of many
variations of the conjugate gradient algorithm with similar
properties. This one uses a different inner product. If $A$ is
Hermitian and positive definite, then the $p$'s are $\transpose A
A$-orthogonal, the $r$'s are $A$-orthogonal, and $x_{n+1}$ is the
global minimizer of $\| r_{n+1} \|$ within $$\eqalign {&x_0 + \Span
{p_0 \;\; p_1 \;\; p_2 \;\; \ldots \;\; p_n} = {}\cr &\quad x_0 + \Span
{r_0 \;\; r_1 \;\; r_2 \;\; \ldots \;\; r_n} = {}\cr &\qquad x_0 + \Span
{r_0 \;\; Ar_0 \;\; A^2 r_0 \;\; \ldots \;\; A^n r_0}.}$$ Remember the
convention that unspecified norms are the $2$-norm. Joubert and
Young [\Joubert ] show from the work of Faber and Manteuffel
[\Fabera ] [\Faberb ] [\Faberc ] that among matrices whose
Hermitian part is positive definite, the conjugate residual
algorithm has the same terminating class as the conjugate
gradient algorithm.

\Proclaim Alternate Conjugate Residual Algorithm, 1951. Iterate
from $x_0$ $$p_n = A p_{n-1} - {\transpose {p_{n-1}} \transpose A
A^2 p_{n-1} \over \transpose {p_{n-1}} \transpose A A p_{n-1}} \,
p_{n-1} - {\transpose {p_{n-2}} \transpose A A^2 p_{n-1} \over
\transpose {p_{n-2}} \transpose A A p_{n-2}} \, p_{n-2} $$ $$x_{n+1}
= x_n + {\transpose {p_n} \transpose A r_n \over \transpose {p_n}
\transpose A A p_n} \, p_n$$ but choose $p_0$ as in the original
conjugate residual algorithm.

This method has been invented by many, but Forsythe, Hestenes
and Rosser [\Forsythe ] appear to be the first [\Craig ] [\Golubd ].
It is the conjugate residual method with a different prescription
whose terminating class is larger. Faber and Manteuffel [\Faberb ]
[\Faberc ] show the method converges for all $y$ and $x_0$ exactly
when either $\transpose {A} = P(A)$ where $P(X)$ is a polynomial of
degree at most $1$, or $P(A) = 0$ where $P(X)$ is a nonzero
polynomial of degree at most $3$. These are weaker conditions
than for the conjugate residual algorithm because $3$ replaces $2$,
and more importantly, the Hermitian part of $A$ need not be
positive definite. And yet the $p$'s are $\transpose A
A$-orthogonal, and $x_{n+1}$ is the global minimizer of $\| r_{n+1}
\|$ within $$\eqalign {&x_0 + \Span {p_0 \;\; p_1 \;\; p_2 \;\; \ldots
\;\; p_n}\cr &\quad {} = x_0 + \Span {r_0 \;\; Ar_0 \;\; A^2 r_0 \;\;
\ldots \;\; A^n r_0}}$$ just like the conjugate residual algorithm.

The relative merits of the two prescriptions are not clear. The
original fails for an indefinite matrix when a direction vector
makes no contribution to the solution. In this unlikely event the
residuals do not change and the subsequent direction vectors lie
within the span of the previous. The alternate version succeeds
because it builds the direction vectors from a self-contained
recurrence. In some sense the alternate prescription is an analytic
continuation of the original. But the theoretically more powerful
prescription amounts to evaluating the Lanczos recurrence, which
may be more sensitive to rounding errors.

\pageinsert {\Picture 1. 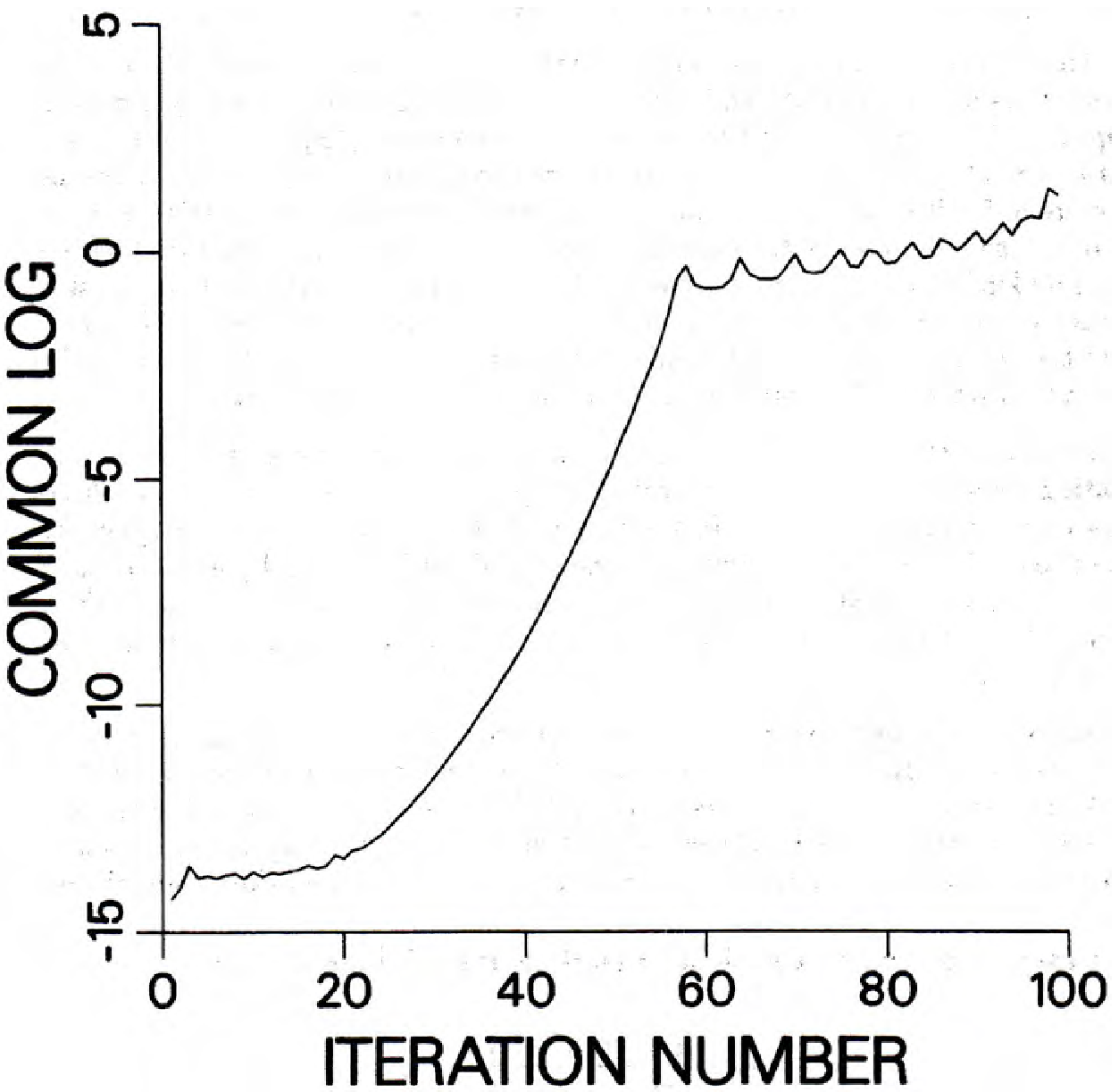 $2$-norm relative errors in the computed
basis vectors of the conjugate gradient algorithm for a system of
order $100$. Appendix~2 and Section~2b explain the calculations.
\par }\endinsert

The effects of rounding errors are a major disappointment for all
these algorithms. Loss of orthogonality and failure to terminate
are the most easily observed symptoms. The conjugate gradient
algorithm should have $\transpose {p_i} A p_j = 0$ for $i \ne j$, but
in practice $\transpose {p_i} A p_j / ( \| p_i \|_A \| p_j \|_A)$ grows
exponentially with $|i-j|$. Figure~1 makes the more difficult
comparison between the numerically computed direction vectors
and those that would be obtained from error-free arithmetic. This
data may be the first of its kind in print. The formulas for the
direction vectors evidently are unstable because they magnify the
small perturbations due to rounding error. Figure~2 shows the
resulting delayed convergence. Greenbaum has a detailed analysis
of the retarded convergence [\Greenbaum ], but no universal,
inexpensive cure is known. 

\pageinsert {\Picture 2. 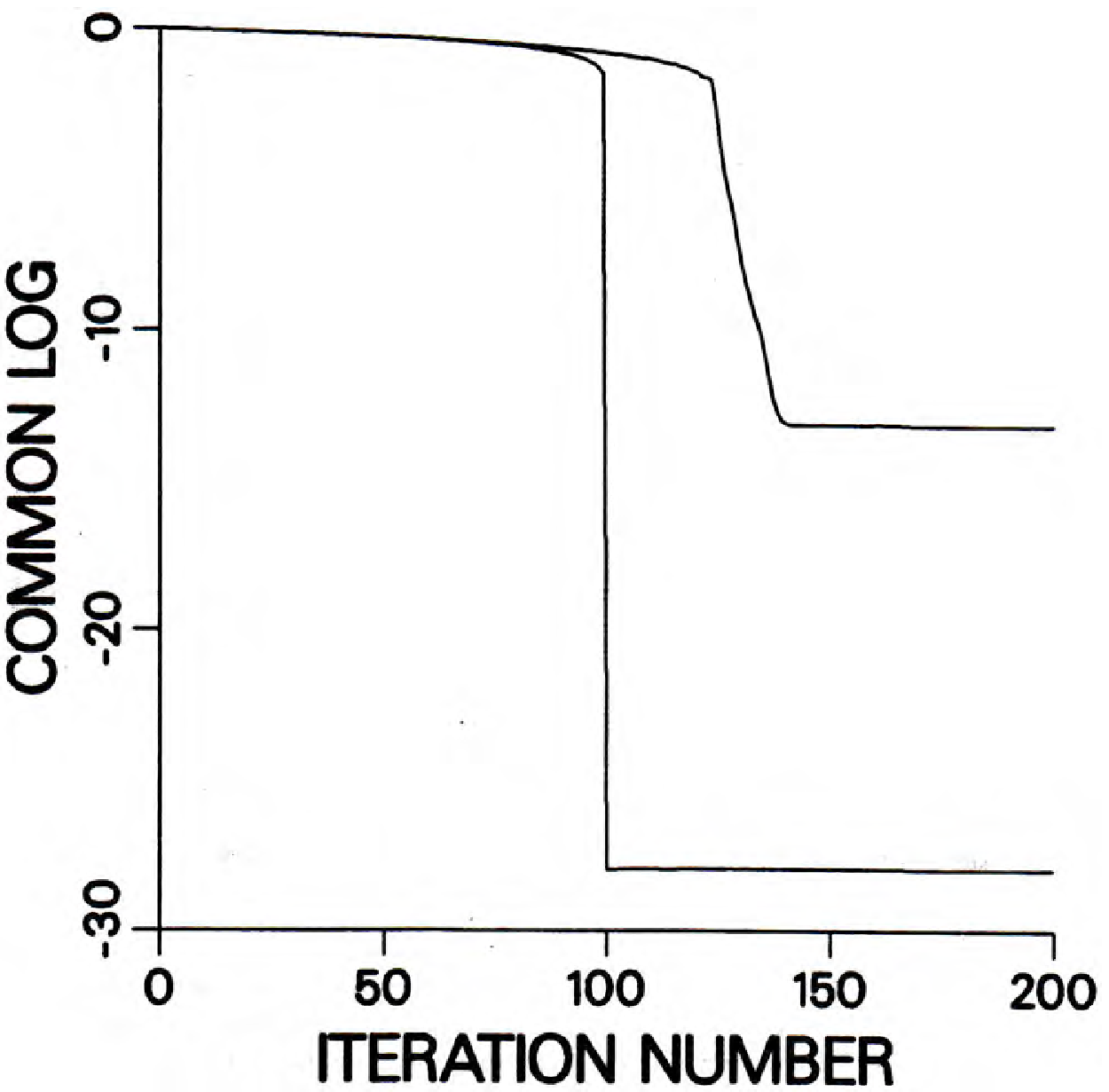 Relative $A$-norm solution errors for the
system of Figure~1. The upper curve is for single precision and the
lower for reorthogonalized double precision. Appendix~2 and
Section~2b explain the calculations.\par }\endinsert

Interest in the conjugate gradient algorithm was intense for a
time. Many terminating algorithms were proposed with enlarged
terminating classes. {\sl Generalized} refers indiscriminately to
these methods for which there is no consistent naming convention.
The generalized conjugate gradient algorithm [\Concusa ]
[\Concusb ] [\Widlund ] established the use of very different inner
products and to some extent prompted the work of Faber and
Manteuffel. Ashby, Manteuffel and Saylor classify
many generalizations of this kind~[\Ashby ]. 

Some terminating generalizations of the conjugate gradient
algorithm have mostly theoretically use. The analyses of Faber and
Manteuffel and of Joubert and Young actually proceed by seeking
conditions under which the original algorithm is equivalent to a
generalized one that uses all the direction vectors to build the
next. The same generalization can be made of the conjugate
residual algorithm. 

\Proclaim Generalized Conjugate Residual Algorithm, 1982. Iterate
from $x_0$. $$p_n = r_n  - \sum_{i=1}^{n} {{\transpose {p_{n-i}}
\transpose A A r_{n} \over \transpose {p_{n-i}} \transpose A A
p_{n-i}} \, p_{n-i}}$$ $$x_{n+1} = x_n + {\transpose {p_n} \transpose
A r_n \over \transpose {p_n} \transpose A A p_n} \, p_n$$

Formulas of this kind may be traced to Arnoldi [\Arnoldi ], but this
algorithm is from Elman [\Elmana ] and Eisenstat, Elman and Schultz
[\Eisenstatb ]. Its terminating class includes all matrices with
positive definite Hermitian parts. The method itself is impractical,
and is outside the bounds of the present survey, because all the
previous direction vectors must be saved. It can be made practical
either by {\sl truncating} the sum, or by {\sl restarting} the
iteration, as follows.

\beginsubsection {2c. Completely Specified, Non-Terminating
Methods}

\Proclaim Orthomin(m), 1976. Iterate from $x_0$. $$p_n = r_n  -
\sum_{i=1}^m {{\transpose {p_{n-i}} \transpose A A r_n \over
\transpose {p_{n-i}} \transpose A A p_{n-i}} \, p_{n-i}}$$ $$x_{n+1}
= x_n + {\transpose {p_n} \transpose A r_n \over \transpose {p_n}
\transpose A A p_n} \, p_n$$

Vinsome [\Vinsome ] invented this {\sl truncated} algorithm and
made non-terminating or suboptimal convergence respectable
again. Elman [\Elmana ] and Eisenstat, Elman and Schultz
[\Eisenstatb ] prove convergence for all $y$ and $x_0$
whenever the Hermitian part of $A$ is positive definite.
Specifically, they prove $$\| r_{n+1} \; \| \le \; \| r_n \| \, \sqrt {1 -
\left[ \, {\min | \lambda (\transpose A + A) | \over 2 \| A \|} \,
\right]^2} \; < \; \| r_n \|.$$ They also show each set of $m+1$
consecutive direction vectors is $\transpose A
A$-or\-thog\-o\-nal, and $x_{n+1}$ minimizes $\| r_{n+1} \|$ within
$$x_{{n-m} \wedge 0} + \Span {p_{n-m} \;\; \ldots \;\; p_{n-2} \;\;
p_{n-1} \;\; p_n} .$$ The notation is rigorously correct because
things with negative subscripts vanish and the wedge in $x_{{n-m}
\wedge 0}$ means the maximum of $n-m$ and $0$. It is an open
question why the suboptimal convergence results ignore $m$ and
exclude the Hermitian indefinite case for which the alternate
conjugate residual algorithm has terminating convergence.

The truncated algorithm is expected to have dependable
convergence for many matrices rather than terminating
convergence for a few. As $m$ increases the terminating class
grows beyond that of the conjugate residual method, but only
through the addition of matrices having at most $m^2$ distinct
eigenvalues [\Faberb ] [\Faberc ] [\Joubert ]. Termination also
occurs, of course, for $m$ so impractically large that the method
reverts to the generalized conjugate residual method.

As with the conjugate residual algorithm, there is an alternate
version that generates the direction vectors independently of the
residuals. Saad and Schultz survey many equivalent prescriptions
[\Saadb ]. Saad develops some of these algorithms himself, and
names the class {\sl incomplete orthogonalization methods} [\Saad
]. Jea and Young [\Jea ] also develop a broad class of methods and
also provide extensive references to other work. In their
terminology, the original conjugate residual algorithm is an {\sl
orthomin} and the alternate is an {\sl orthodir}, both with $m = 1$
and with specific choices for inner products and the like. Faber and
Manteuffel [\Fabera ] [\Faberb ] [\Faberc ] actually treat the
alternate prescriptions, but Joubert and Young [\Joubert ] and
Ashby, Manteuffel and Saylor [\Ashby ] carefully observe the
distinction.

\Proclaim Gcr(k), Restarted Generalized Conjugate Residual
Meth\-od, 1979. It\-er\-ate from $x_0$, and for each $x_n$ obtain
$x_{n+1}$ by building the sequence $$x_n = x_{(0)} \quad x_{(1)}
\quad x_{(2)} \quad \ldots \quad x_{(k+1)} = x_{n+1}$$ by iterating
the generalized conjugate residual algorithm from $x_{(0)}$ to
$x_{(k+1)}$ $$p_{(j)} = r_{(j)} - \sum_{i=1}^{j} {\transpose
{p_{(i-1)}} \transpose A A r_{(j)} \over \transpose {p_{(i-1)}}
\transpose A A p_{(i-1)}} \, p_{(i-1)}$$ $$x_{(j+1)} = x_{(j)} +
{\transpose {p_{(j)}} \transpose A r_{(j)} \over \transpose {p_{(j)}}
\transpose A A p_{(j)}} \, p_{(j)}$$ in which the subscripts in
parentheses indicate dependence on $n$.

The idea of {\sl restarting} an algorithm contrasts with truncating
in orthomin($m$). It is the theory that jiggling the ignition
recharges the battery and may be due to several people.
Luenberger [\Luenbergerb ] restarts the conjugate gradient
algorithm to circumvent numerical difficulties, while Eisenstat,
Elman, Schultz and Sherman [\Eisenstata ] [\Eisenstatb ] [\Elmana ]
restart the generalized conjugate residual method to conserve
memory space. They show that gcr($k$) minimizes  $\| r_{n+1} \|$
within $$\eqalign {&x_n + \Span {r_n \;\; Ar_n \;\; A^2 r_n \;\; \ldots
\;\; A^k r_n}}$$ and converges when the Hermitian part of $A$ is
positive definite. 

Like truncated orthomin, restarted gcr is expected to have
dependable convergence for many matrices rather than
terminating convergence for a few. And again there is an alternate
prescription. This one additionally makes the direction vectors
orthogonal with respect to the Euclidean inner product rather than
the $\transpose A A$ inner product.

\Proclaim Gmres(k), Restarted Generalized Minimum Residual
Algorithm, 1983. It\-er\-ate from $x_0$, and for each $x_n$ build
the orthonormal sequence $${r_n \over \| r_n \|} = p_{(1)} \quad
p_{(2)} \quad p_{(3)} \quad \ldots \quad p_{(k+1)}$$ from Arnoldi's
recurrence equations $$h_{(j+1,\,j)} \, p_{(j+1)} = A \, p_{(j)} -
\sum_{i=1}^j h_{(i,\,j)} \, p_{(i)}$$ with appropriately chosen
$h_{(i,\,j)}$'s, and then chose  $\alpha_{(j)}$'s to minimize $$\left\| \;
\| r_n \| \, e_1 - [ h_{(i,\,j)} ] \, [ \alpha_{(j)} ] \; \right\|$$ in which
$e_1$ is the first column of an identity matrix and $[ h_{(i,\,j)} ]$ is
the $(k+1) \times k$ matrix of recurrence coefficients, and finally
construct $$x_{n+1} = x_n + \sum_{j=1}^k \alpha_{(j)} \, p_{(j)} .$$
The subscripts in parentheses indicate dependence on $n$.

Saad and Schultz [\Saada ] [\Saadc ] developed this most widely
used version of restarted gcr. It uses an alternate prescription to
generate the normalized direction vectors, with the precise
choice of $h_{(i,\,j)}$'s clear and mercifully omitted. The
coefficients for $x_n$ are selected from the small matrix that
describes the action of $A$ on the orthonormal basis. Nevertheless,
gmres($k$) minimizes $\| r_{n+1} \|$ within $$\eqalign {&x_n + \Span
{r_n \;\; Ar_n \;\; A^2 r_n \;\; \ldots \;\; A^{k-1} r_n}}$$ just like
gcr($k{-}1$). Both prescriptions sometimes have difficulty
solving the least squares problem, and other prescriptions have
been proposed [\Walkera ] [\Walkerb ].

Note the multiple names. Gcr($k{-}1$) is gmres($k$). The first
generalizes the conjugate residual algorithm and the second
generalizes the minimum residual algorithm. So there is yet
another line of development which leads to the same methods.
But this is too broad a subject for discussion.

\beginsection {3. Simplification}

The completely specified algorithms in the survey can be reduced to
simpler but equivalent form. Here, all inessential notation is
removed to leave what may be the vital core. This naive approach
leads to useful generalizations in subsequent sections, and even to
useful implementations. The simplified algorithms are strikingly
similar. Each chooses its next solution from a small {\sl selection
space}, and uses a minimization problem as the {\sl selection
criterion}. 

Another form of the conjugate gradient algorithm discards the
direction vectors and reveals it to be a 2nd order method. $$x_{n+1}
= \alpha_n r_n + \beta_n x_n + (1 - \beta_n) x_{n-1}$$ This version
has several sources. One builds the coefficients recursively, and is
attributed to Engeli, Ginsberg, Rutishauser and Stiefel [\Engeli ] by
[\Youngd ], and to Rutishauser alone by [\Reid ]. Hestenes
[\Hestenesb ] cites a form called {\sl paratan} with a geometric
interpretation and explicit coefficient formulas [\Shah ]. $$\alpha_n
= {\| r_{n-1} \|^2 \| r_n \|^2 \over \| r_{n-1} \| ^2 ( \transpose r_n A r_n
) + \| r_n \|^2 ( \transpose r_{n-1} A r_n )}$$ $$\beta_n = {\| r_{n-1}
\|^2 ( \transpose r_n A r_n ) \over \| r_{n-1} \| ^2 ( \transpose r_n A
r_n ) + \| r_n \|^2 ( \transpose r_{n-1} A r_n )}$$ These formulas
too can be discarded because theorems say they make $\| e_{n+1}
\|_A$ globally minimal, and therefore locally minimal. In this way a
simple minimization criterion concisely replaces many elaborate
formulas. This interpretation succinctly characterizes both the
conjugate gradient algorithm and its cousin.

\Proclaim Simplest Conjugate Gradient Algorithm. Iterate from
$x_0$. $$\vcenter {\centerline {minimize $\| e_{n+1} \|_A$ over
$\Span {r_n \;\; x_n \;\; x_{n-1}}$} \smallskip \nobreak \centerline {so
coefficients of $x_n$ and $x_{n-1}$ sum to $1$}}$$

\Proclaim Simplest Conjugate Residual Algorithm. Iterate from
$x_0$. $$\vcenter {\centerline {minimize $\| r_{n+1} \|_2$ over $\Span
{r_n \;\; x_n \;\; x_{n-1}}$} \smallskip \nobreak \centerline {so
coefficients of $x_n$ and $x_{n-1}$ sum to $1$}}$$

The starting point for a simpler version of orthomin($m$) $$p_n
= r_n - \sum_{i=1}^m {{\transpose {p_{n-i}} \transpose A A r_{n}
\over \transpose {p_{n-i}} \transpose A A p_{n-i}} \, p_{n-i}}$$
$$x_{n+1} = x_n + {\transpose {p_n} \transpose A r_n \over
\transpose {p_n} \transpose A A p_n} \, p_n$$ is the work of Elman
[\Elmana ] and Eisenstat, Elman and Schultz [\Eisenstatb ]. They show
if the Hermitian part of $A$ is positive definite, then $x_{n+1}$
minimizes $\| r_{n+1} \|$ within $$x_{{n-m} \wedge 0} + \Span {p_{n-m}
\;\; \ldots \;\; p_{n-2} \;\; p_{n-1} \;\; p_n}$$ and the coefficient in the
formula for $x_{n+1}$ can't vanish. The notation is rigorously correct
because things with negative subscripts vanish and the wedge in
$x_{{n-m} \wedge 0}$ means the maximum of $n-m$ and $0$. The
formula for $p_n$ means $r_n$ can replace $p_n$ inside the span. The
formula for $x_{n+1-j}$ means $(x_{n+1-j} - x_{n-j})$ can replace
$p_{n-j}$. With these substitutions the affine space becomes
$$x_{{n-m} \wedge 0} + \hbox {span} \, \left\{ \; \tableau {(x_{n+1-m}
- x_{n-m})& \ldots& (x_{n-1} - x_{n-2})& (x_n - x_{n-1})&
r_n\cr} \; \right\}$$ from which $(x_{n+1-j} - x_{n-j})$ vanishes if $n-j
< 0$. The $x$-coordinates of each vector in the span sum to zero, and
the affine space adds $x_{{n-m} \wedge 0}$, so the $x$-coordinates of
each vector in the affine space sum to $1$. 

\Proclaim Simplest Orthomin(m). Iterate from $x_0$ $$\vcenter
{\centerline {minimize $\| r_{n+1} \|$ over $\Span {r_n \;\; x_n \;\;
x_{n-1} \;\; x_{n-2} \;\; \ldots \;\; x_{n-m}}$} \smallskip \nobreak
\centerline {so coefficients of $x_n \;\; x_{n-1} \;\; x_{n-2} \;\; \ldots
\;\; x_{n-m}$ sum to $1$}}$$

This definition of orthomin($m$) may be new, but the simplest form
of gcr($k{-}1$)/\discretionary{}{}{}gmres($k$) is well known and needs no derivation. Table~1 allows
side-by-side comparison of the simplest versions of all these
methods for the first time.

\midinsert {\table 1. {Simplest prescriptions for the survey's
completely specified algorithms. Section~3 provides further
explanation.} {\vbox {\halign {#\hfil\cr \noalign {\medskip \hrule
\medskip} conjugate gradient\cr \noalign {\smallskip} \qquad
minimize $\| e_{n+1} \|_A$\cr \qquad over $\Span {r_n \;\, x_n \;\,
x_{n-1}}$\cr \qquad so coefficients of $x_n$ and $x_{n-1}$ sum to
$1$.\cr \noalign {\medskip \hrule \medskip} conjugate residual\cr
\noalign {\smallskip} \qquad minimize $\| r_{n+1} \|_2$\cr \qquad
over $\Span {r_n \;\, x_n \;\, x_{n-1}}$\cr \qquad so coefficients of
$x_n$ and $x_{n-1}$ sum to $1$\cr \noalign {\medskip \hrule
\medskip} orthomin($m$)\cr \noalign {\smallskip} \qquad minimize
$\| r_{n+1} \|_2$\cr \qquad over $\Span {r_n \;\, x_n \;\, x_{n-1} \;\,
x_{n-2} \;\, \ldots \;\, x_{n-m}}$\cr \qquad so coefficients of $x_n \;\,
x_{n-1} \;\, x_{n-2} \;\, \ldots \;\, x_{n-m}$ sum to $1$\cr \noalign
{\medskip \hrule \medskip}
gcr($k{-}1$)/\discretionary{}{}{}gmres($k$)\cr \noalign {\smallskip}
\qquad minimize $\| r_{n+1} \|_2$\cr \qquad over $\Span {x_n \;\, r_n
\;\, Ar_n \;\, A^2 r_n \;\, \ldots \;\, A^{k-1} r_n}$\cr \qquad so
coefficient of $x_n$ equals $1$\cr \noalign {\medskip \hrule }}}}\par
}\endinsert

These versions are proposed as archetypes for study and use. Each
prescription is a simple linear recurrence with coefficients that
repetitively solve a simple minimization problem. Confusion and
duplication of effort are unlikely because identity and functionality
are clear at a glance. Straightforward solution of the minimization
problems affords easy comparison and substitution of methods.

\beginsection {4. Inhomogeneous Methods}

New methods can be derived by further simplifying the algorithms of
Table~1. The resulting algorithms apparently cannot be analyzed by
traditional theory, nor is a new theory offered here. These
algorithms do simplify Section~5's presentation of more important
generalizations.

A common feature of all the algorithms in Table~1 is the constraint
that the $x$-coefficients sum to $1$. That is, the next solution
equals a linear combination of previous solutions and other things, in
which the coefficients of the previous solutions sum to $1$. This
constraint is called the {\sl consistency condition}, but {\sl
homogeneity condition} more accurately describes its use.
With it, the formula for the next solution can be multiplied by $A$
and subtracted from $y$ to make an homogeneous recurrence for the
residuals, and this can be multiplied by $A^{-1}$ to make a similar
recurrence for the errors. If the recurrence formulas are used to
produce polynomials rather than vectors, then all the residuals and
errors can be obtained formally, by multiplying the initial residual
and error by these so-called residual polynomials evaluated at
the matrix $A$. 

The entire convergence theory of iterative methods rests on
residual polynomials. Convergence to the solution of $Ax = y$
depends on both $A$ and $y$, but the homogeneity condition allows a
separation of variables in which the entries of $A$ are more
prominent than the entries of $y$. Convergence is equivalent to the
residual polynomials having small values at the matrix eigenvalues.
The incompletely specified algorithms in Section~2a need
parameters that make the polynomials small independent of $y$. The
completely specified algorithms in Sections~2b and~2c choose
parameters that make the polynomials small in norms weighted by
the entries of $y$. In both cases, convergence depends strongly on
$A$ and weakly on $y$.

When the algorithms are stated so simply as in Table~1, however,
there is clearly no reason to impose homogeneity. It is a theoretical
convenience for convergence analysis that is superfluous to the
algorithms. Completely new algorithms can be derived by
removing the constraint. Table~2 presents these even simpler, {\sl
inhomogeneous} algorithms. Henceforth, the original algorithms are
called {\sl homogeneous}.

\midinsert {\table 2. {Inhomogeneous, simplest prescriptions for the
survey's completely specified methods. Section~4 provides further
explanation.} {\vbox {\halign {#\hfil\cr \noalign {\hrule \medskip}
un-conjugate gradient\cr \noalign {\smallskip} \qquad minimize $\|
e_{n+1} \|_A$\cr \qquad over $\Span {r_n \; x_n \; x_{n-1}}$\cr
\noalign {\medskip \hrule \medskip} un-conjugate residual\cr
\noalign {\smallskip} \qquad minimize $\| r_{n+1} \|_2$\cr \qquad
over $\Span {r_n \; x_n \; x_{n-1}}$\cr \noalign {\medskip \hrule
\medskip} un-orthomin($m$)\cr \noalign {\smallskip} \qquad
minimize $\| r_{n+1} \|_2$\cr \qquad over $\Span {r_n \; x_n \; x_{n-1}
\; x_{n-2} \; \ldots \; x_{n-m}}$\cr \noalign {\medskip \hrule
\medskip} un-gcr($k{-}1$)/gmres($k$)\cr \noalign {\smallskip}
\qquad minimize $\| r_{n+1} \|_2$\cr \qquad over $\Span {x_n \; r_n \;
Ar_n \; A^2 r_n \; \ldots \; A^{k-1} r_n}$\cr \noalign {\medskip \hrule
}}}}\par }\endinsert

\pageinsert {\Picture 3. 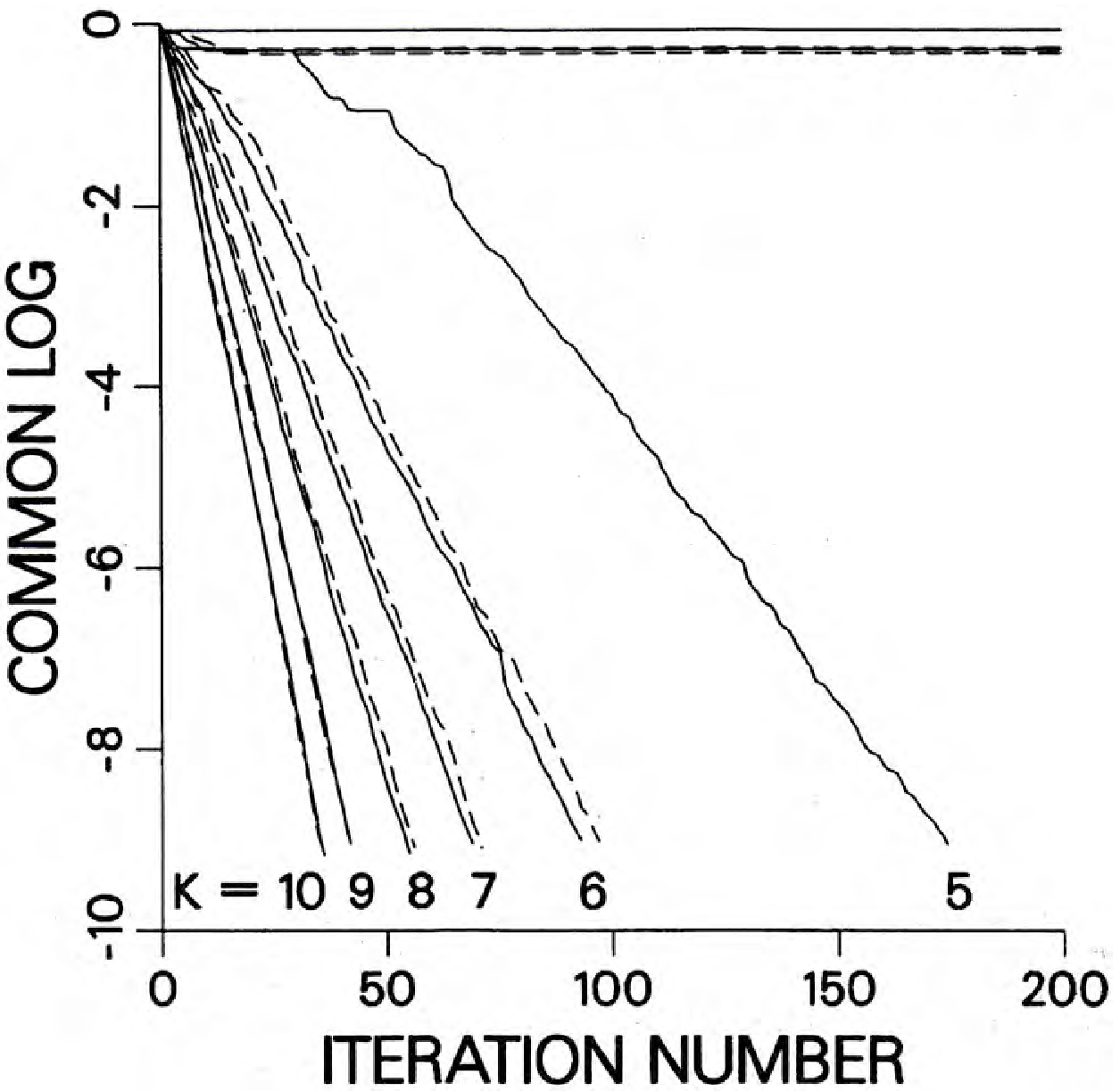 $2$-norm relative residuals for homogeneous
(dashed) and inhomogeneous (solid)
gcr\/$(k{-}1)$/\discretionary{}{}{}gmres\/$(k)$, $k = 1$, $2$,
$\ldots\,$, $10$, applied to one system. The two methods perform
alike except for $k = 5$ when the original, homogeneous method
stagnates and the new, inhomogeneous method converges. Appendix
2 and Section~4 explain the calculations. \par }\endinsert

Figure~3 shows that the new, inhomogeneous methods may converge
when the old, homogeneous methods do not. The selection criteria
evidently find smaller minima when the selection spaces grow by
removing the homogeneity constraint. This explanation is too
simple, however, because it does not characterize the new
convergence rate. This difficult question is not addressed here
beyond the following comments. First, terminating algorithms
already make globally optimal choices, so removing the constraint
should not change them, at least in exact arithmetic. Second, if the
recurrence coefficients of non-terminating algorithms converge to
constant values, then the recurrence formulas must be homogeneous
in the limit. In particular, algorithms with constant coefficients must
be homogeneous.

\beginsection {5. Operator Coefficient Methods}

This paper's major observation is that new iterative methods can be
derived by combining the algorithms of Table~2. The easiest way to
join the algorithms is to amass their selection spaces, as follows. 

The basis vectors naturally fit into a tableau. Those of the conjugate
gradient algorithm and the conjugate residual algorithm occupy a
corner, those of or\-tho\-min($m$) add a row, and those of
gcr($k{-}1$)/\discretionary{}{}{}gmres($k$) fill a column. \def \hide
{\phantom} $$\tableau { x_n& x_{n-1}& x_{n-2}& \;\; \ldots \;\; &
x_{n-m}\cr r_n& \hide {r_{n-1}}& \hide {r_{n-2}}& \hide {\;\; \ldots
\;\;}& \hide {r_{n-k}}\cr Ar_n& \hide {Ar_{n-1}}& \hide {Ar_{n-2}}&
\hide {\;\; \ldots \;\;}& \hide {Ar_{n-k}}\cr \Vdots \;\;& \hide {\Vdots
\;\;}& \hide {\; \Vdots \;\;\;}& \hide {\; \Ddots \;\;}& \hide {\Vdots
\;\;}\cr A^{k-1}r_n& \hide {A^{k-1}r_{n-1}}& \hide {A^{k-1}r_{n-2}}&
\hide {\;\; \ldots \;\;}& \hide {A^{k-1}r_{n-k}}\cr }$$
Orthomin($m$) apparently gains its advantage over the conjugate
residual algorithm by keeping more old solutions. It is likely the
vectors of gcr($k{-}1$)/\discretionary{}{}{}gmres($k$) can be kept
with some advantage too. The tableau does have room for many
more. \def \hide {} $$\tableau {x_n& x_{n-1}& x_{n-2}& \;\; \ldots \;\;
& x_{n-m}\cr r_n& \hide {r_{n-1}}& \hide {r_{n-2}}& \hide {\;\; \ldots
\;\;}& \hide {r_{n-m}}\cr Ar_n& \hide {Ar_{n-1}}& \hide {Ar_{n-2}}&
\hide {\;\; \ldots \;\;}& \hide {Ar_{n-m}}\cr \Vdots \;\;& \hide {\Vdots
\;\;}& \hide {\; \Vdots \;\;\;}& \hide {\; \Ddots \;\;}& \hide {\Vdots
\;\;}\cr A^{k-1}r_n& \hide {A^{k-1}r_{n-1}}& \hide {A^{k-1}r_{n-2}}&
\hide {\;\; \ldots \;\;}& \hide {A^{k-1}r_{n-m}}\cr}$$ The remainder of
this paper demonstrates that better iterative methods can be
created by placing some or all of these vectors into the selection
spaces. The following definition of the new methods involves a
change of notation because the tableau loses one column.

\Proclaim Oc(k,$\;$m), Operator Coefficient Methods of Degree k and
Order m. Be\-gin from $x_0$ and optionally from $x_{-1}$, $x_{-2}$,
$\ldots\,$,  $x_{1-m}$, and choose $x_n$ to $$\hbox {minimize $\| r_n
\|_{\hbox {whatever}}$ or whatever}$$ from among \def \hide {}
$$\hbox {span} \; \left\{ \tableau {x_{n-1}& x_{n-2}& \;\; \ldots \;\; &
x_{n-m}\cr \hide {r_{n-1}}& \hide {r_{n-2}}& \hide {\;\; \ldots \;\;}&
\hide {r_{n-m}}\cr \hide {Ar_{n-1}}& \hide {Ar_{n-2}}& \hide {\;\;
\ldots \;\;}& \hide {Ar_{n-m}}\cr \hide {\;\; \Vdots \;\;}& \hide {\;\;
\Vdots \;\;}& \hide {\; \Ddots \;\;}& \hide {\;\; \Vdots \;\;}\cr \hide
{A^{k-1}r_{n-1}}& \hide {A^{k-1}r_{n-2}}& \hide {\;\; \ldots \;\;}&
\hide {A^{k-1}r_{n-m}}\cr } \; \right\}$$ and, in the homogeneous
case, choose $x_n$ so the $x$-coefficients sum to $1$.

The definition above introduces a name for a generic class of old and
new algorithms. Three aspects need further explanation. First, it isn't
necessary to employ all the vectors in the span. Some old
algorithms do not. Second, the selection criteria is unspecified
because there are so many possibilities. Some are explored in later
sections. Third, since the recurrence formula has order $m$, it is
possible to begin from $m$ initial guesses, $x_{0}$, $x_{-1}$,
$x_{-2}$, $\ldots\,$,  $x_{1-m}$. In this case the operator coefficient
method is not a Krylov space method.

Operator coefficient methods include many known iterative
algorithms. For example, gcr($k{-}1$)/\discretionary{}{}{}gmres($k$)
is an homogeneous \ocm (k 1) method minimizing the $2$-norm of
the residual. Orthomin($m$) is an homogeneous \ocm (1 m{+}1)
method that also minimizes the $2$-norm of the residual and has
only the latest residual in the selection space. The conjugate
gradient and conjugate residual algorithms are homogeneous \ocm
(1 2) methods that minimize various norms and also have only the
latest residual in the selection space. It would be interesting to find
a polynomial-based iterative algorithm that is not an operator
coefficient method.

Methods that employ the entire $(k+1) \times m$ tableau may
greatly reduce the matrix-vector multiplications needed to solve
equations to prescribed accuracy. The convergence rate of
truncated orthomin generally improves as the order, $m$, increases.
Orthomin is a 1st degree method, $k = 1$, and similar be\-havior
may be expected for higher degree methods, $k > 1$. Figure~4 shows
convergence significantly improves by increasing $m$ and fixing
$k$. In this case the matrix-vector multiplications for each step are
independent of $m$ and are identical to those of
gcr($k{-}1$)/\discretionary{}{}{}gmres($k$). Thus, convergence
quickens by solving larger minimization problems but by performing
the same matrix-vector multiplications per step. Faster
convergence means fewer steps, and fewer matrix-vector
mul\-ti\-pli\-ca\-tions overall. This subject is discussed again in
Section~7.

\pageinsert {\Picture 4. 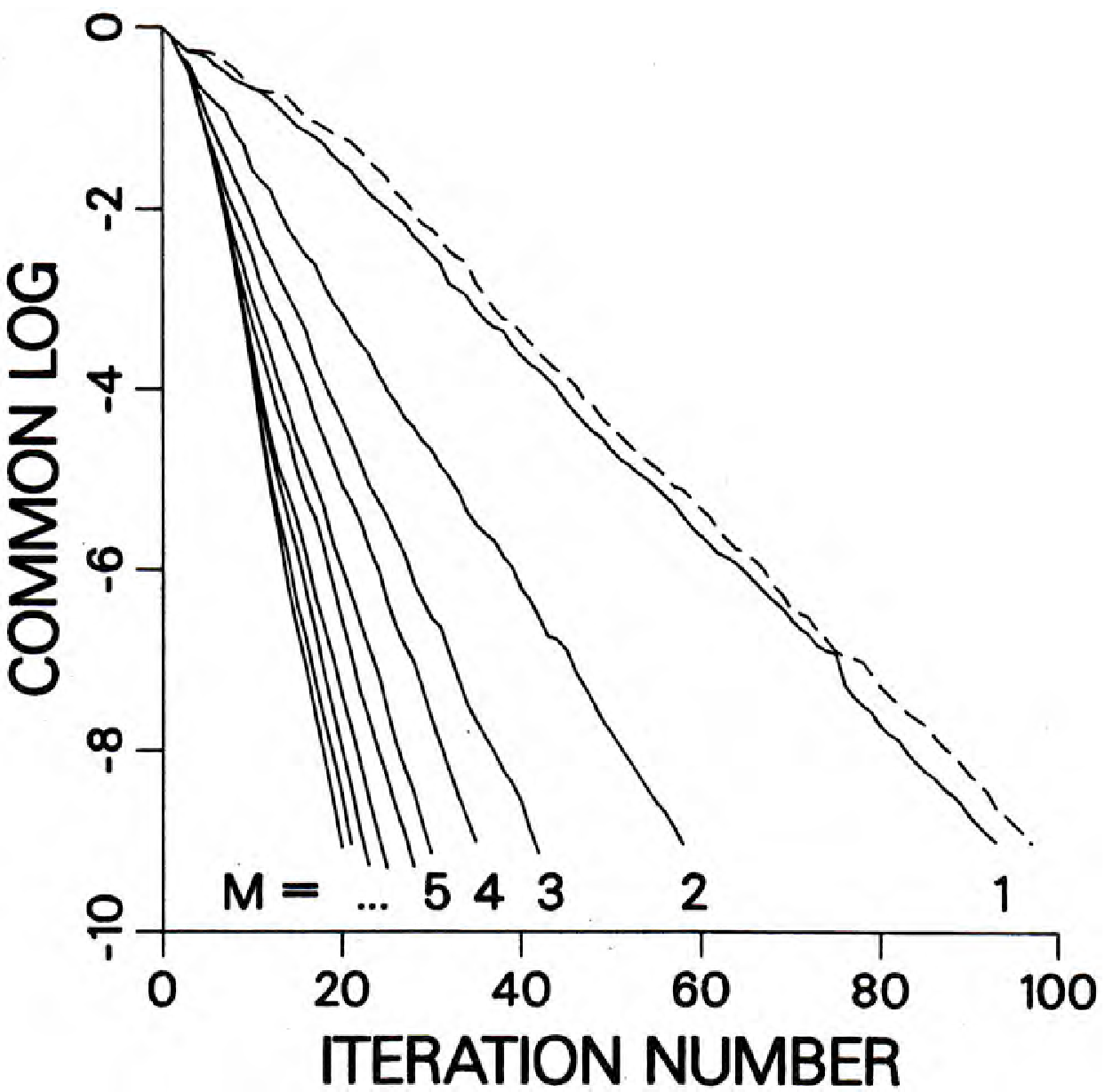 $2$-norm relative residuals for homogeneous
(dashed) gcr\/$(5)$/\discretionary{}{}{}gmres\/$(6)$ and
inhomogeneous (solid) \ocm (6 m), $m = 1$ $2$ $\ldots$ $10$, applied
to the same system. Appendix 2 and Section~5 explain the
calculations.\par }\endinsert

Reducing matrix-vector multiplications is a significant achievement
because they can account for most of the computational work. When
the matrix is randomly sparse, then matrix-vector multiplications
perform random memory accesses which are comparatively slow. If
the matrix is not explicitly known, as in matrix-free solution of
ordinary differential equations [\Brown ], then matrix-vector
multiplications require numerical differentiation of functions whose
evaluation may be very slow. 

Several recurrence formulas are associated with an operator
coefficient method. The selection criterion $$\hbox {minimize $\| r_n
\|_{\hbox {whatever}}$ or whatever}$$ chooses a  {\sl coefficient
tableau} $$\tableau {c_{(0,1)}& c_{(0,2)}& \;\; \ldots \;\; &
c_{(0,m)}\cr c_{(1,1)}& c_{(1,2)}& \;\; \ldots \;\;& c_{(1,m)}\cr
c_{(2,1)}& c_{(2,2)}& \;\; \ldots \;\;& c_{(2,m)}\cr \Vdots \;\;& \;
\Vdots \;\;& \; \Ddots \;\;& \Vdots \;\;\cr c_{(k,1)}& c_{(k,2)}& \;\;
\ldots \;\;& c_{(k,m)}\cr }$$ which produces the next iterate $$x_n =
\sum_{j=1}^m c_{(0,\,j)} x_{n-j} + \sum_{i=1}^k \sum_{j=1}^m
c_{(i,\,j)} A^{i-1} r_{n-j} .$$ Parenthetical subscripts in this and other
formulas indicate dependence on $n$. The next residual can be
obtained by a similar formula $$r_n = \sum_{j=1}^m c_{(0,\,j)} r_{n-j}
- \sum_{i=1}^k \sum_{j=1}^m c_{(i,\,j)} A^i r_{n-j} + f_n$$ $$f_n = y -
\sum_{j=1}^m c_{(0,\,j)} y$$ which can be written as a recurrence
formula of order $m$ $$r_n = P_{(1)} (A) r_{n-1} + P_{(2)}(A) r_{n-2}
+ \cdots + P_{(m)}(A) r_{n-m} + f_n$$ whose coefficients are
operators, that is, are polynomials of degree $k$ $$P_{(j)}(X) =
c_{(0,\,j)} - c_{(1,\,j)} X - c_{(2,\,j)} X^2 - \cdots - c_{(k,\,j)} X^k$$
evaluated at the matrix $A$. Whence the name, {\sl operator
coefficient method of degree $k$ and order $m$}. If $c_{(0,1)} +
c_{(0,2)} + \cdots + c_{(0,m)} = 1$, then the $f_n$'s vanish and the
re\-sid\-u\-als satisfy homogeneous recurrence formulas. Whence
{\sl homogeneous} and {\sl in\-ho\-mo\-gen\-eous} methods. In the
homogeneous case, the residuals can be expressed succinctly in
terms of the initial residuals $$r_n = P_{n,1} (A) r_{0} + P_{n,2} (A)
r_{-1} + P_{n,3} (A) r_{-2} + \cdots + P_{n,m} (A) r_{1-m}$$ by means
of residual polynomials, $P_{n,\,j} (X)$, generated from the
recurrence formulas $$P_{n,\,j} = P_{(1)} P_{n-1,\,j} + P_{(2)}
P_{n-2,\,j} + \cdots + P_{(m)} P_{n-m,\,j}$$ with initial values
$P_{1-j,\,j} = 1$ and others zero. This representation for $r_n$ is
numerically correct, however, only if the recurrence formulas are
stable when $X=A$.

\beginsection {6. Implementations}

The implementation of inhomogeneous operator coefficient
methods that minimize the $2$-norm of the residual is considered
here. This section has two parts. The first analyzes
implementations of older methods, the second describes a
reasonable implementation for all \ocm(k m) methods. Those who
are interested in the new methods may prefer to read the notation
below and to begin at Section~6b.

Inhomogeneous \ocm (k m) methods that minimize the $2$-norm of
the residual perform the following task at step $n$. They choose
$x_n = V_n c_n$, where $c_n$ solves the least squares problem
$$\min \| y - A V_n c_n \|_2,$$ and where $V_n$'s columns are a basis
for the selection space. The {\sl natural basis} for the full-tableau
method is the following. $$\Span {V_n} \; = \; \hbox {span} \, \left\{
\tableau {x_{n-1}& x_{n-2}& \;\; \cdots \;\;& x_{n-m}\cr r_{n-1}&
r_{n-2}& \;\; \cdots \;\;& r_{n-m}\cr \;\; \Vdots \;\;&  \;\; \Vdots
\;\;& \; \Ddots \;\;& \;\; \Vdots \;\;\cr A^{k-1} r_{n-1}& A^{k-1}
r_{n-2}& \;\; \cdots \;\;& A^{k-1} r_{n-m}\cr } \; \right\}$$ 

Any implementation makes three choices. The first is the basis for
the selection space. This basis becomes the columns of $V_n$. The
second choice is the basis for the least squares problem. This might
be the columns of $A V_n$. The third choice is the process to solve
the least squares problem. All the choices affect both numerical
accuracy and computational efficiency. The bases might overlap to
conserve storage, for example, or they might facilitate the solution
process to conserve time. The chief numerical considerations are
the accuracy of the bases and the accuracy of the least squares
solution. A comparative analysis of all the possibilities is beyond
the scope of this paper. Ashby, Manteuffel and Saylor [\Ashby ],
Saad and Schultz [\Saadb ], Walker~[\Walkera ]~[\Walkerb ] and
references therein should be consulted for more implementation
ideas.

\beginsubsection {6a. Some Existing Implementations}

This section analyzes implementations of known methods. It
reverses Section~3's simplification process, and rebuilds
the algorithms with explicit justification for each implementation
detail.

Like more general operator coefficient methods, gcr($k{-}1$)
and gmres($k$) solve a least squares problem, but their selection
space is smaller and affine. They choose $x_n = x_{n-1} + V_n c_n$,
where $c_n$ solves $$\min \| r_{n-1} - A V_n c_n \|_2,$$ and where
$$\Span {V_n} = \hbox {span} \, \left\{ \; \tableau {r_{n-1}\cr
Ar_{n-1}\cr \Vdots \;\;\cr A^{k-1} r_{n-1}\cr } \; \right\} .$$ 

The gcr implementation of gcr($k{-}1$)/\discretionary{}{}{}gmres($k$) maintains two
separate bases. It uses $A V_n$ for the least squares basis, and it
chooses $V_n$ to be $\transpose A A$-orthogonal. Evaluation of the
inner products during the orthogonalization process requires either
that the columns of $A V_n$ be saved, or that additional
matrix-vector multiplications be performed. Since $A V_n$ is
Euclidean-orthogonal, the normal equations are diagonal and are
easily solved. However, normal equations may solve least squares
problems with accuracy less than best.
 
The gmres version of gcr($k{-}1$)/\discretionary{}{}{}gmres($k$) may be the most
efficient for this method. It chooses an Euclidean-orthonormal basis
for $\Span {V_n}$ that becomes an orthonormal basis for $\,\Span
{V_n} + \,\Span {A V_n}$ by the inclusion of one more vector. That is,
$A V_n = W_n H_n$ where $W_n$ has the columns of $V_n$ plus one,
and where the small matrix $H_n$ is constructed along with the
basis. $H_n$ represents $A$ under an orthonormal change of basis. If
the change of basis can be computed accurately, then $H_n$ can be
used to solve the least squares problem accurately.

It is not clear whether so efficient an implementation is possible
for more general operator coefficient methods. The \ocm (k 1)
methods have the advantage of simplicity because their selection
spaces involve a single group of nested Krylov spaces. $$\Span
{r_{n-1}} = K_{0} \subseteq K_{1} \subseteq \cdots \subseteq
K_{k-1} = \Span {r_{n-1} \;\; Ar_{n-1} \;\; \ldots \;\; A^{k-1}r_{n-1}}$$
More general \ocm (k m) methods have several Krylov spaces and
so may not attain the efficiencies of the 1st order, $m=1$, methods. 

The usual practice with \ocm (k 1) methods is to recursively build 
orthogonal bases for the nested Krylov spaces, $$K_{j+1} = \Span
{p_0 \;\; p_1 \;\; \ldots \;\; p_{j} \;\; p_{j+1}} = \Span {K_j \;\;
p_{j+1}},$$ in which $p_{j+1}$ is orthogonal to $K_{j}$. Restarted gcr
employs $\transpose A A$ orthogonality, while restarted gmres
chooses Euclidean orthogonality. Orthogonality is desired for two
reasons. It is generally believed orthogonal bases provide better
numerical representations for their spans, moreover, orthogonality
can help solve the least squares problem. 

The experience with orthomin($m{-}1$), an homogeneous \ocm (1 m)
method, suggests orthogonality has a third use. It may provide
efficient implementations of high order, $m > 1$, methods. Like
gcr/gmres, orthomin($m{-}1$) chooses $x_n = x_{n-1} + V_n c_n$
where $c_n$ solves $$\min \| r_{n-1} - A V_n c_n \|_2,$$ but in this
case $$\Span {V_n} = \hbox {span} \, \left\{ \; \tableau {r_{n-1}&
(x_{n-1} - x_{n-2})& \ldots& (x_{n-(m-1)} - x_{n-m})\cr} \;
\right\}.$$ Orthomin($m{-}1$) represents this selection space by an
{\sl inventory} of $\transpose A A$-or\-thog\-o\-nal basis vectors,
$$V_n = \left[ \matrix {p_{n-1} \;\; p_{n-2} \;\; \ldots \;\; p_{n-m}}
\right].$$ Each step maintains the basis by discarding the oldest
vector and inserting the residual's component orthogonal to the
others.  As with gcr, $\transpose A A$-orthogonality results in
diagonal normal equations for the least squares problem. The
solution update involves only the newest basis vector because
previous steps account for the others. In this way, each basis vector
spans exactly the difference between a pair of successive
approximate solutions.

Some algorithms of Chronopoulos and Gear [\Chronopoulosa ]
[\Chronopoulosb ] [\Chronopoulosc ] [\Chronopoulosd ] appear to be
more general \ocm (k m) methods that follow the approach taken by
truncated orthomin. They build the natural basis for a Krylov space
of low dimension, say $k$, and then perform an orthogonalization
step to enforce $\transpose A A$-orthogonality among a number of
such spaces, say $m$. Like orthomin, the difference between a pair
of successive approximate solutions lies in a space of low
dimension, in this case $k$, but in the absence of arithmetic error
the new approximate solution is the best in a larger space. An
analysis like the one in Section~3 for orthomin would be needed to
identify the selection spaces in terms of the natural \ocm (k m)
basis.

All implementations that rely on recursively produced, orthogonal
bases can be expected to share the failing of the original conjugate
gradient algorithm. The vectors are not orthogonal in practice and,
as shown by Figure~1 for the conjugate gradient method, they can
be quite different from the intended vectors. The loss of
orthogonality in the basis is readily detected, and obviously affects
the accuracy of the least squares solutions. The loss of accuracy in
the basis vectors is difficult to detect, but surely affects the
essential character of the approximations. Elaborate means such as
reorthogonalization can remedy the orthogonality, but aside from
producing more nearly orthogonal vectors, they have not been
proved to result in better approximations to the underlying Krylov
spaces. It is an open question whether the approximations can be
made consistently better. The natural bases have been observed to
be badly conditioned [\Chronopoulosc ] [\Walkera ], so linear
transformations that make them better conditioned evidently must
be ill-conditioned too, and thus must be difficult to apply
accurately.

The least squares problem can be difficult to solve however it is
formulated. The gcr/orthomin approach may be flawed because it
solves the normal equations, but other methods must contend with
near singularity of the least squares bases. The matrix $H_n$ of the
gmres approach has a $2$-norm condition number no worse than
$A$'s, but by being smaller it may reflect ill-conditioning more.
Alternatively, the natural bases of Krylov spaces can be very
nearly singular. 

It is difficult to concede that any but the best solution method
should be applied to the least squares problem. For least squares
problems in general, ``the only fully reliable way to treat rank
deficiency is to compute the singular value decomposition''
[\Golube ,~p.~170]. Apparently no iterative methods heed this
advice. Yet the singular value decomposition is fairly inexpensive
for the small matrices that appear in restarted gmres, for example,
and may remove some of the difficulties occasionally reported for
this method [\Walkera ] [\Walkerb ].

In summary, existing implementations always employ orthogonal
bases to reduce storage and computation. The savings in storage
appear to be at most a factor of two, as for gmres versus gcr. This
improvement is marginal on present-day computers and should not
govern the choice of implementations. The savings in time may be
more significant. 

The advantages of orthogonal bases must be weighed against 
numerical concerns. The $\transpose A A$-orthogonal bases impose
inferior least squares solution methods. In the presence of rounding
error moreover, it is known that recursively generated bases may
not accurately span the intended spaces. The effects of this on
nonterminating, iterative algorithms are largely unexplored.

\beginsubsection {6b. An Implementation}

The following implementation is generic to all \ocm (k m) methods.
It solves the least squares problems by the singular value
decomposition, the best available method, and uses the natural
bases for the Krylov spaces. This implementation is offered both as
a research tool and as a model of programming simplicity. It has
several advantages. First, the implementation allows easy
substitution of methods, including restarted
gcr/\discretionary{}{}{}gmres and truncated orthomin. Second, it
addresses some numerical difficulties likely to trouble both old and
new methods. Third, it conveniently relies on well-known numerical
procedures found in many scientific computing libraries. This
implementation, applicable to all \ocm (k m) methods and devoid of
programming complications, may be the most appropriate in the
present, early stages of development. 

With the natural basis, all algorithms have the same implementation
but for the choice of basis vectors. Methods such as conjugate
residual and orthomin that don't use the full tableau can be
implemented by simply choosing a subset of the larger basis. The
columns of $A V_n$ for the specific $V_n$ of interest must be
constructed explicitly. Most can be borrowed from previous steps.
Only the vectors $A r_{n-1}$, $A^2 r_{n-1}$, $\ldots\,$, $A^k r_{n-1}$
associated with the most recent solution are new. They require $k$
matrix-vector products. The vector $A x_{n-1}$, which also forms
$r_{n-1}$, requires one more matrix-vector product or can be
obtained recursively. 

The least squares solution process is numerically robust. The
singular value decomposition solves the least squares problem
more accurately, though perhaps more expensively, than
orthomin-like implementations would solve the normal equations.
Errors can enter the least squares basis only through the
matrix-vector multiplications which produce $A V_n$ from $V_n$.

Operator coefficient methods should alleviate the concern that the
natural bases for Krylov spaces are too nearly singular. High order
operator coefficient methods make high degree Krylov spaces
unnecessary. If very high degrees are needed, then
Euclidean-orthogonal bases may be computed in the manner of
Arnoldi, and may be integrated into the computations.

The following steps compute the minimum norm solution of the
least squares problem. They are based on recommendations in the
text of Golub and Van~Loan [\Golube ]. First, the columns of $A V_n$
should be scaled to have unit $2$-norms. This makes $\kappa_2 (A
V_n)$ nearly minimal and improves numerical accuracy. Scaling also
avoids numerical overflow and underflow when a matrix repeatedly
multiplies a vector. Second, Householder transformations should
reduce $A V_n$ to an upper triangular matrix. This reduces the
arithmetic costs when, as here, there are many more rows than
columns. Third, the singular value decomposition of the small, upper
triangular matrix must be computed. Fourth, the minimum norm
solution of the column-scaled least squares problem can be
approximated by applying an approximate pseudoinverse obtained
by ignoring small singular values. Singular values smaller than
machine round-off relative to the largest singular value might be
discarded, or more sophisticated methods might be used to
determine numerical rank. Finally, the unscaled minimum norm
solution, $c_n$, combines the columns of $V_n$ to  produce the next
approximate solution for the \ocm (k m) method, $x_n = V_n c_n$.
Table~3 restates these steps and counts their arithmetic
operations.

\midinsert {\table 3. {An implementation of \ocm (k m) methods
minimizing the $2$-norm of the residual for a selection space basis
of size $t \le (k + 1) m$, with operation counts. Terms
independent of the matrix order $N$ are omitted. Section~6b
provides further explanation.} {\vbox {\openup 0.5\jot \halign
{$\;\;$#\hfil& \hfil #$\;\;$\cr \omit \hfil step\hfil& operations\cr
\noalign {\medskip \hrule \medskip} scale $A V_n$ to unit column
norm,& $3tN$\cr \qquad $V := A V_n D^{-1}$\cr Householder
transformations reduce& $2t^2N$\cr \qquad $V$ to triangular form,
$R := Q V$ \cr orthogonal projection of $y$, $z := Q y$& $4tN$\cr
singular value decomposition of $R$\cr \qquad solves $\min \, \| R D
c_n - z \|_2$\cr assemble next solution, $x_n := V_n c_n$& $(2t -
1)N$\cr \noalign {\medskip \hrule \medskip} optionally assemble
$Ax_n$ and $r_n$& $2tN$\cr \noalign {\medskip \hrule }}}}\par
}\endinsert

With $t$ basis vectors, $t \le (k+1) m$, the implementation needs
the following resources per step of the inhomogeneous
\ocm (k m) method. There are $k$ matrix-vector multiplications,
and one more if residuals are not recursively computed. Memory
space of $2 t$ vectors is needed to store $V_n$ and the Householder
transformations that reduce $A V_n$ to triangular form. Very
compact memory management schemes are possible since the basis
vectors pass from one step to the next but the Householder
transformations do not. Table~3 counts $(2t^2 + 9t - 1) N$
arithmetic operations, from which terms independent of the matrix
order, $N$, have been discarded. The residual calculation requires
either $N$ or $2tN$ operations.

This implementation repeatedly solves a large, dense,
overdetermined, singular, least squares problem. This task is basic
to numerical linear algebra. The Householder reduction and the
singular value decomposition already appear in many computing
libraries, and solution methods tuned to specialized computer
architectures are being developed. The implementation therefore
improves, in a sense automatically, with advances to numerical
software and hardware.

\beginsection {7. Varying Coefficients}

Many operator coefficient methods dynamically select
recurrence coefficients by minimizing the $2$-norm of the
residual. This selection criterion is examined here. Theorem~1
and its Corollary prove convergence for a large class of
matrices distinguished by a simple polynomial relationship.
Moreover, experimental results indicate there is a trade-off
between degree and order. This may allow high order methods to
replace comparatively less economical high degree methods.

\Proclaim Theorem 1. If the Hermitian part $H$ of $P(A)$ is
positive or negative definite for some polynomial $P$ with
degree at most $k$ and $P(0) = 0$, then for every $x_n$ the
affine space $$x_n + \Span {r_n \;\; A r_n \;\; A^2 r_n \;\; \ldots
\;\; A^{k-1} r_n}$$ contains a vector $x_{n+1}$ with $\| r_{n+1} \|_2
\le \rho \, \| r_n \|_2$ where $$\rho = \sqrt {1 - \left[ \, {\min |
\lambda (H) | \over \| P(A) \|_2} \, \right]^2} \; < \; 1.$$ The
affine space also contains a vector $x_{n+1}$, usually
different from the first, with $\| e_{n+1} \|_2 \le \rho \, \| e_n
\|_2$ (proof appears in Appendix~1).

\Proclaim Corollary to Theorem 1. If the Hermitian part $H$ of
$P(A)$ is positive or negative definite for some polynomial $P$
with degree at most $k$ and $P(0) = 0$, then \ocm (k m) methods
whose selection spaces contain the affine space $$x_n + \Span
{r_n \;\; A r_n \;\; A^2 r_n \;\; \ldots \;\; A^{k-1} r_n}$$ converge
for the selection criteria that minimize the $2$-norm of either
the residual or the error. At each step the norm declines by at
least the factor $$\sqrt {1 - \left[ \, {\min | \lambda (H) | \over \|
P(A) \|_2} \, \right]^2} \; < \; 1.$$

The thesis of Elman [\Elmana ] is the inspiration for Theorem~1.
The Corollary applies to
gcr($k{-}1$)/\discretionary{}{}{}gmres($k$) as well as to more
general methods. However, only the case in which $A$ itself is
positive or negative definite appears to have been published
previously, by Eisenstat, Elman and Schultz [\Eisenstatb ]. Saad
and Schultz mention this case too, and present more detailed
convergence results for diagonalizable matrices~[\Saadc ].

Theorem~1's bound on the convergence rate may be weak
because it is independent of the recurrence order $m$. The
Theorem minimally assumes each residual $r_n$ equals a linear
combination that includes vectors from $$\Span {r_{n-1} \;\; A
r_{n-1} \;\; A^2 r_{n-1} \;\; \ldots \;\; A^k r_{n-1}},$$ but 
the combination also may employ vectors from the larger
space $$\hbox {span} \; \left\{ \tableau {r_{n-1}& r_{n-2}& \;\;
\cdots \;\;& r_{n-m}\cr Ar_{n-1}& Ar_{n-2}& \;\; \cdots \;\;& 
Ar_{n-m}\cr \;\; \Vdots \;\;&  \;\; \Vdots \;\;& \; \Ddots \;\;& \;\;
\Vdots \;\;\cr A^k r_{n-1}& A^k r_{n-2}& \;\; \cdots \;\;& A^k
r_{n-m}\cr } \; \right\} \, .$$ Thus, $r_n$ depends on powers of $A$
up to $A^{km}$. This suggests the convergence rate may vary
with the product $km$. 

Figure~5 provides numerical evidence for this interpretation.
The Figure exhibits level curves of observed convergence rates
as functions of $k$ and $m$ for the convergence histories 
shown in Figure~6. The level curves have the expected
qualitative behavior. In this example, the \ocm (6 1) and \ocm (3
5) methods have essentially the same convergence rate. This
means they achieve the same accuracy in the same number of
steps, but the \ocm (3 5) method requires half the
matrix-vector multiplications. 

\pageinsert {\Picture 5. 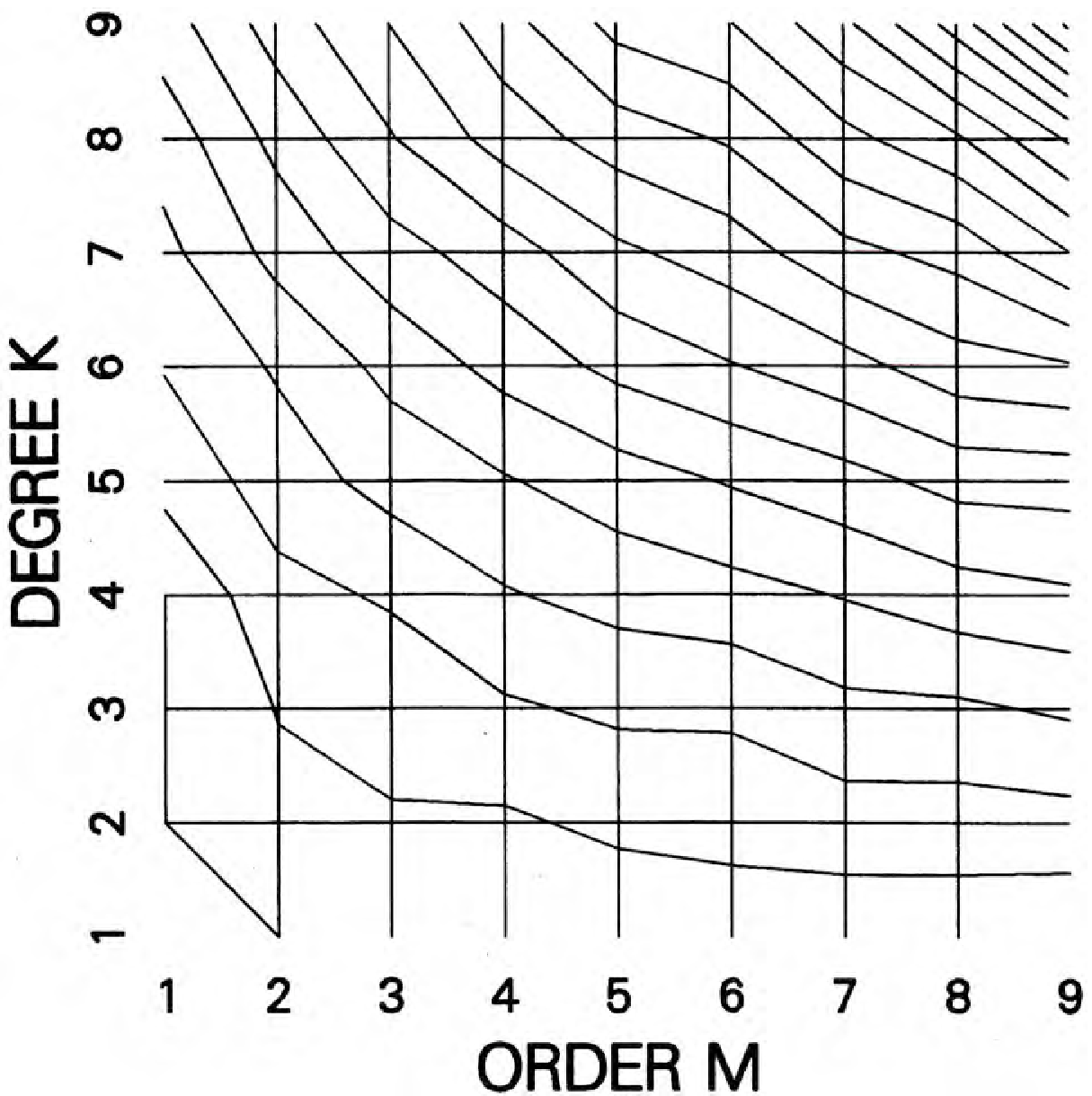 Level curves of observed convergence
rate for inhomogeneous \ocm (k m) as a function of $k$ and $m$
for one system. The curves range in multiplicative steps of
$10^{0.05}$ from $10^{-1}$ at the upper right to $10^0$, no
convergence, at the lower left. The left edge corresponds to
gcr\/$(k{-}1)$/\discretionary{}{}{}gmres\/$(k)$, the bottom edge
to orthomin\/$(m{-}1)$. Figure~6 displays the convergence
histories. Appendix~2 and Section~7 explain the calculations.
\par }\endinsert

\pageinsert {\Picture 6. 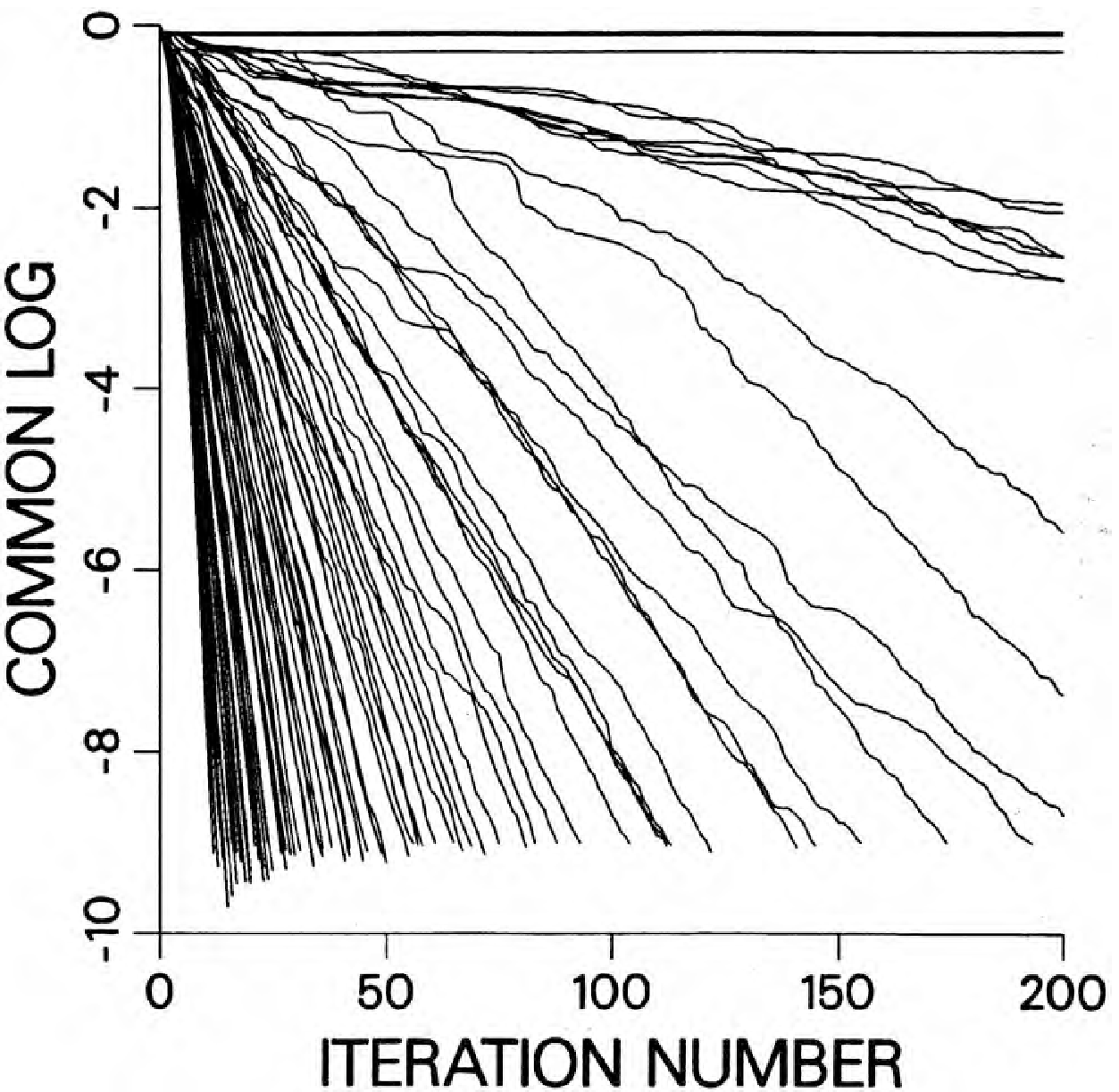 Convergence histories from which the
level curves of Figure~5 are derived. Appendix~2 and Section~7
explain the calculations. \par} \endinsert

This apparent trade-off between $k$ and $m$ may be the most
important aspect of \ocm (k m) methods. It allows the beneficial
effects of larger $k$ in gcr($k{-}1$)/\discretionary{}{}{}gmres($k$) to be realized
more economically with larger $m$. The most efficient choice of
$k$ and $m$ can be expected to change with the computer and
the problem, and with the expense of solving each step's
minimization problem relative to the expense of  performing
matrix-vector multiplications. If the convergence rate $\rho$
did vary only with $km$, for example if $$\rho (k, m) = \rho (
k/\ell, \, \ell \, m ) ,$$ then low degree, high order methods would
be more economical. Minimizing the $2$-norm of the residual
involves solving a least squares problem with a basis of size
$(k+1)m$ whose computation and memory requirements are
roughly constant among \ocm (k m) methods with the same $km$.
But \ocm (k m) performs $k$ or $k+1$ matrix-vector
multiplications per step, and this part of the total cost
decreases with $k$. Considerations of this kind can be expected
for all manner of coefficient choices.

The frequently used minimization criteria of the kind in
Theorem~1 are a powerful but imprecise tool for choosing
coefficients. They may not find coefficients that produce
convergence, and even when they do, they may not produce the
fastest convergence. Moreover, the Theorem's bound may be a
poor estimate for the convergence rate because the selection
criteria minimize norms of vectors, but the bound employs
norms of matrices. When these matrix-vector norms are applied
to matrices, they depend on both the eigenvalues and
eigenvectors, while Theorem~2 in the next section shows
convergence can depend on the eigenvalues alone. 

The following
example illustrates these concerns. The matrix $$\left[ \,
\matrix {\alpha& \beta\cr &\alpha& \beta\cr &&\alpha& \ddots
\cr &&&\ddots \cr} \; \right]$$ has eigenvalue $\alpha$ and its
Hermitian part has eigenvalues between $\hbox {Real} \, (\alpha)
\pm | \beta |$. Among \ocm (1 m) methods, a stationary
Richardson's 1st order method can be made to converge
independent of $\beta$, yet $\beta$ can be chosen so the
Hermitian part is indefinite and Theorem~1's bound is ineffective.

\beginsection {8. Constant Coefficients}

This final section determines exactly when operator coefficient
methods with constant coefficients converge. This information
has several uses. First, it may guide the choice of $k$ and $m$
needed to achieve convergence with coefficients selected by
any means. Second, it suggests ways to select coefficients other
than by the usual minimization criteria of Section~7. Finally, it
proves that some operator coefficient methods are new by
showing they converge when previously known methods do not. 

The Chebyshev iteration and the stationary 2nd order method are
\ocm (1 2) methods that converge only for matrices with
eigenvalues inside ellipses that exclude $0$. It is demonstrated
below that some constant coefficient \ocm (1 2) methods
converge for non-elliptical eigenvalue distributions. These, then,
are new methods.

Only the homogeneous case is possible for constant
coefficients. As remarked in Section~4, with the approximate
solutions converging to $x_*$, and with the residuals converging
to $0$, the sums of $x$-coefficients on both sides of the
recurrence equation must balance.

 \goodbreak
 \Proclaim Theorem 2. A constant coefficient, homogeneous,
operator coefficient method of degree $k$ and order $m$ $$x_n =
\sum_{j=1}^m c_{0,j} x_{n-j} + \sum_{i=1}^k \sum_{j=1}^m c_{i,j}
A^{i-1} r_{n-j}$$ with coefficient tableau $$\matrix {c_{0,1}&
c_{0,2}& \cdots& c_{0,m}\cr c_{1,1}& c_{1,2}& \cdots& c_{1,m}\cr
\vdots& \vdots& \ddots& \vdots\cr c_{k,1}& c_{k,2}& \cdots&
c_{k,m}\cr}$$ converges to a solution of $A x = y$ for all $y$ and
all initial vectors $x_{0}$, $x_{-1}$, $\ldots\,$, $x_{1-m}$ exactly
when, for each eigenvalue $\lambda$ of $A$, the maximum
magnitude $r(\lambda )$ of the roots $X$ of the polynomial
$$P(\lambda, X) = X^m - P_1(\lambda ) X^{m-1} - P_2(\lambda )
X^{m-2} - \cdots - P_m(\lambda ) X^{m-m}$$ with coefficients
given by the columns of the tableau $$P_j (\lambda ) = c_{0,\,j} -
c_{1,\,j} \lambda - c_{2,\,j} \lambda^2 - \cdots - c_{k,\,j}
\lambda^k$$ is strictly less than $1$. Moreover, there is a bound
upon the residuals for all $y$ and all initial vectors $x_{0}$,
$x_{-1}$, $\ldots\,$, $x_{1-m}$ $$\| r_n \| \le \left( \, \| r_{0} \| + \|
r_{-1} \| + \cdots + \| r_{1-m} \| \, \right) Q(n) \, R^n,$$ and if $A$
is nonsingular there is an identical bound upon the errors $$\| e_n
\| \le \left( \, \| e_{0} \| + \| e_{-1} \| + \cdots + \| e_{1-m} \| \,
\right) Q(n) \, R^n.$$ $R$ is the maximum $r(\lambda )$ for all the
eigenvalues of $A$. $Q(n)$ is a polynomial that depends on the
norm, on $A$, and on the coefficient tableau. The norm may be any
consistent matrix-vector norm (proof appears in Appendix~1).

There is some evidence that constant coefficients may work well
in the long run. The Chebyshev iteration's coefficients converge
to values for which the stationary 2nd order method converges
identically [\Manteuffeld ]. To the extent coefficients chosen by
some means do become stationary, Theorem~2 explains the
minimal $k$ and $m$ necessary before the coefficient selection
criteria can make \ocm (k m) methods converge. An entirely
constant coefficient iteration might be useful when many
systems of equations feature the same matrix. The trick is to
find the coefficients, and Theorem~2 is the first step in this
direction.

\pageinsert {\Picture 7. 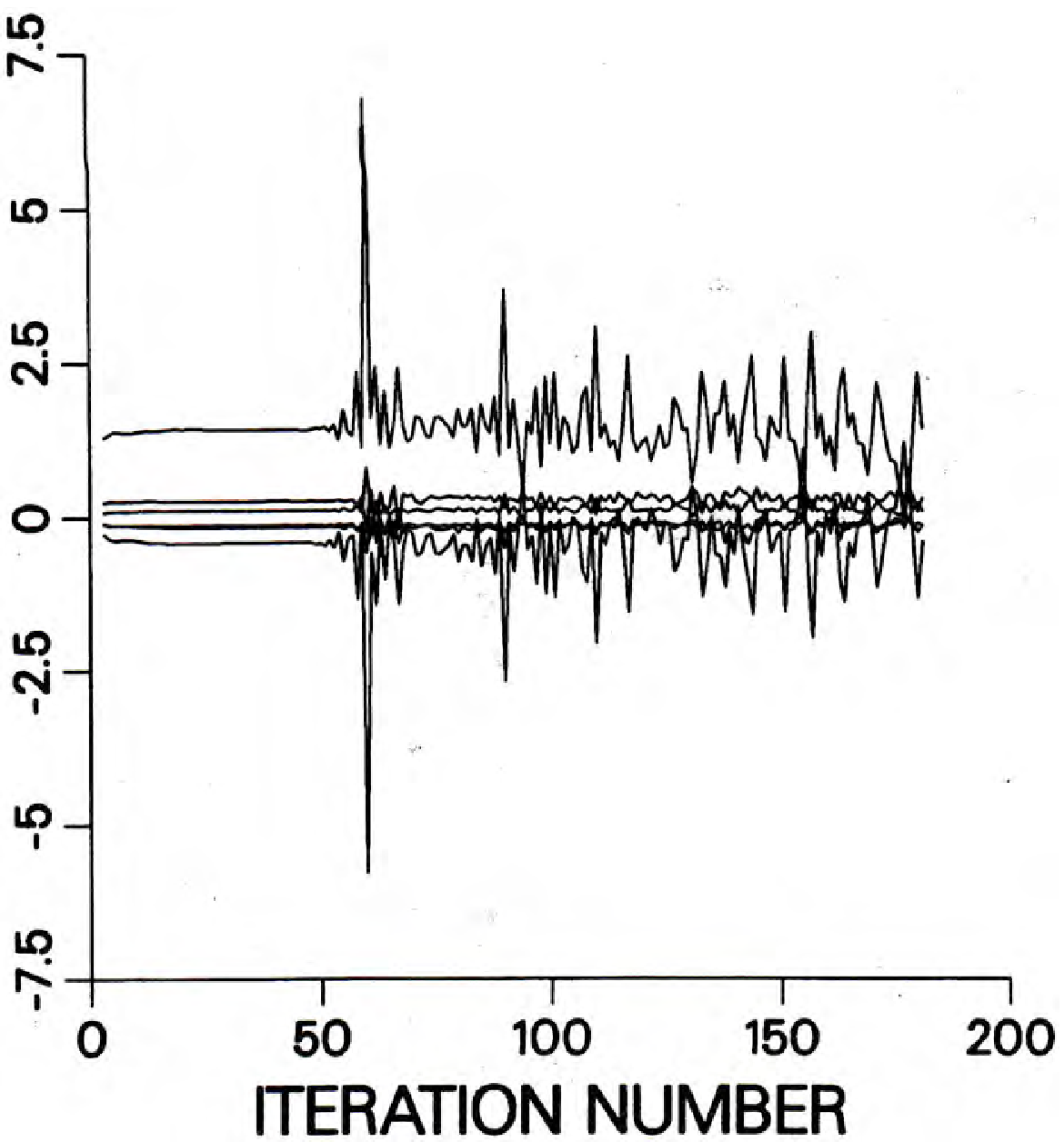 Coefficients of inhomogeneous \ocm (2 2)
minimizing the $2$-norm of the residual for one system.
Appendix~2 and Section~8 explain the calculations. \par}
\endinsert

\pageinsert {\Picture 8. 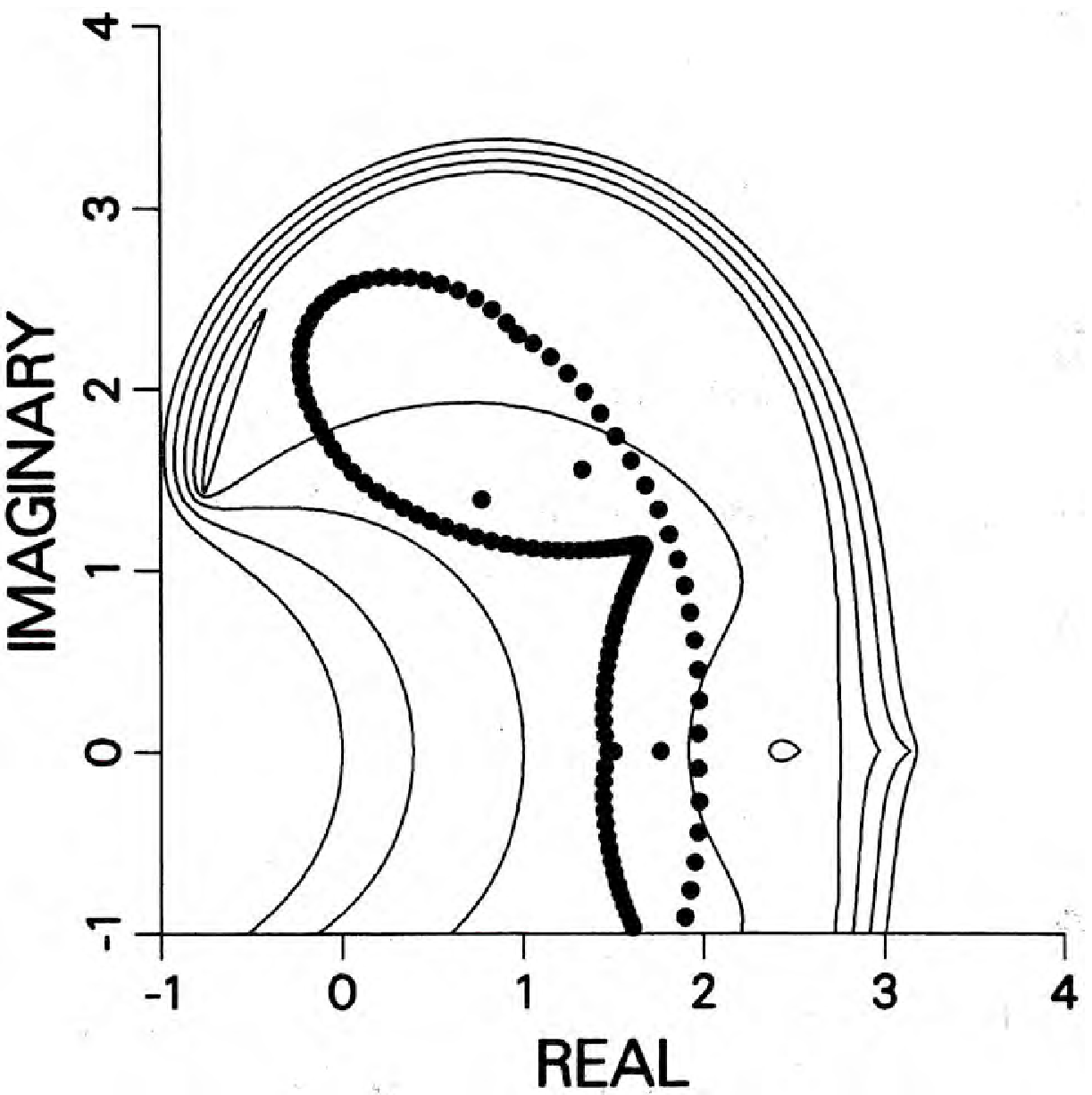 Eigenvalues of the matrix (solid dots)
superimposed on some level curves of the convergence rate
$r(\lambda)$ for the constant coefficient regime of Figure~7. 
The levels start at $1$ on the boundary and decrease in steps of
$0.05$. Appendix~2 and Section~8 explain the calculations. \par}
\endinsert

The following example suggests how constant coefficients might
be found, and clarifies the statement of Theorem~2. Figure~7
shows that coefficients chosen to minimize the residual's
$2$-norm can be nearly constant over several iterations. When
dynamically chosen coefficients remain fixed for a time, then a
constant coefficient iteration with these fixed values may
converge. To see if this is the case here, Figure~8 superimposes
the matrix eigenvalues, as black dots, over some level curves of
the Theorem's {\sl eigenvalue-specific convergence rate}
$$r(\lambda ) = \hbox {maximum $|X|$ of all $X$ for which
$P(\lambda, X) = 0$},$$ where $P(\lambda, X)$ is the Theorem's
polynomial $$P(\lambda, X) = X^m - \sum_{j=1}^m c_{0,\,j}
X^{m-j} + \sum_{i=1}^k \sum_{j=1}^m c_{i,\,j} \lambda^i X^{m-j}$$
and where the polynomial's coefficients $c_{\/i,\,j}$ are taken
from the nearly constant regime of Figure~7. The {\sl
convergence domain} is the set of $\lambda$ for which
$r(\lambda) < 1$. Figure~8 shows all the eigenvalues lie within
the convergence domain, so Theorem~2 guarantees convergence
with these constant coefficients. Figure~9 exhibits the relative
residuals for the original right hand side and one other in a
constant coefficient iteration. 

\pageinsert {\Picture 9. 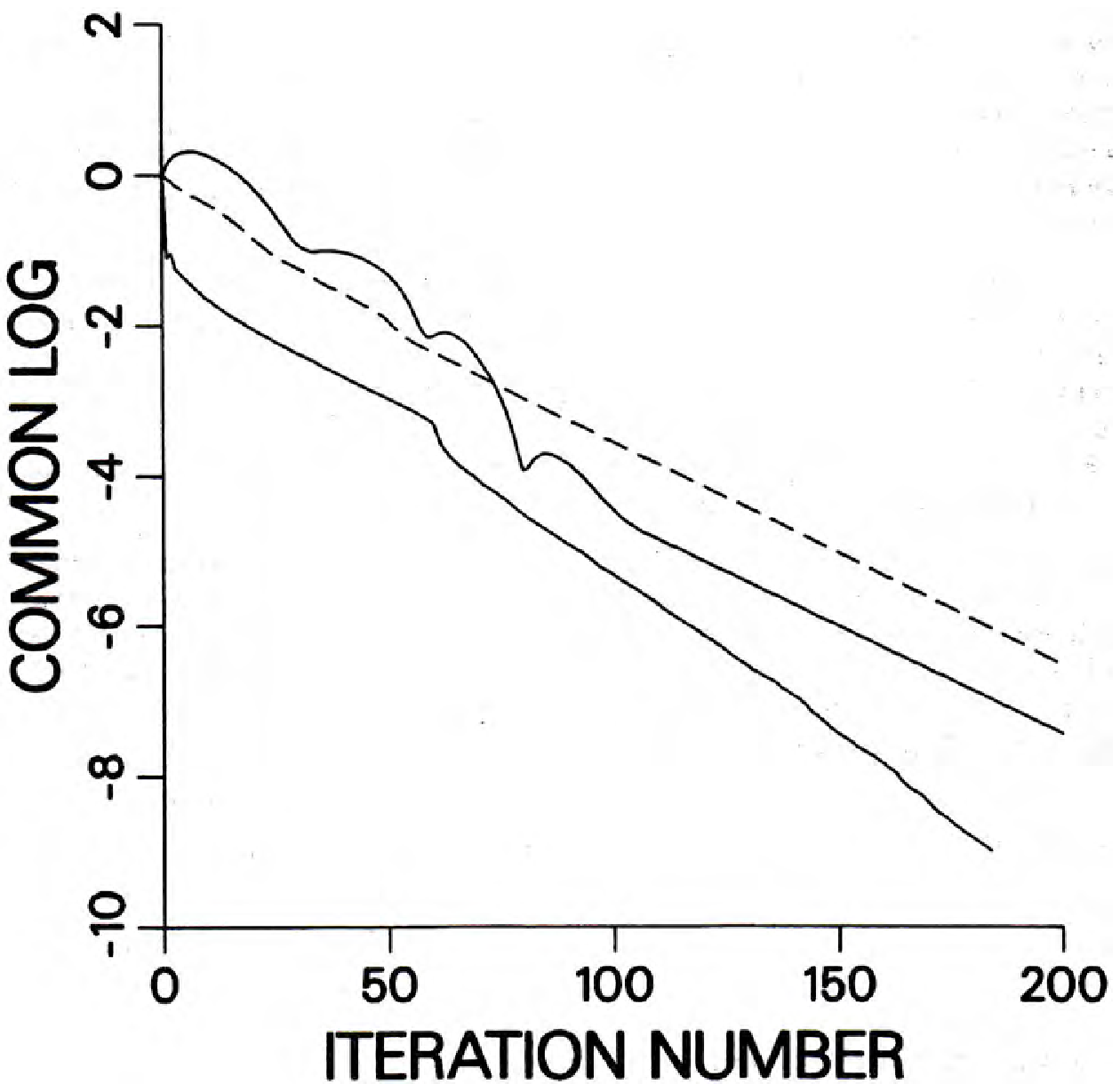 $2$-norm relative residuals for the
iteration of Figure~7 (lower solid line), and for a constant
coefficient iteration with the same right hand side (higher solid
line) and a different right hand side (dashed line). Appendix~2 and
Section~8 explain the calculations. \par} \endinsert

A phenomenon discovered by Trefethen [\Trefethen ] explains
why Figure~9 exhibits slower convergence than Figure~8
predicts. When the matrix eigenvalues are sensitive to
perturbation, then convergence depends on an envelope of
approximating eigenvalues introduced by rounding error.
Theorem~2 must be applied to these approximations to predict
the convergence rate. Nevertheless, convergence is assured in
Figure~9 because a convergent iteration that already accounts
for the envelope suggests the constant coefficients.

For a given matrix even with known eigenvalues or
approximations thereto, it can be  difficult to find any
convergent coefficients let alone optimal ones that minimize
Theorem~2's {\sl convergence rate}, $R = \max \, r(\lambda)$. The
inverse problem of finding eigenvalue domains convergent for
given coefficients is at least numerically straightforward. It
amounts to seeking the $\lambda$ for which all the roots $X$ of
the Theorem's polynomial $P(\lambda, X)$ have magnitude less
than $1$.

\goodbreak
Figure~10 displays the convergence domains for some arbitrary
coefficient tab\-leaux as large as $3 \times 2$, that is, for
methods up to \ocm (2 2). $$\vcenter {\normalbaselines
\mathsurround=0pt \ialign {\hfil$#$\hfil& $\quad #$\hfil\cr
\mathstrut\cr \noalign {\kern-\baselineskip} c_{0,1}& c_{0,2} = 1
- c_{0,1}\cr c_{1,1}& c_{1,2}\cr c_{2,1}& c_{2,2}\cr \mathstrut\cr
\noalign {\kern-\baselineskip}}}$$ Table~4 lists the specific
coefficients. 

\pageinsert {\Picture 10. 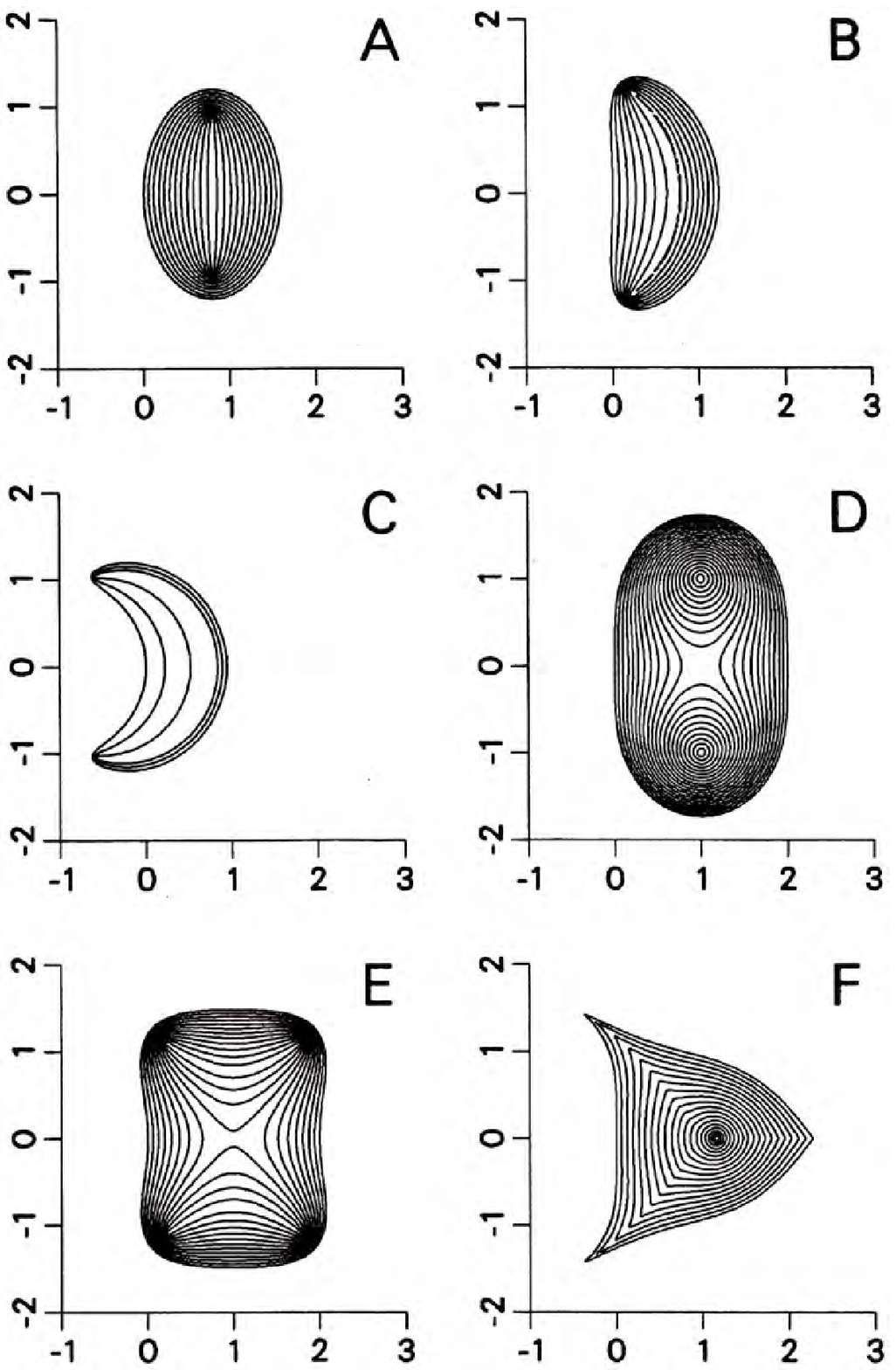 Level curves of convergence rate
$r(\lambda )$ as a function of $\lambda$ in the complex plane for
Table~4's constant coefficients. The levels begin at $1$ on the
boundaries and decrease in steps of $0.05$. Appendix~2 and
Section~8 explain the calculations. \par} \endinsert

\midinsert {\table 4. {Constant coefficient tableaux for the
convergence domains pictured in Figure~10. Section~8 provides
further explanation.} {\vbox {\normalbaselines \halign {\hfil #&
\hfil \quad$#$.& $#$\hfil & \hfil \quad$#$.& $#$\hfil & \hfil
\qquad#& \hfil \quad$#$.& $#$\hfil & \hfil \quad$#$.& $#$\hfil \cr
\noalign {\hrule \bigskip}a)& 0& 8& 0& 2& b)& 0& 8& 0& 2\cr &
1& 0& 0& && 1& 0& -0& 3\cr & 0& & 0& && 0& & 0& \cr \noalign
{\bigskip} c)& 0& 8& 0& 2& d)& 1& 0& 0& \cr & 1& 0& -0& 7&&
1& 0& 0& \cr & 0& & 0& && -0& 5& 0& \cr \noalign {\bigskip} e)&
0& 8& 0& 2& f)& 0& 5& 0& 5\cr & 1& 0& 0& && 1& 0& 0& 2\cr &
-0& 5& 0& && -0& 5& 0& 2\cr \noalign {\bigskip \hrule} }}}\par
}\endinsert

As explained in the survey of Section~2, Richardson's 1st order
method, \ocm (1 1), has a circular convergence domain.
Figure~10a shows the typically elliptic domain of the stationary
2nd order method, \ocm (1 2).

Figures~10b and 10c prove that some operator coefficient
methods are new. \Ocm (1 2) methods require just one
matrix-vector multiplication per step, and yet fully populated $2
\times 2$ tableaux have non-elliptic convergence domains
because the relationship that $P(\lambda, X) = 0$ creates
between $\lambda$ and $X$, $$\lambda = - {X^2 - c_{0,1} X -
c_{0,2} \over c_{1,1} X + c_{1,2}} \, ,$$ generally does not map
circles in $X$ to ellipses in $\lambda$ when $c_{1,2} \ne 0$. The
domain in Figure~10b closely abuts the imaginary axis, a difficult
feat for the Chebyshev iteration's ellipses. The crescent-shaped
domain in Figure~10c actually crosses the axis, an impossible
feat for the Chebyshev iteration's ellipses.

The remaining plots in Figure~10 illustrate the possibilities of
tailoring the domain. Figure~10d is for a stationary
gcr($1$)/\discretionary{}{}{}gmres($2$) method, \ocm (2 1).
Figure~10e shows introducing just the previous iterate, \ocm (2
2), produces a square convergence region. Figure~10e is for a
fully populated, $3 \times 2$ tableau.

\vfil \eject \beginsection {References}

{\frenchspacing

\item {[\Aitken ]} A.~C.~Aitken, {\sl On the Iterative Solution of a
System of Linear Equations}, Proceedings of the Royal Society of
Edinburgh Section A, {\bf 63}:52--60 (1950).

\item {[\Arnoldi ]} W.~E.~Arnoldi, {\sl The Principal of Minimized
Iterations in the Solution of the Matrix Eigenvalue Problem},
Quarterly of Applied Mathematics, {\bf 9}:17--29 (1951).

\item {[\Ashby ]} S.~F.~Ashby, T.~A.~Manteuffel and P.~E.~Saylor,
{\sl A Taxonomy for Conjugate Gradient Methods}, Lawrence
Livermore National Laboratory Report UCRL-98508, Livermore,
1988. Submitted to SIAM Journal on Numerical Analysis.

\item {[\Brown ]} P.~N.~Brown and A.~C.~Hindmarsh, {\sl
Matrix-free methods for stiff systems of ODE's}, SIAM Journal on
Numerical Analysis, {\bf 23}:610--638 (1986).

\item {[\Chronopoulosa ]} A.~T.~Chronopoulos, {\sl A class of
parallel iterative methods implemented on multiprocessors},
University of Illinois Doctoral Dissertation, Department of
Computer Science Report UIUCDCS-R-86-1267, Urbana, 1986.

\item {[\Chronopoulosb ]} A.~T.~Chronopoulos, {\sl $s$-Step
Iterative Methods For (Non)symmetric (In)definite Linear
Systems}, Computer Science Department Report TR~89-34,
University of Minnesota, Minneapolis, 1989.

\item {[\Chronopoulosc ]} A.~T.~Chronopoulos and C.~W.~Gear, {\sl
$s$-step iterative methods for symmetric linear systems},
Journal of Computational and Applied Mathematics, {\bf
25}:153--168 (1989).

\item {[\Chronopoulosd ]} A.~T.~Chronopoulos and C.~W.~Gear, {\sl
On the efficient implementation of preconditioned $s$-step
conjugate gradient methods on multiprocessors with memory
hierarchy}, to appear in Parallel Computing.

\item {[\Concusa ]} P.~Concus and G.~H.~Golub, {\sl A generalized
conjugate gradient method for nonsymmetric systems of linear
equations}, in R.~Glowinski and J.~L.~Lions, Lecture Notes in
Economics and Mathematical Systems {\bf 134}, Springer-Verlag,
Berlin, 1976, 56--65.

\item {[\Concusb ]} P.~Concus, G.~H.~Golub and D.~P.~O'Leary, {\sl A
generalized conjugate gradient method for the numerical solution
of elliptic partial differential equations}, in Sparse Matrix
Computations, edited by J.~R.~Bunch and D.~J.~Rose, Academic
Press, New York, 1976, 309--332.

\item {[\Craig ]} E.~J.~Craig, {\sl The $n$-step iteration
procedures}, Journal of Mathematics and Physics, {\bf 24}:64--73
(1955). This journal continues as Studies in Applied
Mathematics.

\item {[\Eisenstata ]} S.~C.~Eisenstat, H.~Elman, M.~H.~Schultz and
A.~H.~Sherman, {\sl Solving approximations to the convection
diffusion equations}, in Proceedings of the Society of Petroleum
Engineers of the AIME Fifth Symposium on Reservoir Simulation,
Denver, 1979, 127--132.

\item {[\Eisenstatb ]} S.~C.~Eisenstat, H.~C.~Elman and
M.~H.~Schultz, {\sl Variational iterative methods for
nonsymmetric systems of linear equations}, SIAM Journal on
Numerical Analysis, {\bf 20}:345--357 (1983).

\item {[\Elmana ]} H.~C.~Elman, {\sl Iterative Methods for Large,
Sparse, Nonsymmetric Systems of Linear Equations}, Yale
University Doctoral Dissertation, New Haven, 1982.
Available from University Microfilms International Dissertation
Information Service, catalogue number 8222744.

\item {[\Elmanb ]} H.~C.~Elman and R.~L.~Streit, {\sl Polynomial
Iteration for Nonsymmetric Indefinite Linear Systems}, Yale
University Department of Computer Science Research Report
YALEU/DCS/RR-380, New Haven, 1985.

\item {[\Engeli ]} M.~Engeli, T.~Ginsburg, H.~Rutishauser and
E.~Stiefel, {\sl Refined Iterative Methods for Computation of the
Solution and Eigenvalues of Self-Adjoint Boundary Value
Problems}, Mitteilungen aus dem Institut f\"ur angewandte
Mathematik Zurich, {\bf 8}, Birkha\"user Verlag, Basel, 1959.

\item {[\Fabera ]} V.~Faber and T.~Manteuffel, {\sl Necessary
and sufficient conditions for the existence of a conjugate
gradient method}, SIAM Journal on Numerical Analysis, {\bf
21}:352--362 (1984).

\item {[\Faberb ]} V.~Faber and T.~A.~Manteuffel, {\sl Orthogonal
Error Methods}, Los Alamos National Laboratory Report
LA-UR-84-3150, Los Alamos, 1984.

\item {[\Faberc ]} V.~Faber and T.~A.~Manteuffel, {\sl Orthogonal
error methods}, SIAM Journal on Numerical Analysis, {\bf
24}:170--187 (1987).

\item {[\Forsythe ]} G.~E.~Forsythe, M.~R.~Hestenes and
J.~B.~Rosser, {\sl Iterative Methods for Solving Linear Equations},
Bulletin of the American Mathematical Society, {\bf 58}:480
(1951).

\item {[\Frankel ]} S.~P.~Frankel, {\sl Convergence Rates of
Iterative Treatments of Partial Differential Equations},
Mathematical Tables and other Aids to Computation, {\bf
4}:65--75 (1950). This journal continues as Mathematics of
Computation.

\item {[\Goluba ]} G.~H.~Golub, {\sl The use of Chebyshev matrix
polynomials in the iterative solution of linear equations
compared with the method of successive overrelaxation},
University of Illinois Doctoral Dissertation, Urbana, 1959.

\item {[\Golubb ]} G.~H.~Golub and R.~S.~Varga, {\sl
Chebyshev semi-iterative methods, successive overrelaxation
iterative methods, and second order Richardson iterative
methods, Part~I}, Numerische Mathematik, {\bf
3}:147--156 (1961).

\item {[\Golubc ]} G.~H.~Golub and R.~S.~Varga, {\sl
Chebyshev semi-iterative methods, successive overrelaxation
iterative methods, and second order Richardson iterative
methods, Part~II}, Numerische Mathematik, {\bf
3}:157--168 (1961).

\item {[\Golubd ]} G.~H.~Golub and D.~P.~O'Leary, {\sl Some history
of the conjugate gradient and Lanczos algorithms: 1948--1976},
SIAM Review, {\bf 31}:50--102 (1989).

\item {[\Golube ]} G.~H.~Golub and C.~F.~Van Loan, {\sl Matrix
Computations}, Johns Hopkins University Press, Baltimore, 1983.

\item {[\Greenbaum ]} A.~Greenbaum, {\sl Behavior of Slightly
Perturbed Lanczos and Conjugate-Gra\-di\-ent Recurrences},
Linear Algebra and Its Applications, {\bf 113}:7--63 (1989).

\item {[\Hestenesa ]} M.~R.~Hestenes and E.~L.~Stiefel, {\sl
Methods of Conjugate Gradients for Solving Linear Systems},
Journal of Research of the National Bureau of Standards Section
B, {\bf 49}:409--436 (1952).

\item {[\Hestenesb ]} M.~R.~Hestenes, {\sl Conjugate Direction
Methods in Optimization}, Springer-Ver\-lag, New York, 1980.

\item {[\Householder ]} A.~S.~Householder, {\sl The Theory of
Matrices in Numerical Analysis}, Blaisdell, New York, 1964.
Reprinted by Dover, New York, 1975.

\item {[\Jea ]} K.~C.~Jea and D.~M.~Young, {\sl On the
Simplification of Generalized Conjugate-Gradient Methods for
Nonsymmetrizable Linear Systems}, Linear Algebra and Its
Applications, {\bf 52/53}:399--417 (1983).

\item {[\Joubert ]} W.~D.~Joubert and D.~M.~Young, {\sl
Necessary and Sufficient Conditions for the Simplification of
Generalized Conjugate-Gradient Algorithms}, Linear Algebra and
Its Applications, {\bf 88/89}:449--485 (1986).

\item {[\Luenbergera ]} D.~G.~Luenberger, {\sl The conjugate
residual method for constrained minimization problems}, SIAM
Journal on Numerical Analysis, {\bf 7}:390--398 (1970).

\item {[\Luenbergerb ]} D.~G.~Luenberger, {\sl Introduction to
linear and nonlinear programming}, Addison-Wesley, Reading
Massachusetts, 1973.

\item {[\Manteuffela ]} T.~A.~Manteuffel, {\sl An iterative method
for solving nonsymmetric systems with dynamic estimation of
parameters}, University of Illinois Doctoral Dissertation,
Department of Computer Science Report UIUCDCS-R-75-758,
Urbana, 1975.

\item {[\Manteuffelb ]} T.~A.~Manteuffel, {\sl The Tchebychev
iteration for nonsymmetric linear systems}, Numerische
Mathematik, {\bf 28}:307--327 (1977).

\item {[\Manteuffelc ]} T.~A.~Manteuffel, {\sl Adaptive procedure
for estimating parameters for the nonsymmetric Tchebychev
iteration}, Numerische Mathematik, {\bf 31}:183--208 (1978).

\item {[\Manteuffeld ]} T.~A.~Manteuffel, {\sl Optimal parameters
for linear second-degree stationary iterative methods},
SIAM Journal on Numerical Analysis, {\bf 19}:883--839 (1982).

\item {[\Milne ]} L.~M.~Milne-Thomson, {\sl The calculus of finite
differences}, Macmillan, London, 1933. Reprinted by Chelsea, New
York, 1981.

\item {[\Niethammer ]} W.~Niethammer and R.~S.~Varga, {\sl The
Analysis of $k$-Step Iterative Methods for Linear Systems from
Summability Theory}, Numerische Mathematik, {\bf 41}:177--206
(1983).

\item {[\Opfer ]} G.~Opfer and G.~Schober, {\sl Richardson's
Iteration for Nonsymmetric Matrices}, Linear Algebra and Its
Applications, {\bf 58}:343--361 (1984).

\item {[\Reid ]} J.~K.~Reid, {\sl On the Method of Conjugate
Gradients for the Solution of Large Sparse Systems of Linear
Equations}, in Large Sparse Sets of Linear Equations, Academic
Press, New York, 1971, 231--254.

\item {[\Richardson ]} L.~F.~Richardson, {\sl The approximate
arithmetical solution by finite differences of physical problems
involving differential equations, with applications to the stress
in a masonry dam}, Philosophical Transactions of the Royal
Society of London Series A, {\bf 210}:307--357 (1910).

\item {[\Saad ]} Y.~Saad, {\sl Practical use of some Krylov
subspace methods for solving indefinite and nonsymmetric linear
systems}, SIAM Journal on Scientific and Statistical Computing,
{\bf 5}:203--228 (1984).

\item {[\Saada ]} Y.~Saad and M.~H.~Schultz, {\sl GMRES: A
Generalized Minimum Residual Algorithm for Solving
Nonsymmetric Linear Systems}, Yale University Technical Report
254, New Haven, 1983.

\item {[\Saadb ]} Y.~Saad and M.~H.~Schultz, {\sl Conjugate
Gradient-Like Algorithms for Solving Nonsymmetric Linear
Systems}, Mathematics of Computation, {\bf 44}:417--424 (1985).

\item {[\Saadc ]} Y.~Saad and M.~H.~Schultz, {\sl GMRES: A
Generalized Minimum Residual Algorithm for Solving
Nonsymmetric Linear Systems}, SIAM Journal on Scientific and
Statistical Computing, {\bf 7}:856--869 (1986).

\item{[\Shah ]} B.~Shah, R.~Buehler and O.~Kempthorne, {\sl Some
algorithms for minimizing a function of several variables},
Journal of the Society of Industrial and Applied Mathematics, {\bf
12}:74--92 (1964). This journal continues as SIAM Journal on
Applied Mathematics.

\item {[\Stiefela ]} E.~L.~Stiefel, {\sl Relaxationsmethoden
bester strategie zur losung linearer gleichungssystems},
Commentarii Mathematici Helvetici, {\bf 29}:157--179 (1955).

\item {[\Stiefelb ]} E.~L.~Stiefel, {\sl Kernel Polynomials in Linear
Algebra and Their Numerical Applications}, National Bureau of
Standards Applied Mathematics Series, {\bf 49}:1--22 (1958).

\item {[\Trefethen ]} L.~N.~Trefethen, {\sl Non-normal matrices,
``approximate eigenvalues'', and numerical algorithms},
Contributed Paper 172, SIAM Annual Meeting, San Diego, 1989.

\item {[\Vargaa ]} R.~S.~Varga, {\sl A comparison of the
successive overrelaxation method and semi-iterative methods
using Chebyshev polynomials}, SIAM Journal, {\bf 5}:39--46
(1957).

\item {[\Vargab ]} R.~S.~Varga, {\sl Matrix Iterative Analysis},
Prentice-Hall, Englewood Cliffs, 1962.

\item{[\Vinsome ]} P.~K.~W.~Vinsome, {\sl ORTHOMIN, an iterative
method for solving sparse sets of simultaneous linear equations},
in 4th Symposium of the Society of Petroleum Engineers of the
AIME, Los Angeles, 1976, 149--159. 

\item{[\Voevodin ]} V.~V.~Voevodin, {\sl The problem of a
non-self-adjoint generalization of the con\-ju\-gate-gradient
method has been closed}, U.S.S.R.~Computational Mathematics and
Mathematical Physics, {\bf 23}:143--144 (1983).

\item{[\Walkera ]} H.~F.~Walker, {\sl Implementation of the gmres
method using Householder transformations}, SIAM Journal on
Numerical Analysis, {\bf 9}:152--163 (1988).

\item{[\Walkerb ]} H.~F.~Walker, {\sl Implementations of the
gmres method, Computer Physics Communications}, {\bf
53}:311--320 (1989).

\item {[\Widlund ]} O.~Widlund, {\sl A Lanczos method for a class
of non-symmetric systems of linear equations}, SIAM Journal on
Numerical Analysis, {\bf 15}:801--812 (1978).

\item {[\Younga ]} D.~M.~Young, {\sl Iterative methods for solving
partial differential equations of elliptic type}, Harvard
University Doctoral Dissertation, Cambridge, 1950.

\item {[\Youngb ]} D.~M.~Young, {\sl Iterative methods for solving
partial differential equations of elliptic type}, Transactions of
the American Mathematical Society, {\bf 76}:92--111 (1954).

\item {[\Youngc ]} D.~M.~Young, {\sl Iterative Solution of Large
Linear Systems}, Academic Press, New York, 1971.

\item {[\Younge ]} D.~M.~Young, {\sl A historical overview of
iterative methods}, Computer Physics Communications, {\bf
53}:1--17 (1989).

\item {[\Youngd ]} D.~M.~Young and K.~C.~Jea, {\sl Generalized
Conjugate-Gradient Acceleration of Nonsymmetrizable Iterative
Methods}, Linear Algebra and Its Applications, {\bf
34}:159--194 (1980).

}

\vfil \eject \beginsection {Appendix 1. Proofs}

This appendix proves the theorems cited in the text.

\Proclaim Theorem 1. If the Hermitian part $H$ of $P(A)$ is positive
or negative definite for some polynomial $P$ with degree at most
$k$ and $P(0) = 0$, then for every $x_n$ the affine space $$x_n +
\Span {r_n \;\; A r_n \;\; A^2 r_n \;\; \ldots \;\; A^{k-1} r_n}$$
contains a vector $x_{n+1}$ with $\| r_{n+1} \|_2 \le \rho \, \| r_n
\|_2$ where $$\rho = \sqrt {1 - \left[ \, {\min | \lambda (H) | \over \|
P(A) \|_2} \, \right]^2} \; < \; 1.$$ The affine space also contains a
vector $x_{n+1}$, usually different from the first, with $\| e_{n+1}
\|_2 \le \rho \, \| e_n \|_2$.

{\sl Proof.} This proof generalizes the one Elman [\Elmana ] and
Eisenstat, Elman and Schultz [\Eisenstatb ] use to prove convergence
of truncated orthomin. Its application in this context, and the part
about minimizing the error, appear to be new.

If $r_n = 0$ then choose $x_{n+1} = x_n$. Otherwise let
$P(X) = c_1 X + c_2 X^2 + \cdots + c_k X^k$ and choose $$x_{n+1} =
x_n + \alpha (c_1 + c_2 A + \cdots + c_k A^{k-1}) r_n$$ for some
real $\alpha$. This expression is special because $r_{n+1} = [I -
\alpha P(A)] \, r_n$ so $$\| r_{n+1} \|^2 = \| r_n \|^2 - 2 \alpha
\transpose {r_n} H r_n + \alpha^2 \| P(A) r_n \|^2.$$ Now, $\| P(A) r_n
\| \ne 0$ because $P(A) r_n \ne 0$ because $\transpose {r_n} P(A) r_n
= \transpose {r_n} H r_n \ne 0$ so the formula is quadratic in
$\alpha$. The minimum occurs at $\alpha = \transpose {r_n} H r_n / \|
P(A) r_n \|^2$ and is $$\| r_{n+1} \|^2 \; = \; \| r_n \|^2 - \left[ \,
{\transpose {r_n} H r_n \over \| P(A) r_n \|} \, \right]^2 \; \le \;\, \| r_n
\|^2 - \left[ \, {\min | \lambda (H) | \over \| P(A) \|} \; \| r_n \| \;
\right]^2 .$$ There is no need to ponder the size of $${\min | \lambda
(H) | \over \| P(A) \|} \, .$$ It must be  $\le 1$ because $\| r_{n+1} \|^2$
isn't negative. And it must be $\ge 0$ because it is an absolute value
over a norm. But it can't be zero because $H$ is positive or negative
definite so the numerator won't vanish. 

As for minimizing the error, the choice $$x_{n+1} = x_n + \alpha (c_1
+ c_2 A + \cdots + c_k A^{k-1}) r_n$$ also yields $e_{n+1} = [I -
\alpha P(A)] \, e_n$. This vector's norm can be minimized in the same
manner as the residual's, but the minimum most likely occurs for a
different~$\alpha$. How different? \halmos

\Proclaim Theorem 2. A constant coefficient, homogeneous,
operator coefficient method of degree $k$ and order $m$ $$x_n =
\sum_{j=1}^m c_{0,j} x_{n-j} + \sum_{i=1}^k \sum_{j=1}^m c_{i,j}
A^{i-1} r_{n-j}$$ with coefficient tableau $$\matrix {c_{0,1}&
c_{0,2}& \cdots& c_{0,m}\cr c_{1,1}& c_{1,2}& \cdots& c_{1,m}\cr
\vdots& \vdots& \ddots& \vdots\cr c_{k,1}& c_{k,2}& \cdots&
c_{k,m}\cr}$$ converges to a solution of $A x = y$ for all $y$ and
all initial vectors $x_{0}$, $x_{-1}$, $\ldots\,$, $x_{1-m}$ exactly
when, for each eigenvalue $\lambda$ of $A$, the maximum
magnitude $r(\lambda )$ of the roots $X$ of the polynomial
$$P(\lambda, X) = X^m - P_1(\lambda ) X^{m-1} - P_2(\lambda )
X^{m-2} - \cdots - P_m(\lambda ) X^{m-m}$$ with coefficients
given by the columns of the tableau $$P_j (\lambda ) = c_{0,\,j} -
c_{1,\,j} \lambda - c_{2,\,j} \lambda^2 - \cdots - c_{k,\,j}
\lambda^k$$ is strictly less than $1$. Moreover, there is a bound
upon the residuals for all $y$ and all initial vectors $x_{0}$,
$x_{-1}$, $\ldots\,$, $x_{1-m}$ $$\| r_n \| \le \left( \, \| r_{0} \| + \|
r_{-1} \| + \cdots + \| r_{1-m} \| \, \right) Q(n) \, R^n,$$ and if $A$
is nonsingular there is an identical bound upon the errors $$\| e_n
\| \le \left( \, \| e_{0} \| + \| e_{-1} \| + \cdots + \| e_{1-m} \| \,
\right) Q(n) \, R^n.$$ $R$ is the maximum $r(\lambda )$ for all the
eigenvalues of $A$. $Q(n)$ is a polynomial that depends on the norm,
on $A$, and on the coefficient tableau. The norm may be any
consistent matrix-vector norm.

{\sl Proof.} The Theorem is not trivial because it allows defective
$A$, but the proof's main challenges are notation and pruning. There
are nine parts. 

{\sl Part~1.} Much of what is needed to prove the Theorem can be
taken from the literature of finite differences [\Milne ]. The proof
depends on sequences $\{ \tau_n \}$ generated  by constant
coefficient, homogeneous, $m^{th}$ order, linear recurrence
formulas $$\tau_n = \sum_{j=1}^m P_j \tau_{n-j}.$$ The sequences
begin from initial values $\tau_0$, $\tau_{-1}$, $\ldots\,$,
$\tau_{1-m}$ and the formulas apply when $n \ge 1$. Any recurrence
sequence can be expressed as a linear combination of {\sl
fundamental sequences} of the form $\{ n^{i-1} \rho^n \}$ in which
$\rho$ is a root of multiplicity at least $i$ of the {\sl characteristic
polynomial} $$X^m - \sum_{j=1}^m P_j X^{m-j}.$$ The coefficients in
the linear combination for $\{ \tau_n \}$ can be obtained by solving
an $m \times m$ system of linear equations. Each column of the
coefficient matrix $M$ contains the initial values of a different
fundamental sequence, and the column on the right side of the
equations contains the initial values of $\{ \tau_n \}$. This formula
for the coefficients leads to the following bound $$| \tau_n | \le
n^{m-1} r^n \| M^{-1} \|_1 \sum_{j=1}^{m} |\tau_{1-j}|$$ in which $n
\ge 1$ and $r$ is the magnitude of the largest root $\rho$. The power
$n^{m-1}$ only occurs in the worst case when the characteristic
polynomial has a single root of multiplicity $m$. 

The bound is needed for recurrence sequences whose coefficients
$P_j$ are polynomials of a parameter $\lambda$. In this case, the
characteristic polynomial is the $P(\lambda, X)$ in the statement of
the Theorem. A technical detail occurs when $P_m$ vanishes for
some $\lambda$, because then the recurrence has order less than
$m$. The bound remains valid for a smaller matrix $M$ and a
smaller order $m$, but it is notationally convenient to retain the
original $m$. Thus, for specific polynomials $P_j (\lambda)$ there is
a bound $$| \tau_n | \le n^{m-1} r(\lambda)^n q (\lambda)
\sum_{j=1}^m |\tau_{1-j}|$$ in which $r(\lambda)$ is the magnitude of
the largest root of the characteristic polynomial, and $q(\lambda)$
is a number with a rather complicated definition. The $r(\lambda)$
and $q(\lambda)$ depend on the recurrence coefficients, that is, on
the coefficient tableau of the operator coefficient method.

{\sl Part~2.} The recurrence formulas generate another collection
of fundamental sequences, here denoted $\{ \pi_{n,\,k} \}$ for $k$
from $1$ to $m$. These sequences have only one nonzero initial value
apiece, namely $\pi_{1-k,\,k} = 1$ and others $0$. When they
represent any other recurrence sequence in the manner of Part~1,
then the coefficients in the linear combination are just the other's
initial values. $$\tau_n = \sum_{k=1}^m \tau_{1-k} \pi_{n,\,k}$$ In
this way a homogeneous sequence has a closed-form representation
terms of its initial values.

{\sl Part~3.} Adding an extra term $\sigma_\ell$ to $\tau_{\,\ell}$
amounts to beginning a new recurrence sequence with initial value
$\sigma_\ell$. The contribution to $\tau_{\,\ell +1}$ is $\pi_{1,\,1}
\sigma_\ell$, the contribution to $\tau_{\,\ell + 2}$ is $\pi_{2,\,1}
\sigma_\ell$, and so on. The sequence of multipliers is $\{ \pi_{n,\,1}
\}$. One sequence suffices independent of $\ell$ because the
recurrence formulas have constant coefficients. In combination with
Part~2 therefore, a nonhomogeneous recurrence $$\tau_n =
\sigma_n + \sum_{j=1}^m P_j \tau_{n-j}$$ has a closed-form
representation $$\tau_n = \sum_{\ell=1}^n \pi_{n-\ell,\,1}
\sigma_\ell + \sum_{k=1}^m \pi_{n,\,k} \tau_{1-k}.$$ 

{\sl Part~4.} When the recurrence coefficients $P_1$, $P_2$,
$\ldots\,$, $P_m$ are polynomials of $\lambda$ as they are here, then
so are the $\pi_{n,\,k}$'s. Their recurrence formula $$\pi_{n,\,k} =
\sum_{j=1}^m P_j \pi_{n-j,\,k}$$ can be differentiated $s$ times
$$\eqalign {\pi_{n,\,k}^{(s)}& = \sum_{i=0}^s \sum_{j=1}^m \left( {s
\atop i} \right) P_j^{(i)} \pi_{n-j,\,k}^{(s-i)}\cr &= \sum_{j=1}^m  P_j
\pi_{n-j,\,k}^{(s)} + \sum_{i=1}^s \sum_{j=1}^m \left( {s \atop i}
\right) P_j^{(i)} \pi_{n-j,\,k}^{(s-i)}\cr}$$ to reveal that the
derivatives satisfy nonhomogeneous recurrence formulas. The
superscripts in parentheses denote derivatives of various orders
with respect to $\lambda$. The derivative sequences have zeroes for
initial values, that is, only the nonhomogeneous terms participate in
Part~3's closed-form expansion. $$\pi_{n,\,j}^{(s)} =
\sum_{\ell=1}^n \sum_{i=1}^s \sum_{j=1}^m \left( {s \atop i} \right)
\pi_{n-\ell,1}^{\vphantom {(i)}} P_j^{(i)} \pi_{n-j,\,k}^{(s-i)}$$ This
represents each derivative entirely in terms of lower-order
derivatives.

{\sl Part~5.} The bounds of Part~1 can be applied first to the
polynomials $\pi_{n,\,k}$ of Part~2 and then to the derivative
formulas of Part~5 to obtain $$\left| \pi_{n,\,k}^{(s)} (\lambda)
\right| \le q_s (n, \lambda) \, r(\lambda)^n$$ where $q_s$ is a
polynomial of $n$. Part~1 forces the choice $$q_0 (n, \lambda) =
n^{m-1} q(\lambda ),$$ and the others can be constructed recursively
by $$q_s (n, \lambda) = n \, q_{0} (n, \lambda) \, q_{s-1} (n, \lambda)
\, m \, 2^s \max_{1 \le i \le s, \; 1 \le j \le m} \left| P_j^{(i)} (\lambda )
\right| .$$ Only the polynomial dependence on $n$ is important. The
bounds aren't sharp and needn't be. They have been chosen to
increase monotonically with $n$ and $s$ to ease the following 
derivation. $$\eqalign {\left| \pi_{n,\,k}^{(s)} (\lambda) \right| &=
\Big| \sum_{\ell=1}^n \sum_{i=1}^s \sum_{j=1}^m \left( {s \atop i}
\right) \pi_{n-\ell,1}^{\vphantom {(i)}} (\lambda ) P_j^{(i)} (\lambda)
\pi_{n-j,\,k}^{(s-i)} (\lambda) \Big| \cr &\le \; \sum_{\ell=1}^n
\sum_{i=1}^s \sum_{j=1}^m \left( {s \atop i} \right) q_0 (n, \lambda )
\left| P_j^{(i)} (\lambda ) \right| q_{s-1} (n, \lambda)\cr &\le \; n \,
q_0 (n, \lambda ) \, q_{s-1} (n, \lambda) \sum_{i=1}^s \sum_{j=1}^m
\left( {s \atop i} \right) \left| P_j^{(i)} (\lambda ) \right| \cr &\le \;
q_s (n, \lambda)\cr}$$

{\sl Part~6.} At this point its customary to cite an unimpeachable
source for the definition and existence of a Jordan decomposition,
$J = U A \, U^{-1}$. This allows individual Jordan blocks to be
considered. Manteuffel [\Manteuffela ] [\Manteuffelb ] observes that
evaluating a polynomial $\pi (X)$ at a Jordan block amounts to
differentiating the polynomial. If $\lambda$ is the block's
eigenvalue, there results an upper triangular, Toeplitz matrix with
$\pi^{(i)} (\lambda ) / i\,!$ on the $i^{\,th}$ superdiagonal, as in the
example below. \def \space{\vphantom {\hbox {$\Big|$}}} $$\pi \left(
\left[ \matrix {\space \lambda& 1\cr \space & \lambda& 1\cr \space
&& \lambda\cr} \right] \right) = \left[ \matrix {\space \pi^{(0)}
(\lambda ) / 0!& \pi^{(1)} (\lambda ) / 1!& \pi^{(2)} (\lambda ) / 2!\cr
\space &\pi^{(0)} (\lambda ) / 0!& \pi^{(1)} (\lambda ) / 1!\cr &&
\space \pi^{(0)} (\lambda ) / 0!\cr} \right]$$ In this way the bounds of
Part~5 can be applied to each Jordan block and then combined for
all blocks to obtain bounds $$\| \pi_{n,\,k} (A) \| = \| U^{-1} \pi_{n,\,k}
(J) \, U \| \le Q (n) R^n.$$ $R$ is the largest $r(\lambda)$ for all the
matrix eigenvalues. $Q(n)$ is a polynomial that depends on the norm,
on $A$, and on the coefficient tableau. The norm should be a
consistent matrix-vector norm because the next step applies it to
both.

{\sl Part~7.} Remembering the original use of the $\pi_{n,\,k}$'s in
Part~2 and the bounds in Part~6, the terms of vector sequences $\{
t_n \}$ generated from initial values $t_0$, $t_{-1}$, $\ldots\,$,
$t_{1-m}$ by a constant coefficient, homogeneous, $m^{th}$~order,
linear recurrence formula $$t_n = \sum_{j=1}^m P_j(A) t_{n-j}, $$
whose coefficients are polynomials of $A$, are expressed by $$t_n =
\sum_{k=1}^m \pi_{n,\,k} (A) t_{1-k}$$ and are bound by $$\| t_n \| \;
\le \; Q(n) R^n \sum_{k=1}^m \left \| t_{1-k} \right\|.$$

{\sl Part~8.} The preliminaries are finished and the theorem's proof
begins here. Since the operator coefficient method is homogeneous
there is a residual recurrence $$r_n = \sum_{j=1}^m P_j(A) r_{n-j}$$
and if $A$ is nonsingular this can be multiplied by $A^{-1}$ to obtain
an error recurrence $$e_n = \sum_{j=1}^m P_j(A) e_{n-j}$$ and from
Part~7 both sequences satisfy the bounds in the statement of the
Theorem. 

{\sl Part~9.} If $A$ is singular then the Theorem's polynomial for
$\lambda = 0$ is $$P(0, X) = X^m - c_{0,1} X^{m-1} - c_{0,2} X^{m-2} -
\cdots - c_{0,m} X^{m-m}$$ which has $1$ as a root because $c_{0,1}
+ c_{0,2} + \cdots + c_{0,m} = 0$ in the homogeneous case, so $1 =
r(0) \le R$. Thus, if $R < 1$ then $0$ is not an eigenvalue, $A$ is
nonsingular, and a solution exists. The bound in Part~8 shows the
errors converge to zero since $$\lim_n Q (n) R^n = 0$$ whenever $Q$
is a polynomial of $n$. The Theorem's convergence criterion is
therefore sufficient. If $1\le R$, then from the Jordan block of the
eigenvalue for which the Theorem's polynomial has a root
of magnitude $R$, it is possible to construct some $y$ and initial
vectors $x_0$, $x_{-1}$, $\ldots\,$, $x_{1-m}$ so the corresponding
residuals satisfy $\| r_n \| = \| r_0 \| R^n$. The Theorem's
convergence criterion therefore is necessary. \halmos

\vfil \eject \beginsection {Appendix 2. Figure Explanations}

This appendix explains the numerical experiments reported in the
Figures. All calculations are performed by a Cray XMP with unit
roundoff $3{.}5 \times 10^{-15}$. The initial guess $x_0$ for all
iterations is $0$.

\Proclaim Figure 1. $2$-norm relative errors in the computed basis
vectors of the conjugate gradient algorithm for a system of order
$100$.

\Proclaim Figure 2. Relative $A$-norm solution errors for the
system of Figure~1. The upper curve is for single precision and the
lower for reorthogonalized double precision. 

For Figures~1 and 2 the system $A x = y$ has diagonal $A$ with
entries $1^2$, $2^2$, $\ldots\,$, $100^2$ and uniform $y$ with entries
$1$, $1$, $\ldots\,$, $1$. The single precision  conjugate gradient
algorithm appears in Section~2 and has explicit, not recursive,
residuals. Its deviation from what would be obtained from exact,
infinite precision calculations is measured by comparison with a
double precision version that includes full orthogonalization in
which $$p_n = r_n - \sum_{j=1}^{n} {\transpose {p_{n-j}} A r_n
\over \transpose {p_{n-j}} A p_{n-j}} \, p_{n-j}$$ replaces $$p_n =
r_n - {\transpose {p_{n-1}} A r_n \over \transpose {p_{n-1}} A
p_{n-1}} \, p_{n-1} .$$ Figure~1 actually plots $\| p_n^{(1)} - p_n^{(2)}
\|_2 / \| p_n^{(2)} \|_2$ and Figure~2 plots both $\| e_n^{(1)} \|_A / \|
x_* \|_A$ and $\| e_n^{(2)} \|_A / \| x_* \|_A$. The superscripts
distinguish single and double precision values.

\Proclaim Figure 3. $2$-norm relative residuals for homogeneous
(dashed) and inhomogeneous (solid)
gcr\/$(k{-}1)$/\discretionary{}{}{}gmres\/$(k)$, $k = 1$, $2$,
$\ldots\,$, $10$, applied to one system. The two methods perform
alike except for $k = 5$ when the original, homogeneous method
stagnates and the new, inhomogeneous method converges.

For Figures~3, 4, 5 and 6 the system $A x = y$ resembles one used
by Elman and Streit [\Elmanb ]. The matrix is a preconditioned
discretization of $$\left( \, - {\partial^2 \over \partial x \partial y}
+ \alpha \, {\partial \over \partial x} + \beta \, {\partial \over
\partial y} - \gamma \, \right) u = f$$ for real-valued $u$ on $[0,1]
\times [0,1]$ with zero Dirichlet boundary data. The finite difference
discretization has a $33 \times 33$ uniform grid with $961$ interior
unknowns related by $5$-point approximations to second
derivatives and by centered approximations to first derivatives.
The preconditioned system is $A_2^{-1} A_1 x = y$ where $A_1$ is
the discrete operator for $\alpha = 50$, $\beta = 100$, $\gamma =
250$ and where $A_2$ is the discrete operator for $\alpha = \beta =
\gamma = 0$. A stabilized block cyclic reduction method performs
the multiplication by $A_2^{-1}$ [\Golube ]. The entries of $y$ are
uniformly and randomly distributed between $-1$ and $1$.

Gcr($k{-}1$)/\discretionary{}{}{}gmres($k$) is implemented in the
simplest possible manner suggested by Table~1. The singular value
decomposition solves the least squares problem $$\hbox {minimize
$\| r_{n+1} \| \;$ over $\; x_n + \Span {r_n \;\; Ar_n \;\; A^2 r_n \;\;
\ldots \;\; A^{k-1} r_n}$}$$ by projecting $r_n$ into $\Span {A r_n \;\;
A^2 r_n \;\; \ldots \;\; A^k r_n}$. The inhomogeneous version solves
$$\hbox {minimize $\| r_{n+1} \| \;$ over $\; \Span {x_n \;\; r_n \;\;
Ar_n \;\; A^2 r_n \;\; \ldots \;\; A^{k-1} r_n}$}$$ by projecting $y$
into $\Span {A x_n \;\; A r_n \;\; A^2 r_n \;\; \ldots \;\; A^k r_n}.$

\Proclaim Figure 4. $2$-norm relative residuals for homogeneous
(dashed) gcr$\,(5)$/\discretionary{}{}{}gmres$\,(6)$ and
inhomogeneous (solid) \ocm (6 m), $m = 1$ $2$ $\ldots$ $10$,
applied to the same system.

The system of equations and the implementation of homogeneous
gcr($k{-}1$)/\discretionary{}{}{}gmres($k$) are described with
Figure~3. The implementation of the inhomogeneous \ocm (k m)
methods is described in Section~6b. When $m=1$ it is identical to
Figure~3's in\-homogeneous
gcr($k{-}1$)/\discretionary{}{}{}gmres($k$).

\Proclaim Figure 5. Level curves of observed convergence
rate for inhomogeneous \ocm (k m) as a function of $k$ and $m$
for one system. The curves range in multiplicative steps of
$10^{0.05}$ from $10^{-1}$ at the upper right to $10^0$, no
convergence, at the lower left. The left edge corresponds to
gcr\/$(k{-}1)$/\discretionary{}{}{}gmres\/$(k)$, the bottom edge
to orthomin\/$(m{-}1)$. Figure~6 displays the convergence
histories. Figure~6 displays the convergence histories.

\Proclaim Figure 6. Convergence histories from which the level
curves of Figure~5 are derived. 

The system of equations is the one for Figures~3 and 4. Section~6b
describes the implementation of the inhomogeneous \ocm (k m)
methods. Least squares linear fits to the last half of each curve in
Figure~6's logarithmic scale produce the observed convergence
rates. The level curves are obtained by extending the fits to the
logarithmic data bilinearly throughout the cells of the $k \times m$
grid. 

\Proclaim Figure 7. Coefficients of inhomogeneous \ocm (2 2)
minimizing the $2$-norm of the residual for one system.

\Proclaim Figure 8. Eigenvalues of the matrix (solid dots)
superimposed on some level curves of the convergence rate
$r(\lambda)$ for the constant coefficient regime of Figure~7.  The
levels start at $1$ on the boundary and decrease in steps of $0.05$.

\Proclaim Figure 9. $2$-norm relative residuals for the iteration of
Figure~7 (lower solid line), and for a constant coefficient iteration
with the same right hand side (higher solid line) and a different right
hand side (dashed line).

For Figures~7, 8 and 9 the matrix is Toeplitz and banded with $-
1$'s on the first superdiagonal and with $1$'s on the main and first
$3$ subdiagonals. Careless programming resulted in a matrix
of order $201$. For Figure~7 the entries of the right side all equal
$1$. The implementation of the inhomogeneous \ocm (2 2) method is
described in Section~6b. Figure~7 plots the $6$ coefficients in the
$3 \times 2$ tableaux. The coefficients for the first few tableaux are
not shown because they are larger than the others. The coefficients
for the next several iterations vary on the order of one percent
from $$\vcenter {\normalbaselines \halign {\hfil$#$.& $#$\hfil & \hfil
\quad$#$.& $#$\hfil\cr 1& 421& -0&421\cr 0& 261& -0& 172\cr -0&
130& 0& 102\cr}}$$ which are the constant coefficients used for
Figures~8 and 9. These Figures are explained by their captions and
by Section~8. The other right hand side for Figure~9 has entries
uniformly and randomly distributed between $-1$ and $1$.

\Proclaim Figure 10. Level curves of convergence rate $r(\lambda )$
as a function of $\lambda$ in the complex plane for Table~4's
constant coefficients. The levels begin at $1$ on the boundaries and
decrease in steps of $0.05$.

Like Figure~8, the curves are formed by plotting the contours of
$$r(\lambda ) = \hbox {maximum $|X|$ of all $X$ for which
$P(\lambda, X) = 0$}$$ using the quadratic formula to solve
$P(\lambda, X) = 0$ when $m=2$.

\vfil \eject

 \newdimen\fullhsize
 \baselineskip = 12pt
 \fullhsize = 5.0in
 \hoffset = 0.8125in
 \hsize=2.45in
 \parskip = 6pt plus 1pt minus 1pt
 \topskip =18pt


 \def\fullline{\hbox to\fullhsize}

 \let\lr=L \newbox\leftcolumn
 \output={\if L\lr
      \global\setbox\leftcolumn=\columnbox \global\let\lr=R
   \else \doubleformat \global\let\lr=L\fi
   \ifnum\outputpenalty>-20000 \else\dosupereject\fi}

 \def\doubleformat{\shipout\vbox{\makeheadline
      \fullline{\box\leftcolumn\hfil\columnbox}
   \makefootline}
   \advancepageno}

 \def\columnbox{\leftline{\pagebody}}

 \def\makeheadline{\vbox to 0pt{\vskip-22.5pt
   \hbox to\fullhsize{\vbox to 8.5pt{}\the\headline}\vss}\nointerlineskip}

 \def\makefootline{\baselineskip=24pt
   \hbox to\fullhsize{\the\footline}}

 \headline={
 \ifnum \pageno=1{\tenrm \hfil Invited for Publication by SIAM Journal on
Matrix Analysis and Applications\hfil}\else {
 \ifnum \pageno=2{\hfil}\else {
 \ifnum \pageno=4{\hfil}\else {
 \ifodd \pageno {\sevenrm \hfil \ \hfil}\else {
 \hfil \sevenrm \ \hfil}\fi}\fi}\fi}\fi}
 
 \def\thismonth{\ifcase \month \or JANUARY\or FEBRUARY\or
MARCH\or APRIL\or MAY\or JUNE\or JULY\or AUGUST\or
SEPTEMBER\or OCTOBER\or NOVEMBER\or DECEMBER\fi}

 \def\footremark{{\ }}
 \def\pagenumber{{\tenrm \the \pageno}}
 \footline={\ifnum \pageno=1{\noindent \hfil \footremark \hfil {\tenrm
1/2}}\else {\ifnum \pageno=2{\hfil}\else {\ifnum \pageno=3{\noindent
\hfil \footremark \hfil {\tenrm 3/4}}\else {\ifnum \pageno=4{\hfil}\else
{\ifodd \pageno {\noindent \hfil \footremark \hfil \pagenumber}\else
{\noindent \pagenumber \hfil \footremark \hfil}\fi}\fi}\fi}\fi}\fi}

 \def \eol {\hfil \break}
 \def \Eblock#1\par{\noindent\frenchspacing\helvetica#1\filbreak}
 \def \Lblock#1\par{\noindent\frenchspacing\helvetica#1\filbreak}
 \def \Sblock{\frenchspacing\helvetica\vskip 6pt plus 1pt minus 1pt}

\line {UNLIMITED RELEASE \eol}
\vglue 6pt plus 1pt minus 1pt 
\line {INITIAL DISTRIBUTION \eol}

{\pretolerance = 10000

\Eblock{I. K. Abu-Shumays\eol
Bettis Atomic Power Lab.\eol
Box 79\eol
West Mifflin, PA 15122}

\Eblock{Loyce M. Adams\eol
Univ. of Washington\eol
Dept. of Applied Math.\eol
Seattle, WA 98195}

\Eblock{Peter W. Aitchison\eol
Univ. of Manitoba\eol
Applied Math. Dept.\eol
Winnipeg R3T 2N2\eol
Manitoba, Canada}

\Eblock{Leena Aittoniemi\eol
Tech. Univ. Berlin\eol
Comp. Sci. Dept.\eol
Franklinstrasse 28-29\eol
D-1000 Berlin 10\eol
West Germany}

\Eblock{Fernando L. Alvarado\eol
Univ. of Wisconsin\eol
Electrical and Comp. Eng.\eol
1425 Johnson Dr.\eol
Madison, WI 53706}

\Eblock{Patrick Amestoy\eol
CERFACS\eol
42 ave g Coriolis\eol
31057 Toulouse\eol
France}

\Eblock{Ed Anderson\eol
Argonne National Lab.\eol
Math. and Comp. Sci. Div.\eol
9700 S Cass Av.\eol
Argonne, IL 60439}

\Eblock{Johannes Anderson\eol
Brunel Univ.\eol
Dept. of Math. and Statistics\eol
Uxbridge Middlesex UB8 3PH\eol
United Kingdom}

\Eblock{Peter Arbenz\eol
ETH-Zentrum\eol
Inst. fur Informatik\eol
CH-8092 Zurich\eol
Switzerland}

\Eblock{Mario Arioli\eol
CERFACS\eol
42 ave g Coriolis\eol
31057 Toulouse\eol
France}

\Eblock{Cleve Ashcraft\eol
Yale Univ.\eol
Dept. of Comp. Sci.\eol
P. O. Box 2158, Yale Station\eol
New Haven, CT 06520}

\Eblock{Owe Axelsson\eol
Katholieke Univ.\eol
Dept. of Math.\eol
Toernooiveld\eol
6525 ED Nijmegen\eol
The Netherlands}

\Eblock{Zhaojun Bai\eol
New York Univ.\eol
Courant Institute\eol
251 Mercer Street\eol
New York, NY 10012}

\Eblock{R. E. Bank\eol
Univ. of California, San Diego\eol
Dept. of Math.\eol
La Jolla, CA  92093}

\Eblock{Jesse L. Barlow\eol
Pennsylvania State Univ.\eol
Dept. of Comp. Sci.\eol
Univ. Park, PA 16802}

\Eblock{Chris Bischof\eol
Argonne National Lab.\eol
Math. and Comp. Sci. Div.\eol
9700 S Cass Av.\eol
Argonne, IL 60439}

\Eblock{R. H. Bisseling\eol
Koninklijke/Shell\eol
Laboratorium Amsterdam\eol
P. O. Box 3003\eol
1003 AA Amsterdam\eol
The Netherlands}

\Eblock{Ake Bjorck\eol
Linkoping Univ.\eol
Dept. of Math.\eol
S-581 83 Linkoping\eol
Sweden}

\Eblock{Petter E. Bjorstad\eol
Univ. of Bergen\eol
Thormonhlensgt. 55\eol
N-5006 Bergen\eol
Norway}

\Eblock{Daniel Boley\eol
Stanford Univ.\eol
Dept. of Comp. Sci.\eol
Stanford, CA 94305}

\Eblock{Randall Bramley\eol
Univ. of Illinois\eol
305 Talbot Lab.\eol
104 S Wright Street\eol
Urbana, IL 61801}

\Eblock{Richard A. Brualdi\eol
Univ. of Wisconsin\eol
Dept. of Math.\eol
480 Lincoln Dr.\eol
Madison, WI 53706}

\Eblock{James R. Bunch\eol
Univ. of California, San Diego\eol
Dept. of Math.\eol
La Jolla, CA  92093}

\Eblock{Angelica Bunse-Gerstner\eol
Univ. Bielefeld\eol
Fakultat fur Math.\eol
postfach 8640\eol
D-4800 Bielfeld 1\eol
West Germany}

\Eblock{Ralph Byers\eol
Univ. of Kansas\eol
Dept. of Math.\eol
Lawrence, KA 66045-2142}

\Eblock{Jose Castillo\eol
San Diego State Univ.\eol
Dept. of Math. Sciences\eol
San Diego, CA  92182}

\Eblock{Tony F. Chan\eol
Univ. of California, Los Angles\eol
Dept. of Math.\eol
Los Angeles, CA 90024}

\Eblock{S. S. Chow\eol
Univ. of Wyoming\eol
Dept. of Math.\eol
Laramie, WY 82071}

\Eblock{Anthony T. Chronopoulos\eol
Univ. of Minnesota\eol
Dept. of Comp. Sci.\eol
200 Union Street SE\eol
Minneapolis, MN 55455}

\Eblock{Eleanor C. H. Chu\eol
Univ. of Waterloo\eol
Dept. of Comp. Sci.\eol
Waterloo N2L 3G1\eol
Ontario, Canada}

\Eblock{Len Colgan\eol
South Australian Inst. of Tech.\eol
Math. Dept.\eol
The Levels, 5095, South Australia\eol
Australia}

\Eblock{Paul Concus\eol
Univ. of California, Berkeley\eol
Lawrence Berkeley Lab.\eol
50A-2129\eol
Berkeley, CA 94720}

\Eblock{W. M. Coughran\eol
AT\&T Bell Laboratories\eol
600 Mountain Av.\eol
Murray Hill, NJ 07974-2070}

\Eblock{Jane K. Cullum\eol
IBM T. J. Watson Res. Center\eol
P. O. Box 218\eol
Yorktown Heights, NY 10598}

\Eblock{Carl De Boor\eol
Univ. of Wisconsin\eol
Center for Math. Sciences\eol
610 Walnut Street\eol
Madison, WI 53706}

\Eblock{John E. De Pillis\eol
Univ. of California, Riverside\eol
Dept. of Math. and Comp. Sci.\eol
Riverside, CA 92521}

\Eblock{James W. Demmel\eol
New York Univ.\eol
Courant Institute\eol
251 Mercer Street\eol
New York, NY 10012}

\Eblock{Julio C. Diaz\eol
Univ. of Tulsa\eol
Dept. of Math. and Comp. Sci.\eol
600 S College Av.\eol
Tulsa, OK 74104-3189}

\Eblock{David S. Dodson\eol
Convex Computer Corp.\eol
701 N Plano Rood\eol
Richardson, TX 75081}

\Eblock{Jack Dongarra\eol
Univ. of Tennessee\eol
Comp. Sci. Dept.\eol
Knoxville, TN 37996-1300}

\Eblock{Iain S. Duff\eol
Harwell Lab.\eol
Comp. Sci. and Systems Div.\eol
Oxfordshire OX11 0RA\eol
United Kingdom}

\Eblock{Pat Eberlein\eol
State Univ. of New York\eol
Dept. of Comp. Sci.\eol
Buffalo, NY 14260}

\Eblock{W. Stuart Edwards\eol
Univ. of Texas\eol
Center for Nonlinear Dynamics\eol
Austin, TX 78712}

\Eblock{Louis W. Ehrlich\eol
Johns Hopkins Univ.\eol
Applied Physics Lab.\eol
Johns Hopkins Road\eol
Laurel, MD 20707}

\Eblock{Michael Eiermann\eol
Univ. Karlsruhe\eol
Inst. fur Prakticshe Math.\eol
Postfach 6980\eol
D-6980 Karlsruhe 1\eol
West Germany}

\Eblock{Stanley C. Eisenstat\eol
Yale Univ.\eol
Dept. of Comp. Sci.\eol
P. O. Box 2158, Yale Station\eol
New Haven, CT 06520}

\Eblock{Lars Elden\eol
Linkoping Univ.\eol
Dept. of Math.\eol
S-581 83 Linkoping\eol
Sweden}

\Eblock{Howard C. Elman\eol
Univ. of Maryland\eol
Dept. of Comp. Sci.\eol
College Park, MD 20742}

\Eblock{T. Ericsson\eol
Univ. of Umea\eol
Inst. of Information Processing\eol
S-901 87 Umea\eol
Sweden}

\Eblock{Albert M. Erisman\eol
Boeing Computer Services\eol
565 Andover Park West\eol
Mail Stop 9C-01\eol
Tukwila, WA 98188}

\Eblock{D. J. Evans\eol
Univ. of Technology\eol
Dept. of Computer Studies\eol
Leicestershire LE11 3TU\eol
United Kingdom}

\Eblock{Vance Faber\eol
Los Alamos National Lab.\eol
Group C-3\eol
Mail Stop B265\eol
Los Alamos, NM 87545}

\Eblock{Bernd Fischer\eol
Univ. Hamburg\eol
Inst. fur Angew. Math.\eol
D-2000 Hamburg 13\eol
West Germany}

\Eblock{Geoffrey Fox\eol
California Inst. of Tech.\eol
Mail Code 158-79\eol
Pasadena, CA 91125}

\Eblock{Paul O. Frederickson\eol
NASA Ames Res. Center\eol
RIACS, Mail Stop 230-5\eol
Moffett Field, CA 94035}

\Eblock{Roland W. Freund\eol
NASA Ames Res. Center\eol
RIACS, Mail Stop 230-5\eol
Moffett Field, CA 94035}

\Eblock{Robert E. Funderlic\eol
North Carolina State Univ.\eol
Dept. of Comp. Sci.\eol
Raleigh, NC 27650}

\Eblock{Ralf Gaertner\eol
Max Plank Society\eol
Fritz Haber Inst.\eol
Faradayweg 4-6\eol
D-1000 Berlin 33\eol
West Germany}

\Eblock{Patrick W. Gaffney\eol
Bergen Scientific Center\eol
Allegaten 36\eol
N-5000 Bergen\eol
Norway}

\Eblock{Kyle A. Gallivan\eol
Univ. of Illinois\eol
305 Talbot Lab.\eol
104 S Wright Street\eol
Urbana, IL 61801}

\Eblock{Dennis B. Gannon\eol
Indiana Univ.\eol
Dept. of Comp. Sci.\eol
Bloomington, IN 47405-6171}

\Eblock{Kevin E. Gates\eol
Univ. of Washington\eol
Dept. of Applied Math.\eol
Seattle, WA 98195}

\Eblock{David M. Gay\eol
AT\&T Bell Laboratories\eol
600 Mountain Av.\eol
Murray Hill, NJ 07974-2070}

\Eblock{C. William Gear\eol
Univ. of Illinois\eol
Dept. of Comp. Sci.\eol
1304 W Springfield Av.\eol
Urbana, IL 61801}

\Eblock{George A. Geist\eol
Oak Ridge National Lab.\eol
Math. Sciences Section\eol
Building 9207-A\eol
P. O. Box 2009\eol
Oak Ridge, TN 37831}

\Eblock{J. Alan George\eol
Univ. of Waterloo\eol
Academic Vice President and Provost\eol
Needles Hall\eol
Waterloo N2L 3G1\eol
Ontario, Canada}

\Eblock{Adam Gersztenkorn\eol
Amoco Production Company\eol
Geophysical Res. Dept.\eol
4502 E 41st Street\eol
P. O. Box 3385\eol
Tulsa, OK 74102}

\Eblock{John R. Gilbert\eol
Xerox Palo Alto Research Center\eol
3333 Coyote Hill Road\eol
Palo Alto, CA 94304}

\Eblock{Albert Gilg\eol
Siemens AG\eol
Corporate Res. and Tech.\eol
Otto Hahn Ring 6\eol
D-8000 Munchen 83\eol
West Germany}

\Eblock{Roland Glowinski\eol
Univ. of Houston\eol
Dept. of Math.\eol
4800 Calhoun Road\eol
Houston, TX 77004}

\Eblock{Gene H. Golub\eol
Stanford Univ.\eol
Dept. of Comp. Sci.\eol
Stanford, CA 94305}

\Eblock{W. Gragg\eol
Naval Postgraduate School\eol
Dept. of Math.\eol
Mail Code 53ZH\eol
Monterey, CA 93943-5100}

\Eblock{Anne Greenbaum\eol
New York Univ.\eol
Courant Institute\eol
251 Mercer Street\eol
New York, NY 10012}

\Eblock{Roger Grimes\eol
Boeing Computer Services\eol
Mail Stop 7L-21\eol
P. O. Box 24346\eol
Seattle, WA 98124-0346}

\Eblock{W. D. Gropp\eol
Yale Univ.\eol
Dept. of Comp. Sci.\eol
P. O. Box 2158, Yale Station\eol
New Haven, CT 06520}

\Eblock{Fred G. Gustavson\eol
IBM T. J. Watson Res. Center\eol
P. O. Box 218\eol
Yorktown Heights, NY 10598}

\Eblock{Martin H. Gutknecht\eol
ETH-Zentrum\eol
IPS, IFW D 25.1\eol
CH-8092 Zurich\eol
Switzerland}

\Eblock{Louis Hageman\eol
Bettis Atomic Power Lab.\eol
Box 79\eol
West Mifflin, PA 15122}

\Eblock{Charles A. Hall\eol
Univ. of Pittsburgh\eol
Dept. of Math. and Statistics\eol
Pittsburgh, PA 15260}

\Eblock{Sven J. Hammarling\eol
Numerical Algorithms Group Ltd.\eol
Wilkinson House\eol
Jordan Hill Road\eol
Oxford OX2 8DR\eol
United Kingdom}

\Eblock{Michael T. Heath\eol
Oak Ridge National Lab.\eol
Math. Sciences Section\eol
Building 9207-A\eol
P. O. Box 2009\eol
Oak Ridge, TN 37831-8083}

\Eblock{Don E. Heller\eol
Shell Development Company\eol
Bellaire Res. Center\eol
P. O. Box 481\eol
Houston, TX 77001}

\Eblock{Nicholas J. Higham\eol
Cornell Univ.\eol
Dept. of Comp. Sci.\eol
Ithaca, NY 14853}

\Eblock{Mary Hill\eol
USGS-WRD\eol
Denver Federal Center\eol
Mail Stop 413\eol
P. O. Box 25046\eol
Lakewood, CO 80225}

\Eblock{Mike Holst\eol
Univ. of Illinois\eol
Dept. of Comp. Sci.\eol
1304 W Springfield Av.\eol
Urbana, IL 61801}

\Eblock{James M. Hyman\eol
Los Alamos National Lab.\eol
Group T-7\eol
Mail Stop B284\eol
Los Alamos, NM 87545}

\Eblock{Ilse Ipsen\eol
Yale Univ.\eol
Dept. of Comp. Sci.\eol
P. O. Box 2158, Yale Station\eol
New Haven, CT 06520}


\Eblock{S. Lennart Johnsson\eol
Yale Univ.\eol
Dept. of Comp. Sci.\eol
P. O. Box 2158, Yale Station\eol
New Haven, CT 06520}

\Eblock{J. Jerome\eol
Northwestern Univ.\eol
2033 Sheridan Road\eol
Evanston, IL 60208}

\Eblock{Tom Jordan\eol
Los Alamos National Lab.\eol
Group C-3\eol
Mail Stop B265\eol
Los Alamos, NM 87545}

\Eblock{Wayne Joubert\eol
Univ. of Texas\eol
Center for Numerical Anlysis\eol
Moore Hall 13.150\eol
Austin, TX 78712}

\Eblock{E. F. Kaasschieter\eol
Inst. of Applied Geoscience DGV-TNO\eol
P. O. Box 285\eol
Schoemakerstraat 97\eol
2600 AG Delft\eol
The Netherlands}

\Eblock{W. M. Kahan\eol
Univ. of California, Berkeley\eol
Dept. of Math.\eol
Berkeley, CA 94720}

\Eblock{Bo Kagstrom\eol
Univ. of Umea\eol
Inst. of Information Processing\eol
S-901 87 Umea\eol
Sweden}

\Eblock{Chandrika Kamath\eol
Digital Equipment Corp.\eol
HLO2-3/M08\eol
77 Reed Road\eol
Hudson, MA 01749}

\Eblock{Hans Kaper\eol
Argonne National Lab.\eol
Math. and Comp. Sci. Div.\eol
9700 S Cass Av.\eol
Argonne, IL 60439}

\Eblock{Linda Kaufman\eol
AT\&T Bell Laboratories\eol
600 Mountain Av.\eol
Murray Hill, NJ 07974-2070}

\Eblock{Mark D. Kent\eol
Stanford Univ.\eol
Dept. of Comp. Sci.\eol
Stanford, CA 94305}

\Eblock{T. Kerkoven\eol
Univ. of Illinois\eol
Dept. of Comp. Sci.\eol
1304 W Springfield Av.\eol
Urbana, IL 61801}

\Eblock{David Keyes\eol
Yale Univ.\eol
Dept. of Mechanical Eng.\eol
P. O. Box 2159 Yale Station\eol
New Haven, CT 06520}

\Eblock{David R. Kinkaid\eol
Univ. of Texas\eol
Center for Numerical Anlysis\eol
Moore Hall 13.150\eol
Austin, TX 78712}

\Eblock{Virgina Klema\eol
Massachusetts Inst. of Tech.\eol
Statistics Center\eol
Cambridge, MA 02139}

\Eblock{Steven G. Kratzer\eol
Supercomputing Res. Center\eol
17100 Science Dr.\eol
Bowie, MD 20715-4300}

\Eblock{Edward J. Kushner\eol
Intel Scientific Computers\eol
15201 NW Greenbrier Pkwy.\eol
Beaverton, OR 97006}

\Eblock{John Lavery\eol
NASA Lewis Res. Center\eol
Mail Stop 5-11\eol
Cleveland, OH 44135}

\Eblock{Kincho H. Law\eol
Stanford Univ.\eol
Dept. of Civil Eng.\eol
Stanford, CA 94305-4020}

\Eblock{Charles Lawson\eol
Jet Propulsion Lab.\eol
Applied Math. Group\eol
Mail Stop 506-232\eol
4800 Oak Grove Dr.\eol
Pasadena, CA 91109}

\Eblock{Yannick Le Coz\eol
Rensselaer Polytechnic Inst.\eol
Electrical Eng. Dept.\eol
Troy, NY 12180}

\Eblock{Steve L. Lee\eol
Univ. of Illinois\eol
Dept. of Comp. Sci.\eol
1304 W Springfield Av.\eol
Urbana, IL 61801}

\Eblock{Steve Leon\eol
Southeastern Massachusetts Univ.\eol
Dept. of Math.\eol
North Dartmouth, MA 02747}

\Eblock{Michael R. Leuze\eol
Oak Ridge National Lab.\eol
Math. Sciences Section\eol
P. O. Box 2009\eol
Oak Ridge, TN 37831}

\Eblock{John G. Lewis\eol
Boeing Computer Services\eol
Mail Stop 7L-21\eol
P. O. Box 24346\eol
Seattle, WA 98124-0346}

\Eblock{Antonios Liakopoulos\eol
System Dynamics Inc.\eol
1211 NW 10th Av.\eol
Gainesville, FL 32601}

\Eblock{Joseph W. H. Liu\eol
York Univ.\eol
Dept. of Comp. Sci.\eol
North York M3J 1P3\eol
Ontario, Canada}

\Eblock{Peter Lory\eol
Tech. Univ. Munchen\eol
Dept. of Math.\eol
POB 20 24 20\eol
D-8000 Munchen 2\eol
West Germany}

\Eblock{Robert F. Lucas\eol
Supercomputing Res. Center\eol
17100 Science Dr.\eol
Bowie, MD 20715-4300}

\Eblock{Franklin Luk\eol
Cornell Univ.\eol
School of Electrical Eng.\eol
Ithaca, NY 14853}

\Eblock{Tom Manteuffel\eol
Univ. of Colorado\eol
Computational Math. Group\eol
1100 14th Street\eol
Denver, CO 80202}

\Eblock{Robert Mattheij\eol
Eindhoven Univ. of Tech.\eol
Dept. of Math.\eol
5600 MB Eindhoven\eol
The Netherlands}

\Eblock{Paul C. Messina\eol
California Inst. of Tech.\eol
Mail Code 158-79\eol
Pasadena, CA 91125}

\Eblock{Gerard A. Meurant\eol
Centre d'Etudes de Limeil\eol
Service Mathematiques Appliquees\eol
Boite Postale 27\eol
94190 Villeneuve St. Georges\eol
France}

\Eblock{Carl D. Meyer\eol
North Carolina State Univ.\eol
Dept. of Math.\eol
Raleigh, NC 27650}

\Eblock{Ignacy Misztal\eol
Univ. of Illinois\eol
Dept. of Animal Sciences\eol
1207 W Gregory Dr.\eol
Urbana, IL 61801}

\Eblock{Gautam Mitra\eol
Brunel Univ.\eol
Dept. of Math. and Statistics\eol
Uxbridge Middlesex UB8 3PH\eol
United Kingdom}

\Eblock{Hans Mittleman\eol
Arizona State Univ.\eol
Dept. of Math.\eol
Tempe, AZ 85287}

\Eblock{Cleve Moler\eol
Ardent Computers\eol
550 Del Ray Av.\eol
Sunnyvale, CA 94086}

\Eblock{Rael Morris\eol
IBM Corp.\eol
Almaden Res. Center\eol
Dept. K08/282\eol
650 Garry Road\eol
San Jose, CA 95120-6099}

\Eblock{Edmond Nadler\eol
Wayne State Univ.\eol
Dept. of Math.\eol
Detroit, MI 48202}

\Eblock{N. Nandakumar\eol
Univ. of Nebraska\eol
Dept. of Math. and Comp. Sci.\eol
Omaha, NB 68182}

\Eblock{Olavi Nevanlinna\eol
Helsinki Univ. of Tech.\eol
Inst. of Math.\eol
SF-02150 Espoo\eol
Finland}

\Eblock{Esmond G. Y. Ng\eol
Oak Ridge National Lab.\eol
Math. Sciences Section\eol
Building 9207-A\eol
P. O. Box 2009\eol
Oak Ridge, TN 37831}

\Eblock{Viet-nam Nguyen\eol
Pratt \& Whitney\eol
Dept. of Eng. and Comp. Applications\eol
1000 Marie Victorian\eol
Longueuil J4G 1A1\eol
Quebec, Canada}

\Eblock{Nancy Nichols\eol
Reading Univ.\eol
Dept. of Math.\eol
Whiteknights Park\eol
Reading RG6 2AX\eol
United Kingdom}

\Eblock{W. Niethammer\eol
Univ. Karlsruhe\eol
Inst. fur Praktische Math.\eol
Engelerstrasse 2\eol
D-7500 Karlsruhe\eol
West Germany}

\Eblock{Takashi Nodera\eol
Keio Univ.\eol
Dept. of Math.\eol
3-14-1 Hiyoshi Kohoku\eol
Yokohama 223\eol
Japan}

\Eblock{Wilbert Noronha\eol
Univ. of Tennessee\eol
310 Perkins Hall\eol
Knoxville, TN 37996}

\Eblock{Bahram Nour-Omid\eol
Lockheed Palo Alto Res. Lab.\eol
Comp. Mechanics Section\eol
Org. 93-30 Bldg. 251\eol
3251 Hanover Street\eol
Palo Alto, CA 94304}

\Eblock{Dianne P. O'Leary\eol
Univ. of Maryland\eol
Dept. of Comp. Sci.\eol
College Park, MD 20742}

\Eblock{Julia Olkin\eol
SRI International\eol
Building 301, Room 66\eol
333 Ravenswood Av.\eol
Menlo Park, CA 94025}

\Eblock{Steve Olson\eol
Supercomputer Systems Inc\eol
1414 W Hamilton Av.\eol
Eau Claire, WI 54701}

\Eblock{Elizabeth Ong\eol
Univ. of California, Los Angles\eol
Dept. of Math.\eol
Los Angeles, CA 90024}

\Eblock{James M. Ortega\eol
Univ. of Virginia\eol
Dept. of Applied Math.\eol
Charlottesville, VA 22903}

\Eblock{Christopher C. Paige\eol
McGill Univ.\eol
School of Comp. Sci.\eol
3480 University\eol
Montreal H3A 2A7\eol
Quebec, Canada}

\Eblock{M. C. Pandian\eol
IBM Corp.\eol
Numerically Intensive Computing\eol
Dept. 41U/276\eol
Neighborhood Road\eol
Kingston, NY 12401}

\Eblock{Roy Pargas\eol
Clemson Univ.\eol
Dept. of Comp. Sci.\eol
Clemson, SC 29634-1906}

\Eblock{Beresford N. Parlett\eol
Univ. of California, Berkeley\eol
Dept. of Math.\eol
Berkeley, CA 94720}

\Eblock{Merrell Patrick\eol
Duke Univ.\eol
Dept. of Comp. Sci.\eol
Durham, NC 27706}

\Eblock{Victor Pereyra\eol
Weidlinger Associates\eol
Suite 110\eol
4410 El Camino Real\eol
Los Altos, CA 94022}

\Eblock{Barry W. Peyton\eol
Oak Ridge National Lab.\eol
Math. Sciences Section\eol
Building 9207-A\eol
P. O. Box 2009\eol
Oak Ridge, TN 37831}

\Eblock{Robert J. Plemmons\eol
North Carolina State Univ.\eol
Dept. of Math.\eol
Raleigh, NC 27650}

\Eblock{Jim Purtilo\eol
Univ. of Maryland\eol
Dept. of Comp. Sci.\eol
College Park, MD 20742}

\Eblock{Giuseppe Radicati\eol
IBM Italia\eol
Via Giorgione 159\eol
00147 Roma\eol
Italy}

\Eblock{A. Ramage\eol
Univ. of Bristol\eol
Dept. of Math.\eol
University Walk\eol
Bristol BS8 1TW\eol
United Kingdom}

\Eblock{Lothar Reichel\eol
Bergen Scientific Centre\eol
Allegaten 36\eol
N-5007 Bergen\eol
Norway}

\Eblock{John K. Reid\eol
Harwell Lab.\eol
Comp. Sci. and Systems Div.\eol
Oxfordshire OX11 0RA\eol
United Kingdom}

\Eblock{John R. Rice\eol
Purdue Univ.\eol
Dept. of Comp. Sci.\eol
Layfayette, IN 47907}

\Eblock{Jeff V. Richard\eol
Science Applications International\eol
Mail Stop 34\eol
10260 Campus Point Dr.\eol
San Diego, CA 92121}

\Eblock{Donald J. Rose\eol
Duke Univ.\eol
Dept. of Comp. Sci.\eol
Durham, NC 27706}

\Eblock{Vona Bi Roubolo\eol
Univ. of Texas\eol
Center for Numerical Anlysis\eol
Moore Hall 13.150\eol
Austin, TX 78712}

\Eblock{Axel Ruhe\eol
Chalmers Tekniska Hogskola\eol
Dept. of Comp. Sci.\eol
S-412 96 Goteborg\eol
Sweden}

\Eblock{Youcef Saad\eol
NASA Ames Res. Center\eol
RIACS, Mail Stop 230-5\eol
Moffett Field, CA 94035}

\Eblock{P. Sadayappan\eol
Ohio State Univ.\eol
Dept. of Comp. Sci.\eol
Columbus, OH 43210}

\Eblock{Joel Saltz\eol
Yale Univ.\eol
Dept. of Comp. Sci.\eol
P. O. Box 2158, Yale Station\eol
New Haven, CT 06520}

\Eblock{Ahmed H. Sameh\eol
Univ. of Illinois\eol
305 Talbot Lab.\eol
104 S Wright Street\eol
Urbana, IL 61801}

\Eblock{Michael Saunders\eol
Stanford Univ.\eol
Dept. of Operations Research\eol
Stanford, CA 94305}

\Eblock{Paul E. Saylor\eol
Univ. of Illinois\eol
Dept. of Comp. Sci.\eol
1304 W Springfield Av.\eol
Urbana, IL 61801}

\Eblock{Mark Schaefer\eol
Texas A \& M Univ.\eol
Dept. of Math.\eol
College Station, TX 77843}

\Eblock{U. Schendel\eol
Freie Univ. Berlin\eol
Inst. fur Math.\eol
Arnimallee 2-6\eol
D-1000 Berlin 33\eol
West Germany}

\Eblock{Robert S. Schreiber\eol
NASA Ames Res. Center\eol
RIACS, Mail Stop 230-5\eol
Moffett Field, CA 94035}

\Eblock{Martin Schultz\eol
Yale Univ.\eol
Dept. of Comp. Sci.\eol
P. O. Box 2158, Yale Station\eol
New Haven, CT 06520}

\Eblock{David St. Clair Scott\eol
Intel Scientific Computers\eol
15201 NW Greenbrier Pkwy.\eol
Beaverton, OR 97006}

\Eblock{Jeffrey S. Scroggs\eol
NASA Langley Res. Center\eol
ICASE, Mail Stop 132-C\eol
Hampton, VA 23665}

\Eblock{Steven M. Serbin\eol
Univ. of Tennessee\eol
Dept. of Math.\eol
Knoxville, TN 37996-1300}

\Eblock{Shahriar Shamsian\eol
MacNeal Schwendler Corp.\eol
175 S Madison\eol
Pasadena, CA 91101}

\Eblock{A. H. Sherman\eol
Scientific Computing Associates Inc.\eol
Suite 307\eol
246 Church Street\eol
New Haven, CT 06510}

\Eblock{Kermit Sigmon\eol
Univ. of Florida\eol
Dept. of Math.\eol
Gainesville, FL 32611}

\Eblock{Horst D. Simon\eol
NASA Ames Res. Center\eol
Mail Stop 258-5\eol
Moffett Field, CA 94035}

\Eblock{Richard F. Sincovec\eol
NASA Ames Res. Center\eol
RIACS, Mail Stop 230-5\eol
Moffett Field, CA 94035}

\Eblock{Dennis Smolarski\eol
Univ. of Santa Clara\eol
Dept. of Math.\eol
Santa Clara, CA 95053}

\Eblock{Mitchell D. Smooke\eol
Yale Univ.\eol
Dept. of Mechanical Eng.\eol
P. O. Box 2159 Yale Station\eol
New Haven, CT 06520}

\Eblock{P. Sonneveld\eol
Delft Univ. of Tech.\eol
Dept. of Math. and Informatics\eol
P. O. Box 356\eol
2600 AJ Delft\eol
The Netherlands}

\Eblock{Danny C. Sorensen\eol
Rice Univ.\eol
Dept. of Math. Sciences\eol
Houston, TX 77251-1892}

\Eblock{Alistair Spence\eol
Univ. of Bath\eol
School of Math. Sciences\eol
Bath BA2 7AY\eol
United Kingdom}

\Eblock{Trond Steihaug\eol
Statoil, Forus\eol
P. O. Box 300\eol
N-4001 Stravanger\eol
Norway}

\Eblock{G. W. Stewart\eol
Univ. of Maryland\eol
Dept. of Comp. Sci.\eol
College Park, MD 20742}

\Eblock{G. Strang\eol
Massachusetts Inst. of Tech.\eol
Dept. of Math.\eol
Cambridge, MA 02139}

\Eblock{Jin Su\eol
IBM Corp.\eol
Dept. 41UC/276\eol
Neighborhood Road\eol
Kingston, NY 12498}

\Eblock{Uwe Suhl\eol
IBM Corp.\eol
Math. Sciences Dept.\eol
Bergen Gruen Strasse 17-19\eol
D-1000 Berlin 38\eol
West Germany}

\Eblock{Daniel B. Szyld\eol
Duke Univ.\eol
Dept. of Comp. Sci.\eol
Durham, NC 27706}

\Eblock{Hillel Tal-Ezer\eol
Brown Univ.\eol
Div. of Applied Math.\eol
Box F\eol
Providence, RI 02912}

\Eblock{Wei Pai Tang\eol
Univ. of Waterloo\eol
Dept. of Comp. Sci.\eol
Waterloo N2L 3G1\eol
Ontario, Canada}

\Eblock{G. D. Taylor\eol
Colorado State Univ.\eol
Math. Dept.\eol
Fort Collins, CO 80523}

\Eblock{Charles Tong\eol
Univ. of California, Los Angeles\eol
Dept. of Comp. Sci.\eol
Los Angeles, CA 90024}

\Eblock{Toru Toyabe\eol
Hitachi Ltd.\eol
Central Research Lab.\eol
7th Dept.\eol
Kokobunji, Tokyo 185\eol
Japan}

\Eblock{Eugene A. Trabka\eol
Eastman Kodak Company\eol
Research Laboratories\eol
Rochester, NY 14650}

\Eblock{L. N. Trefethen\eol
Massachusetts Inst. of Tech.\eol
Dept. of Math.\eol
Cambridge, MA 02139}

\Eblock{Donato Trigiante\eol
Univ. di Bari\eol
Inst. di Matematica\eol
Campus Universitario\eol
70125 Bari\eol
Italy}

\Eblock{Donald G. Truhlar\eol
Univ. of Minnesota\eol
Minnesota Supercomputer Inst.\eol
1200 Washington Av. South\eol
Minneapolis, MN 55415}

\Eblock{Kathryn L. Turner\eol
Utah State Univ.\eol
Dept. of Math.\eol
Logan, UT 84322-3900}

\Eblock{Charles Van Loan\eol
Cornell Univ.\eol
Dept. of Comp. Sci.\eol
Ithaca, NY 14853}

\Eblock{Henk Van der Vorst\eol
Delft Univ. of Tech.\eol
Dept. of Math. and Informatics\eol
P. O. Box 356\eol
2600 AJ Delft\eol
The Netherlands}

\Eblock{James M. Varah\eol
Univ. of British Columbia\eol
Dept. of Comp. Sci.\eol
Vancouver, V6T 1W5\eol
British Columbia, Canada}

\Eblock{Richard S. Varga\eol
Kent State Univ.\eol
Dept. of Math.\eol
Kent, OH 44242}

\Eblock{Anthony Vassiliou\eol
Mobil Res. and Development Corp.\eol
Dallas Reseach Lab.\eol
13777 Midway Road\eol
Dallas, TX 75244}

\Eblock{Robert G. Voigt\eol
NASA Langley Res. Center\eol
ICASE, Mail Stop 132-C\eol
Hampton, VA 23665}

\Eblock{Eugene L. Wachspress\eol
Univ. of Tennessee\eol
Dept. of Math.\eol
Knoxville, TN 37996-1300}

\Eblock{Robert C. Ward\eol
Oak Ridge National Lab.\eol
Math. Sciences Section\eol
Building 9207-A\eol
P. O. Box 2009\eol
Oak Ridge, TN 37831}

\Eblock{Daniel D. Warner\eol
Clemson Univ.\eol
Dept. of Math. Sciences\eol
Clemson, SC 29634-1907}

\Eblock{Andrew J. Wathen\eol
Univ. of Bristol\eol
Dept. of Math.\eol
University Walk\eol
Bristol BS8 1TW\eol
United Kingdom}

\Eblock{Bruno Welfert\eol
Univ. of California, San Diego\eol
Dept. of Math.\eol
La Jolla, CA  92093}

\Eblock{Mary J. Wheeler\eol
Univ. of Houston\eol
Dept. of Math.\eol
Houston, TX 77204-3476}

\Eblock{Andrew B. White\eol
Los Alamos National Lab.\eol
Group C-3\eol
Mail Stop B265\eol
Los Alamos, NM 87545}

\Eblock{Jacob White\eol
Massachusetts Inst. of Tech.\eol
Dept. of Electrical Eng.\eol
Cambridge, MA 02139}

\Eblock{Olaf Widlund\eol
New York Univ.\eol
Courant Institute\eol
251 Mercer Street\eol
New York, NY 10012}

\Eblock{Harry A. Wijshoff\eol
Univ. of Illinois\eol
305 Talbot Lab.\eol
104 S Wright Street\eol
Urbana, IL 61801}

\Eblock{Roy S. Wikramaratna\eol
Winfrith Petroleum Tech.\eol
Dorchester\eol
Dorset DT2 8DH\eol
United Kingdom}

\Eblock{David S. Wise\eol
Indiana Univ.\eol
Dept. of Comp. Sci.\eol
Bloomington, IN 47405-6171}

\Eblock{P. H. Worley\eol
Oak Ridge National Lab.\eol
Math. Sciences Section\eol
Building 9207-A\eol
P. O. Box 2009\eol
Oak Ridge, TN 37831}

\Eblock{Margaret H. Wright\eol
AT\&T Bell Laboratories\eol
600 Mountain Av.\eol
Murray Hill, NJ 07974-2070}

\Eblock{Steven J. Wright\eol
Argonne National Lab.\eol
Math. and Comp. Sci. Div.\eol
9700 S Cass Av.\eol
Argonne, IL 60439}

\Eblock{Kuo W. Wu\eol
Cray Research Inc\eol
1333 Northland Dr.\eol
Mendota Heights, MN 55120-1095}

\Eblock{Chao W. Yang\eol
Cray Research Inc\eol
Math. Software Group\eol
1408 Northland Dr.\eol
Mendota Heights, MN 55120-1095}


\Eblock{Gung-Chung Yang\eol
Univ. of Illinois\eol
305 Talbot Lab.\eol
104 S Wright Street\eol
Urbana, IL 61801}

\Eblock{Elizabeth Yip\eol
Boeing Aerospace Corp.\eol
Mail Stop 8K-17\eol
P. O. Box 3999\eol
Seattle, WA 98124-2499}

\Eblock{David P. Young\eol
Boeing Computer Services\eol
Mail Stop 7L-21\eol
P. O. Box 24346\eol
Seattle, WA 98124-0346}

\Eblock{David M. Young\eol
Univ. of Texas\eol
Center for Numerical Anlysis\eol
Moore Hall 13.150\eol
Austin, TX 78712}

\Eblock{Earl Zmijewski\eol
Univ. of California, Santa Barbara\eol
Comp. Sci. Dept.\eol
Santa Barbara, CA 93106}

\Eblock{Qisu Zou\eol
Kansas State Univ.\eol
Dept. of Math.\eol
Manhattan, KS 66506}

 \Lblock
{S. F. Ashby, LLNL, L-316\eol
 J. B. Bell, LLNL, L-316\eol
 J. H. Bolstad, LLNL, L-16\eol
 P. N. Brown, LLNL, L-316\eol
 R. C. Chin, LLNL, L-321\eol
 M. R. Dorr, LLNL, L-316\eol
 G. W. Hedstrom, LLNL, L-321\eol
 A. C. Hindmarsh, LLNL, L-316\eol
 D. S. Kershaw, LLNL, L-471\eol
 A. E. Koniges, LLNL, L-561\eol
 L. R. Petzold, LLNL, L-316\eol
 G. H. Rodrigue, LLNL, L-306\eol
 J. A. Trangenstein, LLNL, L-316}

 \vfil \break
{\Sblock \settabs \+\hglue 0.5in& \hglue0.25in& \hglue 0.25in\cr 
 \+1420&W. J. Camp\cr\filbreak
 \+1421&R. J. Thompson\cr\filbreak
 \+1422&R. C. Allen, Jr.\cr\filbreak
 \+1423&E. F. Brickell\cr\filbreak
 \+1424&G. G. Weigand\cr\filbreak
 \+1553&W. L. Hermina\cr\filbreak
 \+2113&H. A. Watts\cr\filbreak
 \+8000& J. C. Crawford\cr
 \+&&Attn: &E. E. Ives, 8100\cr
 \+&&&P. L. Mattern, 8300\cr
 \+&&&R. C. Wayne, 8400\cr
 \+&&&P. E. Brewer, 8500\cr
 \+8200& R. J. Detry\cr
 \+&&Attn: &C. W. Robinson, 8240\cr
 \+&&&R. C. Dougherty, 8270\cr
 \+&&&R. A. Baroody, 8280\cr
 \+8230& W. D. Wilson\cr
 \+&&Attn: &M. H. Pendley, 8234\cr
 \+&&&D. L. Crawford, 8235\cr
 \+&&&P. W. Dean, 8236\cr\filbreak
 \+8231&W. E. Mason\cr\filbreak
 \+8233&J. M. Harris\cr\filbreak
 \+8233&R. E. Cline\cr\filbreak
 \+8233&J. F. Grcar (40)\cr\filbreak
 \+8233&J. F. Lathrop\cr\filbreak
 \+8233&J. C. Meza\cr\filbreak
 \+8233&R. Y. Lee\cr\filbreak
 \+8233&R. A. Whiteside\cr\filbreak
 \+8241&K. J. Perano\cr\filbreak
 \+8243&L. A. Bertram\cr\filbreak
 \+8245&R. J. Kee\cr\filbreak
 \+8245&G. H. Evans\cr\filbreak
 \+8245&W. G. Houf\cr\filbreak
 \+8245&A. E. Lutz\cr\filbreak
 \+8245&W. S. Winters\cr\filbreak}

 {\Sblock \settabs \+\hglue 0.5in& \hglue0.25in& \hglue 0.25in\cr 
 \+8535&Publications/Tech. Lib. Processes, 3141\cr\filbreak
 \+3141&Technical Library Process Division (3)\cr\filbreak
 \+8524-2&Central Technical Files (3)\cr\filbreak}

\end